\documentclass{amsart}
\usepackage{amscd}
\usepackage{amssymb}
\usepackage{mathrsfs} 

\input xy
\xyoption{all}

\def\ip#1{\left< #1 \right>}
\def\ad#1{{\scriptstyle \left[ \frac{1}{#1} \right]}}

\def\smsh{\wedge}
\def\eq{{\text{eq}}}
\def\nil{\operatorname{nil}}
\def\chr{\operatorname{char}}
\def\ord{\operatorname{ord}}

\def\ker{\operatorname{ker}}

\def\cMpctf{{\mathcal M}_{\text{pctf}}}

\def\cJ{\mathcal J}

\def\cA{\mathcal A}
\def\cB{\mathcal B}
\def\cC{\mathcal C}

\def\cD{\mathcal D}
\def\cF{\mathcal F}
\def\cG{\mathcal G}
\def\cO{\mathcal O}
\def\cS{\mathcal S}

\def\cI{\mathcal I}

\def\cdh{{\text{cdh}}}
\def\scdh{{\text{scdh}}}
\def\zar{{\text{zar}}}
\def\red{{\text{red}}}
\def\etale{\'etale~}
\def\nor{\mathrm{nor}}

\def\cK{\mathcal K}
\def\cKH{\mathcal KH}

\def\Proj{\operatorname{Proj}}

\def\Hom{\operatorname{Hom}}

\def\Spec{\operatorname{Spec}}
\def\MSpec{\operatorname{MSpec}}
\def\MProj{\operatorname{MProj}}

\def\MSch{{\mathrm{MSch}}}
\def\Top{\operatorname{Top}}
\def\map#1{{\buildrel #1 \over \lra}}

\def\lra{\longrightarrow}
\def\into{\hookrightarrow} 
\def\onto{\twoheadrightarrow}

\newcommand{\Q}{\mathbb{Q}}

\newcommand{\bP}{\mathbb{P}}

\newcommand{\bF}{\mathbb{F}}
\newcommand{\A}{\mathbb{A}}
\newcommand{\bH}{{\mathbb{H}}}

\newcommand{\R}{{\mathbb{R}}}
\newcommand{\F}{\mathbb{F}}

\newcommand{\Z}{\mathbb{Z}}
\newcommand{\N}{\mathbb{N}}

\newcommand{\fp}{{\mathfrak p}}
\newcommand{\fm}{{\mathfrak m}}
\newcommand{\fq}{{\mathfrak q}}

\newcommand{\comment}[1]{}

\numberwithin{equation}{section}

\theoremstyle{plain} 
\newtheorem{thm}[equation]{Theorem}
\newtheorem{cor}[equation]{Corollary}
\newtheorem{lem}[equation]{Lemma}
\newtheorem{prop}[equation]{Proposition}

\theoremstyle{definition}
\newtheorem{defn}[equation]{Definition}
\newtheorem{ex}[equation]{Example}
\newtheorem{construction}[equation]{Construction}

\theoremstyle{remark}
\newtheorem{rem}[equation]{Remark}

\newtheorem{substuff}{\bf Remark}[equation]
\newtheorem{subrem}[substuff]{Remark} 
\newtheorem{subex}[substuff]{\bf Example} 
\newtheorem{obs}[substuff]{Observation}

\begin{document}

\title{Toric Varieties, Monoid Schemes and $cdh$ descent}
\date{\today}

\author{G. Corti\~nas}
\thanks{Corti\~nas' research was supported by Conicet
and partially supported by grants UBACyT W386, PIP
112-200801-00900, and MTM2007-64704 (Feder funds).}
\address{Dept.\ Matem\'atica-Inst. Santal\'o, FCEyN, Universidad de Buenos Aires,
Ciudad Universitaria, (1428) Buenos Aires, Argentina}
\email{gcorti@dm.uba.ar}

\author{C. Haesemeyer}
\thanks{Haesemeyer's research was partially supported by NSF grant DMS-0966821}
\address{Dept.\ of Mathematics, University of California, Los Angeles CA
90095, USA}
\email{chh@math.ucla.edu}

\author{Mark E. Walker}
\thanks{Walker's research was partially supported by NSF grant DMS-0601666.}
\address{Dept.\ of Mathematics, University of Nebraska - Lincoln,
  Lincoln, NE 68588, USA}
\email{mwalker5@math.unl.edu}

\author{C. Weibel}
\thanks{Weibel's research was supported by NSA and NSF grants.}
\address{Dept.\ of Mathematics, Rutgers University, New Brunswick,
NJ 08901, USA} \email{weibel@math.rutgers.edu}

\begin{abstract}
We give conditions for the Mayer-Vietoris property to hold for
the algebraic $K$-theory of blow-up squares of toric varieties 
and schemes, using the theory of monoid schemes.
These conditions are used to relate algebraic $K$-theory to
topological cyclic homology in characteristic $p$.
To achieve our goals, we develop many notions for monoid schemes
based on classical algebraic geometry,
such as separated and proper maps and resolution of singularities.
\end{abstract}
\subjclass[2000]{19D55, 14M25, 19D25, 14L32}

\keywords{abelian monoids, monoid schemes toric varieties,
algebraic $K$-theory}

\maketitle

\section*{}

The goal of this paper is to prove Haesemeyer's Theorem \cite[3.12]{HKH}
for toric schemes in any characteristic.
It is proven below as Corollary \ref{MainCor}.

\begin{thm} \label{IntroThm}
Assume $k$ is a commutative regular noetherian 
ring containing an infinite field 
and let $\cG$ be a presheaf of spectra defined on the category of schemes
of finite type over $k$. If $\cG$ satisfies the Mayer-Vietoris property
for Zariski covers, finite abstract blow-up squares, and blow-ups along
regularly embedded closed subschemes, then $\cG$ satisfies the 
Mayer-Vietoris property for all abstract blow-up squares of 
toric $k$-schemes obtained from subdividing a fan.
\end{thm}

The application we have in mind is to understand the relationship between
the algebraic $K$-theory $K_*(X)=\pi_*\cK(X)$ and topological cyclic homology
$TC_*(X)=\{\pi_*TC^\nu(X,p)\}$ of a toric scheme over a 
regular ring of characteristic $p$
(and in particular of toric varieties over a field of characteristic $p$).
Thus we consider the presheaf of homotopy fibers $\{\cF^\nu(X)\}$
of the map of pro-spectra from $\cK(X)$ to $\{TC^\nu(X,p)\}$.
Work of Geisser-Hesselholt \cite[Thm.\,B]{GHvanish}, \cite{GHbirel}
shows that this homotopy fiber (regarded as a pro-presheaf of spectra)
satisfies the hypotheses of Theorem \ref{IntroThm} and hence
a slight modification of the proof of
our theorem implies that it satisfies the Mayer-Vietoris property
for all abstract blow-up squares of toric schemes.
We will give a rigorous proof of this in Corollary \ref{cor:pro-hcart} below.

One major tool in our proof will be a theorem of Bierstone-Milman
\cite{BM} which says that the singularities of a toric variety 
(or scheme) can be resolved by a sequence of blow-ups $X_C\to X$
along a center $C$ that is a smooth, equivariant closed
subscheme of $X$ along which $X$ is normally flat. If one only had
to consider toric schemes, this would allow one to use
Haesemeyer's original argument to prove Theorem \ref{IntroThm}, since
toric schemes over a regular ring are normal and Cohen-Macaulay. However, examples
show that the blow-up of a toric scheme along a smooth center
(even a point) can be non-normal. Thus, even starting with a
toric scheme, the tower of blow-ups constructed by Bierstone-Milman
will often involve non-normal schemes with a torus action.
The proof of our theorem requires us to work with a larger class
of schemes, one containing all the schemes in this tower.
Beyond this, we need a class of schemes which is
closed under passage to (possibly non-reduced) equivariant closed subschemes,
pullbacks and blow-ups.

It turns out that all these operations may be lifted to the category
of {\em monoid schemes} of finite type, and that the realizations of 
monoid schemes over a commutative regular ring $k$ containing a field 
form a class of schemes with the above-mentioned properties.
The $k$-realization of an affine monoid scheme is a scheme of the form
$\Spec k[A]$, with $A$ an abelian monoid; the $k$-realization of
a monoid scheme (Definition \ref{defn:realizn}) is a scheme over $k$
which is covered by affine open subschemes of this form,
with homomorphisms of the underlying monoids inducing the gluing maps
between these open subschemes.

To achieve our goals, it is easier to work directly with the category of
monoid schemes, and Sections \ref{sec:monoids}--\ref{sec:bc} of
this paper are devoted to a introduction to monoid schemes.
Toric monoid schemes are introduced in Section 4
and the relation to toric varieties is carefully described.
In Sections \ref{sec:realizn} and \ref{sec:smooth}, we prove that the
$k$-realization functor preserves limits and show that many
monoid scheme-theoretic properties translate well into algebraic geometry.
Projective monoid schemes, blow-ups and proper maps are introduced in
Sections \ref{sec:blowup} and \ref{sec:proper}. After introducing the
technical notion of pctf monoid schemes in Section \ref{sec:pctf},
birational maps and resolution of singularities are given in Sections
\ref{sec:birat} and \ref{sec:ROS}.

The last part of this paper (Sections \ref{sec:cdh}--\ref{sec:main})
is devoted to the notion of cohomological descent
(Definition \ref{def:cdh-descent}),
the proof of our Main Theorem \ref{IntroThm} and its application to
algebraic $K$-theory and topological cyclic homology.

As far as the authors are aware, this paper presents the first attempt
at a systematic study of geometric properties of monoid schemes within
the category of monoid schemes, and the relationship of these with the
geometric properties of their realizations. The idea of a monoid
scheme itself goes back at least to Kato \cite{Kato}, and general
definitions were given by Deitmar in \cite{Deitmar} and (under the
name $\mathcal{M}_0$-schemes) by Connes, Consani and Marcolli in
\cite{CC}. Deitmar studies notions of flatness and \'etaleness for
monoid schemes, and introduces discrete valuation monoids.  
New in this paper is our
systematic investigation of separatedness, properness, general
valuation monoids and the valuative criteria, projectivity and blowing
up, and the introduction of a class of monoid schemes (the above
mentioned pctf monoid schemes) with better formal properties than only
those given by fans, yet avoiding the worst pathologies of
non-cancellative monoids.

\section{Monoids}\label{sec:monoids}

Since we know of no suitable reference for the facts we need
concerning monoids and their prime spectra, we begin with a
short expos\'e of this basic material.

Unless otherwise stated, a {\em monoid} in this paper is a pointed
abelian monoid; i.e., an abelian monoid object in the 
symmetric monoidal category of pointed sets with smash product 
as monoidal product. More explicitly, a monoid is a pointed set $A$ with
basepoint $0$, equipped with a pairing $\mu: A \smsh A \to A$ (written
$\mu(a,b)=ab$) that is associative and commutative and has an identity
element $1$. The basepoint is unique because it is characterized by 
the property that $0a=0$ for all $a\in A$.
For example, if $R$ is a commutative ring, then
forgetting addition gives a monoid $(R, \times)$ of this
type.
Sometimes $+$ notation is used for $\mu$, for example in
applications to toric varieties; in these cases we write $0$ for the
identity element, and $\infty$ for the basepoint.

We can convert any unpointed abelian monoid $B$
into a pointed abelian monoid $B_*$ by adjoining a basepoint.
Neither the zero monoid $\{0\}$ nor the monoid $\{0,1,t\}$ with
$t^2=0$ are of this form.

A morphism of monoids is a map of pointed sets
preserving the multiplicative identity and multiplication.
The initial monoid is $S^0=\{0,1\}$ with $1\cdot1=1$,
and the initial map $\iota_A: S^0 \to A$ is such that the
identity on $A$ equals the composition
\[\xymatrix{\
A \ar[r]^{\cong} &S^0\smsh A \ar[r]^{\iota_A\smsh id}
                 &  A\smsh A \ar[r]^{\mu}& A.}
\]

\medskip
\paragraph{\em Localization}
A \emph{multiplicatively closed} subset $S\subset A$ is a subset 
contaning $1$ and closed under multiplication.
Given a multiplicatively closed subset $S$ of $A$,
the {\it localization} $S^{-1}A$ consists of equivalence classes of
fractions of the form $\frac{a}{s}$ with $a \in A$ and $s \in S$. As
usual, $\frac{a}{s} = \frac{a'}{s'}$ if and only if $as's'' =
a'ss''$ for some $s'' \in S$, and the operation in $S^{-1}A$ is 
given by multiplication of fractions. There is a canonical monoid
homomorphism $A\to S^{-1}A$ sending $a$ to $\frac{a}{1}$, and $a,b
\in A$ are mapped to the same element of $S^{-1}A$ if and only if
$as = bs$ for some $s \in S$.

An {\it ideal} $I$ in a monoid $A$ is a pointed subset such that
$AI\subseteq I$. If $I \subset A$ is an ideal, $A/I$ is the monoid
obtained by collapsing $I$ to $0$ --- i.e., it is canonically
isomorphic to $(A \setminus I) \cup \{0\}$ with the unique
multiplication rule that makes the canonical surjection $A \onto A/I$
into a morphism of monoids. More generally, any surjective
homomorphism of monoids $A\to B$ is the quotient by a
\emph{congruence}, i.e. an equivalence relation compatible with the
monoid operation. 

Every non-zero monoid $A$ has a unique maximal ideal
(written $\fm_A$), namely the complement of the submonoid of units
$U(A) := \{ a \in A \, | \, ab =1 \text{ for some $b$ }\}$.
We say that a monoid morphism $g: A\to B$ is {\it local} if
$g(\fm_A)\subseteq\fm_B$ or, equivalently, if $g^{-1}(U(B))\subseteq U(A)$.

\goodbreak
A {\it prime ideal} is a proper ideal $\fp$ ($\fp\ne A$) whose 
complement $S=A\backslash\fp$ is closed under multiplication; 
in this case we write $A_\fp$ for the localization $S^{-1}A$.
The {\em dimension} of $A$ is the supremum of the lengths of all
chains of prime ideals, and the {\em height} of $\fp$
is the dimension of $A_\fp$. Since the intersection of an arbitrary 
chain of primes is prime, every prime ideal contains a minimal prime ideal
(by Zorn's Lemma).

\begin{lem}\label{lem:local}
For every multiplicatively closed subset $S$ of $A$ with $0\not\in S$,
there is a prime ideal $\fp$ of $A$ such that $S^{-1}A=A_\fp$.
\end{lem}

\begin{proof}
Since $S^{-1}A$ is a non-zero monoid it has a maximal (proper) ideal $\fm$; the
inverse image of $\fm$ in $A$ is a prime ideal $\fp$. Let $T$ denote
$A\setminus\fp$; then $S\subset T$ and any $t\in T$ is a unit in
$S^{-1}A$. Hence there are homomorphisms $S^{-1}A\to T^{-1}A=A_\fp$
and $T^{-1}A\to S^{-1}A$ covering the identity of $A$. Hence both
composites $S^{-1}A\to S^{-1}A$ and $A_\fp\to A_\fp$ are identity maps,
by the universal property of localization. \comment{ so there is a
$u\in A$ so that $tu=s\in S$. Therefore if $at=0$ in $A$ then $as=0$
and $a/t=au/s$. This shows that the canonical map $S^{-1}A\to
T^{-1}A=A_\fp$ is an isomorphism.}
\end{proof}

We let $\MSpec(A)$ denote the set of prime ideals of $A$; it is a
topological space when equipped with its Zariski topology, in which
closed subsets are those of the form $V(I) = \{ \fp \, | \, I
\subset \fp \}$ for an ideal $I$ of $A$.
The principal open subsets
$$
D(s)=\{\fp\in\MSpec(A) \, | \, s\notin\fp\} = \MSpec(A[1/s])
$$
form a basis for the Zariski topology. The space $\MSpec(A)$ is
quasi-compact, since any open $D(s)$ containing
the unique maximal ideal $\fm_A$ must have $D(s)=\MSpec(A)$.

There is a sheaf of monoids
$\cA$ on $\MSpec(A)$ whose stalk at $\fp$ is $A_\fp$; if $U$ is open
then $\cA(U)$ is the subset of $\prod_{\fp\in U} A_\fp$ consisting of
elements which locally come from some $S^{-1}A$. Explicitly,
\[
\cA(U)=\{a\in\prod_{\fp\in U}A_\fp:
(\forall \fp\in U)(\exists s\notin \fp, x\in A)(\forall \fq\in U) 
s\notin \fq\Rightarrow a_\fq=\frac{a}{s}\}.
\]
In particular, $A=\cA_{\fm_A}$, and $\cA(D(s))=A[1/s]$. 
More generally any ideal $I$ of $A$ determines
a sheaf $\cI$ on $\MSpec A$ by
\begin{equation*}
\cI(U)=\{a\in \cA(U):(\forall\fp\in U)\ \ a_\fp\in A_\fp\cdot I\}
\end{equation*}

\begin{ex}\label{affine}
The free (abelian) pointed monoid on the set $\{ t_1,\dots, t_n\}$
is the multiplicative monoid $F_n$ consisting of all monomials 
in the polynomial ring $\Z[t_1,...,t_n]$ (together with $0$). 
Each of the $2^n$ subsets of $\{ t_1,...,t_n\}$ 
generates a prime ideal $\fp$, and every prime ideal of $F_n$ 
has this form. We write $\A^n$ for $\MSpec(F_n)$.
\end{ex}

If $A\to B$ is a morphism of monoids, then the inverse image of a
prime ideal is a prime ideal, and we have a continuous map
$\MSpec(B)\to\MSpec(A)$. If $I$ is an ideal of $A$ then
$\MSpec(A/I)\to\MSpec(A)$ is a closed injection onto $V(I)$.
If $S$ is multiplicatively closed in $A$ then
either $S^{-1}A = 0$ (in which case $\MSpec S^{-1}A=\emptyset$) or 
$S^{-1}A=A_\fp$ for some $\fp$ (Lemma \ref{lem:local}); in either case
$\iota:\MSpec(S^{-1}A)\to\MSpec(A)$ is an injection onto the
set of primes that are disjoint from $S$.
The restriction $\iota^{-1}(\cA)$ to
this subset is the sheaf of monoids on $\MSpec(A_\fp)$.

Recall that a point $x_1$ of a topological space $X$ is called
\emph{generalization} of a point $x_0$ (and $x_0$ is called a
\emph{specialization} of $x_1$) if $x_0$ is in the closure of
$x_1$. For example, if $\fp,\fq\in\MSpec A$, then $\fp$ generalizes
$\fq$ if and only if $\fp\subset \fq$.

\begin{lem}\label{lem:yesopen}
Let $\fp$ be a prime ideal in a monoid $A$.
Then $\MSpec(A_\fp)\to\MSpec(A)$ is an injection,
closed under generalization,
and the following are equivalent:
\\(i) $\MSpec(A_\fp)$ is open in $\MSpec(A)$.
\\(ii) $\MSpec(A_\fp) =  D(s)$ for some $s\in A$.
\\(iii) There is an $s\in A$ such that  $A_\fp=A[1/s]$.
\end{lem}

\begin{proof}
The first assertion was observed above.
Since $D(s)=\MSpec(A[1/s])$, (iii) is equivalent to  (ii), a special case of (i).
Conversely, suppose that $U=\MSpec(A_\fp)$ is the complement of $V(I)$
for some ideal $I$ of $A$. Then $U=\cup_{s\in I}D(s)$. In particular,
there is an $s$ in $I$ such that $\fp\in D(s)$. But then $U\subseteq D(s)$
and hence $U=D(s)$.
\end{proof}

\begin{ex}\label{ex:notopen}
Let $A$ be the free pointed abelian monoid generated by the infinite set
$\{t_1, t_2, \dots\}$. If $\fp$ is the prime ideal generated by some 
finite subset of the $t_i$'s then $\MSpec(A_\fp)$ cannot be open in 
$\MSpec(A)$. Indeed, if it were open then by Lemma \ref{lem:yesopen} it
would have the form $D(s)$ for some element $s\in A$. But any $s$
involves only a finite number of variables, so the prime ideal
$t_jA$ belongs to $D(s)$ for infinitely many $t_j\not\in\fp$. 
In particular, $D(s)$ cannot be contained in $\MSpec A_\fp$. 
\end{ex}

\begin{lem}\label{lem:poset}
If $A$ is finitely generated as a monoid, then $\MSpec(A)$ is a
finite partially ordered set. If $S$ is a multiplicative subset of $A$,
then $S^{-1}A$ is also finitely generated, and $\MSpec(S^{-1}A)$
is open in $\MSpec(A)$.
\end{lem}

\begin{proof}
Suppose $A$ is generated by $x_1,\dots,x_m$. Then for any prime ideal $\fp$,
the multiplicative subset $S = A\setminus\fp$ is generated by
$\{x_i \,|\, x_i\notin\fp\}$. Indeed, if $s \in S$, then
$s = \prod_i x_i^{e_i}$ with $e_i = 0$ whenever $x_i \in \fp$.
Thus $A$ has at most $2^m$ prime ideals.

By Lemma \ref{lem:local}, we may assume $S =
A\setminus \fp$ for some prime $\fp$. If $s$ is the product of the
generators of $S$, then $A_\fp=A[1/s]$. By Lemma \ref{lem:yesopen},
$\MSpec(A_\fp)$ is open.
\end{proof}

We say $A$ is {\em cancellative} if for $a,b,c\in A$ the conditions 
$ab=ac$ and $a\ne0$ together imply that $b=c$.
In this case, the unpointed monoid $A\setminus \{0\}$ injects into 
its group completion and $\{0\}$ is the unique minimal prime ideal of $A$.
We define the {\em  pointed group completion} of $A$ to be
the pointed monoid $A^+$ obtained by adjoining a basepoint to the usual
group completion of the unpointed monoid $A\setminus \{0\}$.
Note that $A$ is a pointed submonoid of $A^+$, and  that $A^+$ is
the localization $A_{\{0\}}$ of $A$ at the minimal prime ideal.

We say $A$ is {\em torsionfree} if whenever $a^n=b^n$ for $a,b\in A$
and some $n\ge1$, we have $a=b$. The monoid $\{0,\pm1\}$ is
cancellative but not torsionfree. If $A$ is cancellative and
$A^+\setminus\{0\}$ is a torsionfree abelian group, then $A$ is torsionfree.

An element is {\em nilpotent} if $a^n=0$ for some $n$, and
the {\em nilradical} of $A$ is the set $\nil(A)$ of nilpotent elements.
It is easy to prove (using Zorn's lemma as in ring theory), that
$\nil(A)$ is the intersection of the minimal prime ideals of $A$. We say that
$A$ is {\em reduced} if $\nil(A)=0$, and set $A_\red=A/\nil(A)$.

Any closed subset $Z$ of $X=\MSpec(A)$ defines a largest ideal $I$ such
that $Z=V(I)$, and $A/I$ is a reduced monoid.
Indeed, if $Z=V(I_0)$ then $A/I=(A/I_0)_\red$; $I$ is the intersection
of the prime ideals containing $I_0$.
Anticipating Lemma \ref{equivclosure}, we write $\Bar{Z}^\eq$ for
$\MSpec(A/I)$ and call it the {\it equivariant closure} of $Z$ in $X$.
For example, $\Bar{X}^\eq$ is $\MSpec(A_\red)$. Another important special case is when
$Z=\{\fp_1,\dots,\fp_l\}$ is a set of prime ideals of $A$; in this case
$\Bar{Z}^\eq=\MSpec(A/\cap\fp_i)$.

\begin{defn}\label{def:norm}
The {\em normalization} of a cancellative monoid $A$ is
defined to be the submonoid
\begin{equation*}
A_{\nor} = \{\alpha\in A^+ \,|\,\alpha^n\in A \text{ for some $n\geq 1$}\}
\end{equation*}
of $A^+$.
We say that $A$ is {\em normal} if it is cancellative and $A=A_{\nor}$.
The normalization of $S^{-1}A$ is $S^{-1}A_\nor$.
If $A$ is torsionfree then so is $A_\nor$.
\end{defn}

\begin{subrem}\label{normzn-homeo}
If $A$ is cancellative then $\MSpec(A_\nor)\to\MSpec(A)$
is a topological homeomorphism. Indeed, if $\fp$ is a prime ideal of $A$
then $\fp_\nor :=\{ b\in A_\nor\,|\,(\exists n)b^n\in \fp\}$
is a prime
ideal of $A_\nor$ and $\fp =\fp_\nor \cap A$. It is easily seen that
every prime ideal of $A_\nor$ has the form $\fp_\nor$ for some $\fp$.
\end{subrem}

\begin{subrem}\label{norm-modP}
If $A$ is normal and $\fp$ is a prime ideal, then $A/\fp$ is also normal.
Indeed if $x,y\in A$ and $s\in A\setminus \fp$ are such that 
$x^n$ and $s^ny$ are mapped to the same element of $A/\fp$, 
then either $x^n=s^ny$ in $A$ or $x,y\in\fp$. Since $A$ is assumed normal,
it follows that either $x\in\fp$ or there is a $z\in A$ such that 
$x=sz$ in $A$.
\end{subrem}

More generally, let $f:A\to B$ be a morphism of monoids.
We say that $f$ is {\em integral} if for every $b \in B$ there is
an integer $n \geq 1$ such that $b^n$ lies in the image of $A$, and
we say that $f$ is {\it finite} if there exist
$b_1,\dots,b_n\in B$ ($n\ge 1$) such that $B=\bigcup_iAb_i$.
The normalization $A\to A_\nor$ is integral but not always finite.

\begin{lem}\label{integralvfinite}
Let $A\map{f} B$ be a monoid morphism
with $B$ finitely generated over $A$.
\begin{enumerate}
\item[i)] If $f$ is integral, then $f$ is finite.
\item[ii)]   If $f$ is finite and $B$ is cancellative, then $f$ is
  integral.
\end{enumerate}
\end{lem}

\begin{proof}
Choose a surjection $A[t_1,\dots,t_n]\onto B$, with the $t_i$ mapping
onto generators $b_i$ of $B$ over $A$.
If $f$ is integral, then there is an $m$ such that $b_i^m$ is in
the image of $A$ for all $i$; thus every element of $B$ can be written
as a product $f(a)c_j$, where $a\in A$ and $c_j$ is a monomial on the
$b_i$ with exponents $\le m$. This proves i).

Next assume that $f$ is finite and that $B$ is cancellative.
Let $b_1,\dots,b_n\in B$ be such that $B=\bigcup_i Ab_i$.
For each $i$, we choose an index $\pi(i)$ and $a_i \in A$
such that $b_i^2= a_ib_{\pi(i)}$; then $\pi$ is a map from the
finite set $\{1,\dots,n\}$ to itself.
For each fixed $i$, the iterates $\pi^r(i)$ cannot all be distinct,
so there exist $s \geq 1$ and $r \geq 1$ such that $j=\pi^r(i)$ satisfies
$\pi^{s}(j)=j$.
Hence there is an $a\in A$ and $m \geq 1$ such that $b_j^{m}=ab_j$.
Because $B$ is cancellative, this implies that $b_j^{m-1}=f(a)$.
Thus $b_j$ and hence $b_i$ is integral over $A$, as required.
\end{proof}

\begin{subrem}
The hypothesis that $B$ be cancellative in part ii) of Lemma
\ref{integralvfinite} is necessary.  For example, the monoid $B$
generated by $x,y$ subject to $y^2=xy$ contains the free monoid $A$
generated by $x$; the extension $A\subset B$ is finite but not integral.
\end{subrem}

\goodbreak
For a pointed set $X$ and commutative ring $k$, $k[X]$ denotes the
free $k$-module on $X$, modulo the summand indexed by the base point of $X$.
If $A$ is a pointed monoid, $k[A]$ is a ring in the
usual way, with multiplication given by the product rule for $A$. If
$B$ is an unpointed monoid, $k[B_*]$ coincides with the usual
monoid ring for $B$ with $k$ coefficients.
If $I$ is an ideal of the monoid $A$ then
$k[I]$ is an ideal of the ring $k[A]$, and $k[A/I]=k[A]/k[I]$.
If $I$ is prime, $k[I]$ need not be a prime ideal.

The category of pointed monoids has all small colimits. For example,
the coproduct of $A_1$ and $A_2$ is the smash product $A_1\smsh A_2$; 
the maps from $A_1$ and $A_2$ to $A_1\smsh A_2$ send 
$a_1$ to $a_1\smsh1$ and $a_2$ to $1\smsh a_2$.
The functor $A\mapsto k[A]$ preserves colimits since it has a
right adjoint, sending an algebra $R$ to $(R,\times)$, the underlying
multiplicative monoid of $R$; in particular, the natural map 
$k[A_1]\otimes_kk[A_2]\to k[A_1\smsh A_2]$ is an isomorphism.
More generally, the pushout $A_1\smsh_C A_2$ of a diagram
\begin{equation}
\label{pushout} \xymatrix{\
C \ar[r]^f \ar[d]^g & A_2 \ar@{..>}[d] \\ A_1 \ar@{..>}[r] & A_1\smsh_CA_2}
\end{equation}
is the quotient of $A_1\smsh A_2$  by the congruence generated by
$(a_1f(c),a_2) \sim (a_1,g(c)a_2)$.
Note that $k[A_1\smsh_C A_2]\cong k[A_1]\otimes_{k[C]}k[A_2]$.

\begin{lem}\label{primes-in-smash}
Every prime ideal $\fp$ of $A_1\smsh A_2$ has the form
$\fp_1\smsh A_2\cup A_1\smsh\fp_2$ for unique prime ideals 
$\fp_1$ and $\fp_2$. Explicitly, $\fp_i$ is the inverse image of $\fp$ 
under the canonical inclusion $A_i\to A_1\smsh A_2$.
\end{lem}

\begin{proof}
Given a prime ideal $\fp$ of $A_1\smsh A_2$, set $\fp_1=\fp\cap A_1$,
$\fp_2=\fp\cap A_2$ and $\fq=\fp_1\smsh A_2\cup A_1\smsh\fp_2$.
Then $\fq$ is prime because its
complement is $(A_1\backslash\fp_1)\times(A_2\backslash\fp_2)$,
which is multiplicatively closed. Clearly $\fq\subseteq\fp$; to see
that $\fq=\fp$, consider an element $a_1\smsh a_2$ of $\fp$. As $\fp$
is prime, either $a_1\smsh1$ or $1\smsh a_2$ is in $\fp$. In the first
case, $a_1\in\fp_1$ so $a_1\smsh a_2$ is in $\fp_1\smsh A_2\subseteq\fq$;
in the second case, $a_2\in\fp_2$ so $a_1\smsh a_2$ is in
$A_1\smsh\fp_2\subseteq\fq$.
\end{proof}

\goodbreak
\begin{ex}\label{extended}
If $T$ is the free monoid on one element $t$, then $A\smsh T$ is the
analogue of a polynomial ring over $A$, and $k[A\smsh T]=k[A][t]$.
For any prime ideal $\fp$ of $A$ there are exactly two primes of
$A\smsh T$ over $\fp$: the extended prime $\fp\smsh T$ and the prime
generated by $\fp$ and $t$ (i.e., $\fp\smsh T\cup A\smsh\{t^n:n\ge1\})$.
The map $\MSpec(A\smsh T)\to\MSpec(A)$ induced by the canonical
inclusion $A\to A\smsh T$ is both open and closed, because the image of 
$D(at^n)$ is $D(a)$ and the image of $V(I)$ is $V(I\cap A)$.
\end{ex}

\begin{prop}\label{primes-pushout}
Given a pushout diagram \eqref{pushout},
every prime ideal of $A_1\smsh_CA_2$ has the form
$\fp_1\smsh A_2\cup A_1\smsh\fp_2$ for unique prime ideals
$\fp_1$ in $A_1$, $\fp_2$ in $A_2$.

Moreover, the ideal $\fp_1\smsh A_2\cup A_1\smsh\fp_2$ of
$A_1\smsh_C A_2$ is prime if and only if $\fp_1$ and $\fp_2$
have a common inverse image in $C$.
\end{prop}
\goodbreak

\begin{proof}
If $\fp$ is a prime in $A_1\smsh_CA_2$, its inverse image in
$A_1\smsh A_2$ is prime; by Lemma \ref{primes-in-smash} it has
the form $\fp_1\smsh A_2\cup A_1\smsh\fp_2$, where $\fp_i\subset A_i$
are the inverse images of $\fp$. Since $A_1\smsh_CA_2$ is a quotient,
this proves the first assertion; because \eqref{pushout} commutes,
$\fp_1$ and $\fp_2$ have a common inverse image in $C$.

Conversely, suppose that $\fp_1$ and $\fp_2$ have a common inverse image
$\fq$ in $C$, and set
$S_1=A_1\backslash\fp_1$, $S_2=A_2\backslash\fp_2$ and
$I=\fp_1\smsh A_2\cup A_1\smsh\fp_2\subset A_1\smsh_C A_2$.
To see that the ideal $I$ is prime, it suffices to show that the image of
$S_1\times S_2$ in $A_1\smsh_CA_2$ is disjoint from $I$.
Since $\fp_1$ and $\fp_2$ are prime,
$a_1f(c)\in S_1$ if and only if $a_1\in S_1$ and $c\not\in\fq$, while
$g(c)a_2\in S_2$ if and only if $a_2\in S_2$ and $c\not\in\fq$.
It follows that $(a_1f(c),a_2)$ is in $S_1\times S_2$ if and only if
$(a_1,g(c)a_2)$ is. Thus $S_1\times S_2$ is closed under the
equivalence relation defining $A_1\smsh_CA_2$,
and its image in $A_1\smsh_CA_2$ is disjoint from $I$.
\end{proof}

\section{Monoid schemes}\label{sec:msch}

We will need to consider {\it monoid schemes}, sometimes known as
``schemes over the field with one element''.  These are the objects
which result by gluing together spectra of pointed monoids along
open subsets, and will be related to classical schemes in
Section \ref{sec:realizn}.
The theory of monoid schemes was developed by Kato
\cite{Kato}, Deitmar \cite{Deitmar}, Connes-Consani-Marcolli
\cite{CCM}, \cite{CC}, \cite{CC1}, etc. The survey \cite{LPL} by 
L\'opez Pe\~na and Lorscheid gives a nice overview of this notion 
and related ideas (but see Remark \ref{flaws} below).

A {\em monoid space} is a pair $(X,\cA_X)$ consisting of a topological
space $X$ and a sheaf $\cA_X$ of pointed abelian monoids on $X$.
A {\em morphism of monoid spaces} from $(X,\cA_X)$ to $(Y,\cA_Y)$
is given by a continuous map $f: X\to Y$ together with a morphism
of sheaves $f_\#:f^{-1}\!\cA_Y\to \cA_X$ on $X$ (or, equivalently, 
a morphism $f^\#:\cA_Y \to f_*\cA_X$ of sheaves on $Y$) that is 
{\em local} in the sense that the maps on stalks 
$\cA_{Y, f(x)} \to \cA_{X,x}$ are local morphisms of monoids, 
for all $x \in X$. By abuse of notation, we will often simply write
$X$ for the monoid space $(X,\cA_X)$.

The association $A\mapsto\MSpec(A)$ extends to a fully faithful
contravariant functor from monoids to monoid spaces, which we will
call $\MSpec$ by abuse of notation. An {\it affine monoid scheme}
is a monoid space isomorphic to $\MSpec(A)$ for some monoid $A$. 
A {\em monoid scheme} is a monoid space $(X,\cA)$ such that
every point has an open neighborhood $U$ such that $(U,\cA|_U)$ is 
isomorphic to an affine monoid scheme. If $(U,\cA|_U)\cong\MSpec A$ 
we shall often abuse notation and write $U=\MSpec A$. 
A morphism of monoid schemes is just a morphism of the underlying monoid
spaces. The {\em dimension} of a monoid scheme is the largest dimension of
its affine open neighborhoods.

\begin{lem}
Let $(X,\cA)$ be a monoid scheme. For any open $U\subseteq X$,
the monoid space $(U,\cA|_U)$ is a monoid scheme.

The scheme $(U,\cA|_U)$ is called the {\em open subscheme} of $X$ 
associated to $U$.
\end{lem}
\goodbreak

\begin{proof}
If $x\in U$ and $V=\MSpec(A)$ is an affine open neighborhood of $x$ in $X$,
$U\cap V$ is also open. Since $U\cap V$ is the union of basic open subschemes
$D(s)$ of $V$, $x$ has a neighborhood of the form $D(s)$,
and $D(s)=\MSpec(A[1/s])$ is affine.
\end{proof}

\goodbreak

We say that a monoid scheme is cancellative (resp., reduced, normal, ...) if
its stalks are cancellative monoids (resp., reduced, normal, ... monoids), 
or equivalently, if its monoids of sections are cancellative (resp., ...).

\begin{ex}
The projective line $\bP^1$ is obtained by gluing $\MSpec(\{ t^n,n\ge0\}_*)$
and $\MSpec(\{ t^{n},n\le0\}_*)$ along $\MSpec(\{ t^{n},n\in\Z\}_*)$.
This monoid scheme is connected, torsionfree and normal.
\end{ex}

\goodbreak\medskip

\paragraph{\em Partial order, maximal and minimal points}

Recall that the points of any topological space may be partially
ordered by the relation that $x\le y$ if and only if $y$ is in the
closure of $\{x\}$. In this way we can speak of maximal and minimal points.
The maximal points are the closed points; minimal points are also called
\emph{generic} points.
For the topological space $\MSpec(A)$ of a monoid $A$,
we have $\fp\le\fq$ if and only if  $\fp\subseteq \fq$.
Minimal points exist in any monoid scheme because,
as noted before \ref{lem:local},
every prime ideal contains a minimal prime ideal.

\begin{lem}\label{components}
Each cancellative monoid scheme $X$ decomposes as the disjoint union
of (closed and open) monoid subschemes $X_\eta$, each the closure
of a unique minimal point $\eta$ of $X$.
In particular, if $X$ is connected then it has a unique minimal point.
\end{lem}

\begin{proof}
For each minimal point $\eta\in X$, let $X_\eta$ denote the closure of 
$\eta$ in $X$. 
Given $x\in X$, choose an affine neighborhood $U_x=\MSpec(A)$ of $x$.
If $y$ is the point of $X$ corresponding to $\fm_A$, then $A=\cA_y$.
Since $A$ is cancellative, $U_x$ has a unique minimal point $\eta$,
so $U_x\subseteq X_\eta$. It follows that $X_\eta=\cup U_x$ is open
(and closed) in $X$, and that $X$ is the disjoint union of the $X_\eta$.
\end{proof}

\goodbreak
\begin{lem}\label{lem:maxima}
Let $X$ be a monoid scheme and $U\subseteq X$ an open subscheme. Then
the following are equivalent.
\\(i) $U$ is an affine monoid scheme.
\\(ii) $U$ has a unique maximal point.

If $X=\MSpec(A)$, every affine open subscheme is $\MSpec(A_\fp)$
for some $\fp$.
\end{lem}

\begin{proof}
Since monoids have unique maximal ideals, (i) implies (ii).
Conversely, suppose that $U$ has a unique maximal point $x$.
Note that $U = \{y\vert y\leq x\}$ by definition of the order relation.
If $\MSpec(A)$ is an affine open neighborhood of $x$, then
$U\subseteq\MSpec(A)$, so we may assume that $X=\MSpec(A)$.
In this case $U=\MSpec(A_x)$ by Lemma \ref{lem:yesopen}.
\end{proof}

\begin{defn}\label{def:closedimm}
Let $f:Y\to X$ be a map of monoid schemes. We say that $f$ is a
\emph{closed immersion} if it induces
a homeomorphism of $Y$ onto its image (equipped with the subspace topology),
and for every affine open subscheme $U=\MSpec(A)$ of $X$
(i) the open subscheme $V=U\times_X Y$ of $Y$ is affine (possibly empty) and
(ii) the map $\cA_X(U)\to \cA_Y(V)$ is surjective.
A {\em closed subscheme} of a monoid scheme $X$ is an isomorphism
class of closed immersions into $X$. Each closed subscheme is represented 
by a monoid scheme $(Z,\cA_Z)$ such that $f$ is a subspace inclusion 
$Z\subset X$.

A closed immersion $f:Y \to X$ is called {\em equivariant}
if in addition each such $\cA_X(U)\to \cA_Y(V)$ is the quotient by an ideal.
\goodbreak
\end{defn}

The terminology ``equivariant closed immersion'' comes from the theory of
toric varieties: the equivariant closed subschemes of a toric variety
are precisely those closed subschemes that are equivariant for the
action of the underlying torus. We will see in Section \ref{sec:toric}
that a toric variety has an associated toric monoid scheme, and that the
equivariant closed subschemes of the monoid scheme determine
equivariant closed subschemes of the toric variety.

\begin{ex}\label{ex:Zred}
Given a closed subset $Z$ of a monoid scheme $X$,
there is an equivariant reduced closed subscheme $Z_\red$
associated to $Z$, defined by patching; if $X=\MSpec(A)$ and $Z=V(I)$
then $Z_\red=\MSpec(A/I)_\red$.
\end{ex}

\begin{lem}\label{lem:immersion1}
Any surjection of monoids $A \map{\pi} B$ determines a closed
immersion $\MSpec(B) \subseteq \MSpec(A)$.
If $B=A/I$ then it is an equivariant closed subscheme.
\end{lem}

\begin{proof}
Set $Y=\MSpec(B)$ and $X=\MSpec(A)$. The map $\pi^*:Y\to X$
of underlying spaces is injective, since if $\fq_1\ne \fq_2$
then $\pi^{-1}(\fq_1)\ne\pi^{-1}(\fq_2)$. If $a\in A$, the image of
the basic open $D(\pi(a))\subseteq Y$ is $D(a)\cap \pi^*(Y)$. Thus
$Y$ is homeomorphic to $\pi^*(Y)$.

Let $U\subseteq\MSpec(A)$ be an affine open subscheme.
By Lemma \ref{lem:maxima} there is a prime $\fp$ of $A$ such that
$U=\MSpec(A_\fp)$; by Lemma \ref{lem:yesopen}, $U=D(s)$ for some $s$.
Hence $U\cap Y=D(\pi(s))=\MSpec(B[1/s])$, which is affine or empty.
Since $A[1/s]\to B[1/s]$ is onto, $Y\to X$ is a closed immersion.
\end{proof}

\begin{subrem}
A closed subscheme $Y \subset X$ need not determine a closed subset
of the underlying topological space. For example, the diagonal
embedding $\A^1\to\A^2$ is a closed immersion by Lemma
\ref{lem:immersion1}, but it is not topologically closed,
because it takes the generic point of $\A^1$ to the generic point of
$\A^2$ and the maximal point to the maximal point; the intermediate
points are not in the image.
\end{subrem}

\begin{defn} \label{qc_ideal}
If $(X,\cA)$ is a monoid scheme, a sheaf of ideals
$\cI$ is said to be \emph{quasi-coherent} if its restriction
to any affine open subscheme $U$ of $X$ is the sheaf associated 
to the ideal $\cI(U)$ of the monoid $\cA(U)$.
Given any closed immersion $i:Y\to X$,
the inverse image $\cI$ of $0$ under $\cA_X\to i_*\cA_Y$ is quasi-coherent.
Lemma \ref{lem:immersion1} shows that conversely any quasi-coherent sheaf
$\cI$ defines an equivariant closed immersion.
\end{defn}

\begin{lem} \label{equivclosure}
For any monoid scheme
$X$ and any subset $Z$ of the underlying poset, there is an
equivariant closed subcheme $\Bar{Z}^\eq$ of $X$ that contains
$Z$ and is contained in every other equivariant closed subscheme
of $X$ containing $Z$.
We call $\Bar{Z}^\eq$  the {\em equivariant closure} of $Z$ in $X$.

If $U$ is an open subscheme of $X$ then
$\Bar{Z}^\eq\cap U$ is $\overline{Z\cap U}{}^\eq$.
\end{lem}

\begin{proof}
We saw in Section \ref{sec:monoids} that if $Z$ is any subset of
$\MSpec(A)$, there is an equivariant closed subscheme
$\Bar{Z}^\eq=\MSpec(A/I)$ which contains $Z$ (and its closure),
and which is minimal with this property. Indeed, if the closure of $Z$
is $V(I_0)$, then $A/I=(A/I_0)_\red$.
Since $S^{-1}(A/I)=(S^{-1}A/S^{-1}I_0)_\red$, this construction patches
to give a general construction.
\end{proof}

\begin{subrem}\label{rem:equiclose}
If every point in $Z$ has height at least $i$ in $X$ then every point in
$\Bar{Z}^\eq$ has height at least $i$ in $X$. This follows from the
local description of $\Bar{Z}^\eq$.
\end{subrem}

\paragraph{\em Finite type}

We say that a monoid scheme has {\em finite type} if it admits a
finite open cover by affine monoid schemes associated to finitely
generated monoids. These monoid schemes are the analogues of
Noetherian schemes, just as finitely generated monoids are the analogues
of commutative Noetherian rings:
if $A$ is a finitely generated monoid then every ideal is
finitely generated, and $A$ has the ascending chain condition on
ideals. (The usual proof of the Hilbert Basis Theorem works.)

By Lemma \ref{lem:poset}, if $(X,\cA)$ is a monoid scheme of finite type,
then $X$ is a finite poset, with the poset topology.
The sheaf of monoids $\cA$ of a monoid scheme $X$ determines a 
(contravariant) functor $A$ from the poset $X$ to monoids, called the
{\it stalk functor} of $(X,\cA)$, sending $x$ to $\cA_x$.
It is useful to introduce the notion of a monoid poset as a context
for thinking about a stalk functor $A$.

A \emph{monoid poset} is a pair $(Y,B)$ consisting of a poset $Y$ and 
a contravariant functor $B$ from $Y$ to monoids.
There is a category of monoid posets; a morphism
$f:(X,A)\to (Y,B)$ of monoid posets is a poset map $g:X\to Y$ 
$A$ to the downward-closed subset $W(x)=\{y\in X\,|\, y\le x\}$.
There is a morphism of monoid posets
\begin{equation}\label{schematic}
\iota_x:(W(x),A|_{W(x)})\to F(\MSpec A(x))
\end{equation}
whose poset map sends a point $y$ of $W(x)$ to the inverse image 
$\fp_y\in\MSpec A(x)$ of the maximal ideal of $A(y)$ under $A(x)\to A(y)$; 
the maps $A(x)_{\fp_y}\to A(y)$ determine the natural
transformation $A\circ \iota_x \Rightarrow A|_{W(x)}$.
If the morphism \eqref{schematic} is an isomorphism for all $x\in X$, 
we will say that the monoid poset $(X,A)$ is {\em scheme-like} 
and (by abuse of notation) we will call $A$ a {\em stalk functor}.

We say that a monoid poset $(X,A)$ is of {\it finite type} if $X$ is a
finite poset and each $A(x)$ is a finitely generated monoid.
If $X$ is a monoid scheme of finite type, then $F(X)$ is a
monoid poset of finite type.
The following proposition shows that the stalk functor is always
enough to determine a monoid scheme of finite type.

\begin{prop}\label{P:stalk}
The functor $F(X,\cA)=(X,A)$ 
induces an equivalence between the full subcategory of 
monoid schemes of finite type and the full subcategory of 
scheme-like monoid posets $(X,A)$ of finite type. 
\end{prop}
\goodbreak

\begin{proof} 
If $(X,A)$ is a monoid poset, 
we may equip $X$ with the poset topology, and define the sheaf
$\cA$ on $X$ by the formula
\begin{equation*}
\cA(U)=\varprojlim_{x\in U} A(x).
\end{equation*}
Thus $G(X,A)=(X,\cA)$ is a monoid space.
It is clear from the formula for $\cA(U)$ 
that a morphism $(Y,B)\to (X,A)$ of monoid posets 
induces a morphism $G(Y,B)\to G(X,A)$ of monoid spaces.
Thus $G$ is a functor.
Because each $W(x)$ has $x$ as its maximal point, $\cA(W(x))=A(x)$.
Thus $F(G(X,A)$ is isomorphic to $(X,A)$.


If $(X,A)$ is scheme-like of finite type, then
$G(X,A)$ is a monoid scheme of finite type.
Conversely if $X$ is a monoid scheme of finite type
and $U$ is an affine open in $X$, we know by Lemma \ref{lem:maxima} that
there is a unique $x\in X$ such that $U=\MSpec(A(x))$ and hence
$A(x)=\cA(U)$. Given an open $U$ in $X$, any point $y$ in $U$
lies in an affine open $V\subset U$, and $V=\MSpec(A(x))$ for some $x\in U$
with $y\le x$ by Lemma \ref{lem:maxima}. It follows that $GF(X)\cong X$.
\end{proof}


A monoid scheme $(X,\cA)$ of finite type will often be specified by 
its monoid poset, 
{\it viz.,} $(X,A)$. To avoid confusion, we shall use roman
letters for stalk functors and script letters for sheaves.

\begin{rem}\label{rem:stalkf}
The proof of Proposition \ref{P:stalk} shows that any 
scheme-like monoid poset $(X,A)$ can be recovered from 
the monoid space $G(X,A)$ because $FG(X,A)\cong (X,A)$.
If $(X,\cA)$ is an arbitrary monoid scheme with stalk functor $A$, 
then the topology of $X$ may be coarser than the poset topology. 
However the argument of the proof of the proposition shows that 
we can recover $\cA$ from $A$ and the topological space
underlying $X$, using the formula 
$\cA(U)=\varprojlim_{x\in U} A(x)$.
\end{rem}

\goodbreak
\section{Basechange and separated morphisms} \label{sec:bc}

It is useful to simplify constructions using base-change. For this, we need
pullback squares in the category of monoid schemes.

There is a canonical morphism $\nu:X\to\MSpec(\cA(X))$ which is
universal for maps from $X$ to affine monoid schemes. It sends a point
$x$ to the preimage $\nu_x$ of the maximal ideal of $\cA_x$. The sheaf
homomorphism $\nu^\#$ is that induced by the canonical maps
$\cA(X)[1/s]\to \cA(\nu^{-1}(D(s))$. 
The universal property shows that the (contravariant) functor $X\mapsto
\cA_X(X)$ from monoid schemes to monoids is left adjoint to the functor 
$\MSpec$, i.e., that affine monoid schemes are a reflective subcategory
of all monoid schemes. It follows that $\MSpec$ converts 
pushouts of diagrams of monoids to pullbacks of diagrams in the 
category of all monoid schemes. In particular, for any pushout diagram 
of monoids \eqref{pushout}, the induced diagram is cartesian:
\[
\xymatrix{\MSpec(A_1\smsh_C A_2)\ar[d]\ar[r]& \MSpec A_2\ar[d]\\
           \MSpec A_1\ar[r]&\MSpec C.}
\]

\begin{prop}\label{prop:pb}
The pullback $X\times_SY$ of a diagram of monoid schemes
$$\xymatrix{
X\times_SY \ar@{..>}[r]  \ar@{..>}[d] & X \ar[d] \\
Y \ar[r] & S}$$
\\
exists in the category of all monoid schemes.
Its underlying topological space is the pullback
$X \times_S Y$ in the category of topological spaces.
\end{prop}

\begin{proof}
Existence of the pullback $X\times_SY$ is derived from the existence of
pullbacks of affine monoid schemes, just as for usual schemes
(\cite[Thm.\ 3.3]{Hart}).

To prove the assertion about underlying topological spaces, it
suffices to consider the affine case. Using the notation of
\eqref{pushout}, write $P$ for the pullback of $\MSpec(A_1)$ and
$\MSpec(A_2)$ over $\MSpec(C)$ in $\Top$. The canonical map
$f:\MSpec(A_1\smsh_C A_2)\to P$ is a continuous bijection by
Proposition \ref{primes-pushout}. To show that $f$ is a
homeomorphism, it suffices to show that it takes any basic open set
$D(s)$ to an open set of $P$. Write $s=s_1\smsh s_2$; then
$s\not\in\fp$ if and only if $s_1\smsh1, 1\smsh s_2\not\in\fp$. We
saw in Proposition \ref{primes-pushout} that if $\fp$ maps to
$(\fp_1,\fp_2)$ then $\fp=\fp_1\smsh A_2\cup A_1\smsh\fp_2$, and
that $s_1\smsh1\not\in\fp$ (resp., $1\smsh s_2\not\in\fp$) is
equivalent to $s_1\not\in\fp_1$ (resp., $s_2\not\in\fp_2$). This
shows that $f$ takes $D(s)$ to the open set $(D(s_1)\times
D(s_2))\cap P$, as required.
\end{proof}

\begin{subex}
The product $X\times Y$ is just the pullback
when $S$ is the terminal monoid scheme $\MSpec(S^0)$.
\end{subex}
\goodbreak

\begin{subrem}
Let $X$ and $Y$ be monoid schemes of finite type, over a common $S$.
Then the pullback $X\times_SY$ has finite type. Indeed, it has a finite
cover by affine opens of the form $\MSpec(A_1\smsh_C A_2)$, and in
each case $A_1\smsh_C A_2$ is finitely generated because $A_1$ and $A_2$ are.
\end{subrem}

\begin{ex}\label{intersection}
Proposition \ref{prop:pb} shows that given two closed subschemes $Z_1,Z_2$
of $X$,  the pullback $Z_1\times_X Z_2$ is a subscheme whose
underlying topological space is the intersection of the two
subspaces of $X$. More generally, given any family of closed
immersions $Z_i \into X$, we can form the inverse limit
$\lim Z_i \into X$ by patching the inverse limits on each
affine open $\MSpec(A)$, because the colimit of a family of
surjections $A\onto B_i$ exists and is a surjection.
\end{ex}

\goodbreak

\paragraph{\em Separated morphisms}

An important hypothesis in many theorems about monoid schemes, often
overlooked in the literature, is that they be separated.

\begin{defn}\label{def:separated}
A morphism $f:X\to S$ of monoid schemes is \emph{separated} if
the diagonal map $\Delta:X\to X\times_SX$ is a closed immersion.
We say that $X$ is separated if it is separated over $\MSpec(S^0)$
where we recall $S^0 = \{0,1\}$.
\end{defn}

Being separated is local on the base: if $S$ has an open cover $\{ U\}$
then $f$ is separated if and only if each $f^{-1}(U)\to U$ is separated.

\begin{lem}\label{lem:affsep}
If $A\to B$ is a morphism of monoids then
$\MSpec(B)\to \MSpec(A)$ is a separated morphism of monoid schemes.

In particular, closed immersions are separated.
\end{lem}

\begin{proof}
By Proposition \ref{prop:pb}, the diagonal map $\Delta$ corresponds to
the multiplication map $B\smsh_AB\to B$, which is surjective.
By Lemma \ref{lem:immersion1}, $\Delta$ is a closed immersion.
\end{proof}

\goodbreak
\begin{subrem} \label{notproper}
Example \ref{extended} shows that $X\times\A^1\to X$ is separated
and universally closed for every monoid scheme $X$. This shows that
``separated and universally closed" does not provide a good notion
of proper morphism of monoid schemes; we will discuss an appropriate
definition in Section \ref{sec:proper}.
\end{subrem}
\goodbreak

\begin{ex}\label{doubleorigin}
Here is an example of a monoid scheme which is non-separated.
Let $A$ and $B$ each be the free abelian monoid with two generators,
$F_2$ (see Example \ref{affine}).
Let $U$ be the open subset of each of $\MSpec(A)$ and $\MSpec(B)$ given
by removing the unique closed point (associated to the maximal ideal in
each monoid); explicitly $U=\{\langle t_1\rangle, \langle t_2\rangle,\{0\}\}$. Then we may glue $\MSpec(A)$ and $\MSpec(B)$ along $U$ to
form a monoid scheme $X$ of finite type. As a poset,
$X$ has five elements, two of which are maximal ---
the two copies of $\langle t_1,t_2\rangle$ --- and the rest are in $U$.

The $k$-realization of $X$ (defined in \ref{defn:realizn} below) is the
non-separated scheme given by the affine plane with the origin doubled.
\end{ex}

\begin{lem}\label{lem:sep}
A map $f:(X,\cA)\to (S,\cB)$ of monoid schemes
is separated if and only if for every $x_1,x_2$ in $X$ such that
$f(x_1)=f(x_2)$ and such that $\MSpec(\cA_{x_1})$ and $\MSpec(\cA_{x_2})$ are open,
either there is no lower bound for $\{x_1,x_2\}$ in the poset $X$ or
else there is a unique maximal lower bound $x_0=x_1\cap x_2$, and
$\cA_{x_1}\smsh_{\cB_{f(x_1)}} \cA_{x_2}\to \cA_{x_0}$ is onto.
\end{lem}

\begin{proof}
By Proposition \ref{prop:pb} and Lemma \ref{lem:maxima},
an affine open subset of $X\times_SX$
has the form $U=(U_1\times U_2)\cap(X\times_SX)$, where the maximal
point $(x_1,x_2)$ of $U$ determines the affine open subsets
$U_i=\MSpec(\cA_{x_i})$ of $X$. Since $\Delta^{-1}(U)=U_1\cap U_2$,
Proposition \ref{prop:pb} implies that $X\to\Delta(X)$ is a homeomorphism
and that the poset underlying $U_1\cap U_2$ is the subset
$\{z\in X\,|\, z\le x_1, z\le x_2\}$ of lower bounds for $\{x_1,x_2\}$.
If $U_1\cap U_2=\emptyset$, $\{x_1,x_2\}$ has no lower bound.

By Lemma \ref{lem:maxima}, $U_1\cap U_2$ is nonempty affine if and only if
it has a unique maximal element. Thus $\Delta$ is a closed immersion if and only if, in the above situation,
whenever $U_1\cap U_2$ is nonempty it is affine (and hence has a unique
maximal lower bound $x_0$), and
$\cA_{x_1}\!\smsh_C \cA_{x_2}\!\to\!\cA_{x_0}$ is onto,
where $s\!=\!f(x_1)\!=\!f(x_2)$ and $C=\cB_s$.
\end{proof}

\goodbreak

\begin{cor}\label{sepft}
If $X$ is a monoid scheme of finite type with stalk functor $A$, then
$X$ is separated if and only if whenever two points $x_1,x_2$ of $X$ have a
lower bound they have a greatest lower bound $x_1\cap x_2$,
and $A(x_1)\smsh A(x_2)\to A(x_1\cap x_2)$ is onto.
\end{cor}

\begin{proof}
Combine Lemma \ref{lem:sep} and Proposition \ref{P:stalk}.
\end{proof}

\begin{cor} \label{Cor524}
The intersection of two affine open subschemes of a separated monoid
scheme is affine.
\end{cor}

\begin{proof}
Suppose $X$ is a separated monoid scheme, 
with $U_1, U_2$ affine and open in $X$. Let $x_1, x_2$ be the unique
closed points of $U_1, U_2$. If $x_1$ and $x_2$ do not have a
common lower bound in $X$, then $U_1 \cap U_2 = \emptyset$.
Otherwise, by Lemma \ref{lem:sep},
they have a greatest lower bound, 
which is the
unique maximal point of $U_0=U_1\cap U_2$. By Lemma \ref{lem:maxima},
$U_0$ is affine.
\end{proof}

\goodbreak
\section{Toric monoid schemes}\label{sec:toric}

As observed by Kato \cite{Kato} and Deitmar \cite{Deitmar},
the fan associated to a toric variety produces a
monoid scheme. In this section we clarify this correspondence,
using the following definition.

\begin{defn}\label{def:toric}
A {\em toric monoid scheme} is a
separated, connected, torsionfree,
normal monoid scheme of finite type.

Recall that a {\em fan} consists of a free abelian group $N$ of
finite rank (written additively) together with a finite collection
$\Delta$ of strongly convex rational polyhedral cones $\sigma$ in
$N_\R$ (hereafter referred  to as just {\em cones}), satisfying the
conditions that (1) every face of a member of $\Delta$ is also a
member of $\Delta$ and (2) the intersection of any two members of
$\Delta$ is
a face of each. Here a \emph{strongly convex rational polyhedral cone}
is a cone with apex at the origin, 
generated by finitely many elements of $N$,
that contains no lines through the origin.

Note that $\Delta$ is a finite poset under containment; 
we now construct a monoid poset $(\Delta,A)$ and
use Proposition \ref{P:stalk} to define the associated monoid scheme.

\begin{construction}\label{ex:fan}
Given a fan $(N,\Delta)$, set $M=\Hom_\Z(N,\Z)$ and $M_{\R}=M\otimes\R$.
We define a contravariant
functor $A$ from $\Delta$ to monoids (written additively) by
\[
A(\sigma) = (\sigma^\vee \cap M)_*,
\quad \sigma^\vee = \{ m\in M_\R\,|\, m(\sigma)\ge0\}.
\]
Each such monoid is torsion-free, normal and finitely generated 
(Gordon's Lemma). If $\tau$ is a face of $\sigma$, then there is
an $m\in A(\sigma)$ such that $A(\tau)=A(\sigma)[-m]$. 
Hence by Lemma \ref{lem:local} there is a prime ideal $P_\sigma(\tau)$ of
$A(\sigma)$ such that $A(\tau)=A(\sigma)_{P(\tau)}$. 
By Proposition \ref{P:stalk} and Corollary \ref{sepft},
$A$ is the stalk functor of a toric monoid scheme $X(N,\Delta)$,
which by abuse of notation we write as
\[ 
X(\Delta) = (\Delta,A).
\]
Thus any fan $\Delta$ determines a toric monoid scheme in the sense of
Definition \ref{def:toric}.

A morphism of fans, from $(N,\Delta)$ to $(N',\Delta')$, is given by
a group homomorphism $\phi: N \to N'$ such that the image of each cone
in $\Delta$ under the induced map $N_\R \to N'_\R$ is contained in a
cone in $\Delta'$. Such a map of fans induces a poset map
$\Delta\to \Delta'$, sending $\sigma$ to the smallest cone $\sigma'$ in
$\Delta'$ that contains $\phi(\sigma)$, and precomposition with $\phi$
yields a natural transformation
$((\sigma')^\vee \cap M')_* \to(\sigma^\vee \cap M)_*$
of stalk functors. According to Proposition \ref{P:stalk},
this data determines a morphism of monoid schemes:
\[ 
X(\phi):X(\Delta) \to X(\Delta').
\]
If $\phi_1\ne \phi_2$ then $X(\phi_1)\ne X(\phi_2)$,
as $\phi_1^*\ne \phi_2^*$ on some $A(\sigma)$.
Thus we have a faithful functor $X$ from fans to toric monoid schemes.
\end{construction}

\goodbreak

\begin{ex}
For the cone $\sigma$ in the plane spanned by $(0,1)$ and $(1,-2)$,
$A(\sigma)=\sigma^\vee \cap M$ is the submonoid of $\Z^2$ spanned by
$\{(1,0),(1,1),(1,2)\}$. If $\Delta$ is the fan spanned by $\sigma$ and
its faces, then $(\Delta,A)=\MSpec A(\sigma)$.
\end{ex}

\goodbreak

If $(X,A)$ is a toric monoid scheme and $x\in X$, we will write
$M_x$ for the group completion of the unpointed monoid $A(x)\setminus\{0\}$.
Each $M_x$ is a torsionfree abelian group of finite rank. The groups
$M_x$ are all isomorphic, because $X$ has a unique minimal point $\eta$
by Lemma \ref{components}, and $M_x\to M_\eta=A(\eta)\setminus\{0\}$ 
is an isomorphism for all $x$.
\end{defn}

\goodbreak
\begin{thm}\label{why_toric}
The faithful functor $\Delta\mapsto X(\Delta)$ from fans to
toric monoid schemes, defined by Construction \ref{ex:fan}
has the following properties.
\begin{enumerate}
\item
Every toric monoid scheme $(X,A)$ is isomorphic to
$X(N,\Delta)$, where:
\begin{enumerate}
\item[a)] The lattice $N$ is the $\Z$-linear dual of $M=M_\eta$,
where $\eta$ is the unique minimal point of $X$.

\item[b)] The poset $\Delta$ of cones in $N_\R$ is isomorphic to 
the poset underlying $X$.
For each $x \in X$, the cone $\sigma_x$ in $N_\R$ is the dual cone of
the convex hull of $A(x) \setminus \{0\}$ in $M_\R$.
\end{enumerate}

\item For fans $(N,\Delta)$ and $(N',\Delta')$,  a morphism
$f:X(\Delta)\to X(\Delta')$ of monoid schemes is given by a
(necessarily unique) morphism of fans
if and only if $f$ maps the generic (i.e., minimal)
point $\eta$ of $X(\Delta)$ to the generic point $\eta'$ of $X(\Delta')$.
In this case, the map of fans $(N, \Delta) \to (N', \Delta')$ is given
by the $\Z$-linear dual of the group homomorphism 
$$
f_\eta^\#: M' = (A'(\eta') \setminus \{0\})
\to
(A(\eta)  \setminus \{0\})
= M.
$$
\end{enumerate}
\end{thm}

\begin{proof}
Throughout this proof, for a cancellative monoid $A$,  we write $A^o$
for the unpointed monoid $A \setminus \{0\}$, written additively,
and we identify each $A^o(x)$ with a submonoid of $M$. Let $(X,A)$
be a toric monoid scheme. We first show that $(N,\Delta)$ as defined in
the statement is a fan. For $x \in X$, let $\sigma_x^\vee \subset M_\R$
denote the convex hull  of $A^o(x)$ in $M_\R$. Note that this defines a 
cone $\sigma_x=(\sigma_x^\vee)^\vee$ in $N$ via the identification $N=N^{**}$.
The cone $\sigma_x^\vee$ is a rational polyhedral cone because it is 
spanned by a finite set $\{ a_i\}$ of generators of $A^o(x)$; 
the cone $\sigma_x$ is thus also a rational polyhedral cone, and it
 is strongly convex since $A^o(x)^+=M$.

To see that $A(x)=(\sigma_x^\vee\cap M)_*$,
let $b=\sum q_i a_i$ be an element of $\sigma_x^\vee\cap M$, written
as a positive $\Q$-linear combination of the $a_i$. Clearing denominators,
$nb$ is a positive $\Z$-linear combination of the $a_i$ for some
positive integer $n$ and hence is in $A^o(x)$.
Because $A(x)$ is normal, $b$ is in $A^o(x)$, as required.

If $\tau$ is a face of $\sigma_x$, it is defined by the
vanishing of some $m \in \sigma_x^\vee$. Clearing denominators and 
using again that $A(x)$ is normal, we
may assume $m \in A^o(x)$. By definition, $\tau$ is the set of linear
functionals on $M_\R$ that are non-negative on $A^o(x)[-m]$.
By Lemma \ref{lem:local}, $A(x)[-m]$
coincides with
$A(y)$ for some $y \leq x$, and thus the face $\tau$ is the element
$\sigma_y$ of $\Delta$.

If $x,y \in X$, we claim that the intersection $\sigma_x\cap\sigma_y$
is a cone in of $\Delta$.
Since $X$ is separated, $x$ and $y$ have a
unique greatest common lower bound, written $x \cap y$, and the map
$A(x)^o\times A(y)^o\to A^o(x\cap y)$ is surjective, by the additive version of Corollary \ref{sepft}; moreover because $X$ is cancellative, it is an isomorphism.
A linear functional on $M_\R$ is non-negative on $A^o(x)\times A^o(y)$ if and
only if it is non-negative on $A^o(x)$ and $A^o(y)$, and thus
we have the required identity:
\[
\sigma_x \cap \sigma_y = \sigma_{x \cap y}.
\]
Moreover, $\sigma_{x\cap y}$ is a face of both $\sigma_x$ and $\sigma_y$,
because by Lemmas \ref{lem:yesopen} and \ref{lem:poset} there are
$m_1,m_2$ such that
$A^o(\sigma_{x\cap y})=A^o(\sigma_x)[-m_1]=A^o(\sigma_y)[-m_2]$,
This proves that $\Delta$ is a fan.
\goodbreak

By Construction \ref{ex:fan}, the fan $(N,\Delta)$ determines a monoid scheme
$(\Delta, B)$. The bijection $\sigma: X \to \Delta$ ($x \mapsto\sigma_x$)
is order preserving, because if $x < y$ in $X$, then
$A^o(y)\subset A^o(x)\subseteq M$.
By construction, we have a natural isomorphism
$A(x) = (\sigma_x^\vee \cap M)_*=B(\sigma_x)$. 
This proves that $\sigma$ determines an isomorphism of monoid schemes, 
completing the proof of property 1).

Construction \ref{ex:fan}
shows that the condition in property 2) is necessary,
since a morphism of fans sends the zero cone to the zero cone.
Conversely, if $f(\eta) = \eta'$, then $f_\eta^\#$ induces a monoid
map $A'(\eta') = M'_* \to M_* = A(\eta)$; since any such map sends
units to units, it induces a group homomorphism $M' \to M$. Let
$\phi: N \to N'$ be the $\Z$-linear dual of this map.
Since for each $x \in X$, the map $f_x^\#$ is the
restriction of $f_\eta^\#$, it follows that $f=X(\phi)$, as desired.
\end{proof}

\begin{subrem}\label{flaws}
There are differing assertions in the literature
related to Theorem \ref{why_toric}.
Using a different definition of `toric variety'
it is claimed in \cite[Thm. 4.1]{Deitmar} that any connected cancellative
monoid scheme of finite type yields a toric variety, but not every such
``toric variety'' is associated to a fan. For example, $\MSpec$ of the cusp
monoid $C=\{ t^2,t^3,...\}_*$ yields the cusp.
In \cite[2.1]{LPL}, the flawed \cite[Thm. 4.1]{Deitmar} is used to
claim that the functor of Theorem \ref{why_toric} is an equivalence, 
under the weaker hypothesis that $A$ has no torsion; the cusp monoid is also 
a counterexample to the assertion in {\it loc.\,cit.}
\end{subrem}

We conclude this section with a description of separated
normal monoid schemes.
If $X$ is connected and cancellative, with minimal prime $\eta$,
then $M_\eta$ is a finitely generated abelian group. Therefore
there is a non-canonical isomorphism $M_\eta\cong M\times T$,
where $M$ is a free abelian group and $T$ is a finite torsion group.

\begin{prop} \label{toric-by-finite}
Any separated, connected, normal monoid scheme of finite type
decomposes as a cartesian product of monoid schemes
\[
X \cong (X, A) \times \MSpec(T_*),
\]
where $(X, A)$ is a toric monoid scheme and $T$ is a finite abelian group.
\end{prop}

\begin{proof}
If $\MSpec(A)$ is an affine open of $X$ then $A$ is a
submonoid of $A_\eta=(M\times T)_*$; since $A$ is normal, $T_*$ is a
submonoid of $A$. Every element of $A_\eta\backslash \{0\}$ can be
written uniquely as a product $mt$ with $m\in M$ and $t\in T$; since
$t\in A$, if $mt\in A$ then $m\in A\cap M$. Thus if we set
$B=A\cap M_*$ there is a decomposition $A\cong B\smsh T_*$.
In other words, $\MSpec(A)\cong\MSpec(B)\times\MSpec(T_*)$

Since every localization of $A$ has the form $A_\fp=B'\smsh T_*$,
the affine open subsets of $\MSpec(A)$ are are all of the form
$\MSpec(B')\times\MSpec(T_*)$. Gluing these together gives the
decomposition of $X$.
\end{proof}
\goodbreak

Note that the factorization in Proposition \ref{toric-by-finite} is not
unique; it depends upon the choice of isomorphism $A_\eta\cong(M\times T)_*$.

\begin{cor}\label{factor}
If $f:X\to X'$ is a morphism between separated and connected
normal monoid schemes of finite type,
inducing an isomorphism $f^*:\cA'_{\eta'}\to\cA_\eta$ of group completions,
then $f$ is isomorphic to the product of a morphism $X(\Delta)\to X(\Delta')$
of toric monoid schemes and an isomorphism $\MSpec(T_*)\to\MSpec(T'_*)$.
\end{cor}

\begin{proof}
By assumption, $f$ maps the generic point $\eta$ of $X$ to the
generic point $\eta'$ of $X'$.
Choosing a decomposition $\cA_\eta\cong(M\times T)_*$, we have an
implicitly defined decomposition $\cA'_{\eta'}\cong(M\times T)_*$.
Then for each $x\in X$ the decompositions $\cA_x\cong B_x\smsh T_*$,
$\cA'_{f(x)}\cong B'_{f(x)}\smsh T_*$ of Proposition \ref{toric-by-finite}
satisfy $f^*(B'_{f(x)})\subseteq B_x\subseteq M_*$.
Therefore the map $\cA'_{f(x)}\to\cA_x$ factors as a product of
$f^*(B'_{f(x)})\subseteq B_x$ and $T_*\cong T_*$, for each $x$.
The result follows.
\end{proof}

\begin{subrem}\label{not-product}
Not every morphism $(X,A)\times\MSpec(T_*)\to(X',A')\times\MSpec(T'_*)$
between connected normal monoid schemes
of finite type will factor as a cartesian product of maps
$(X,A)\to(X',A')$ and $\MSpec(T_*)\to\MSpec(T'_*)$.
For example, this fails for the canonical $\MSpec((\Z/n)_*)\to\MSpec(\Z_*)$.
However, such a map determines both a toric map $(X,A)\to(X',A')$
and a map $\MSpec(T_*)\to\MSpec(T'_*)$.
\end{subrem}
\goodbreak

\medskip\goodbreak
\section{Realizations of monoid schemes}\label{sec:realizn}

In this section we fix a commutative ring $k$. If $A$ is a monoid,
the ring $k[A]$ gives rise to a scheme $\Spec(k[A])$, which is
called the {\it $k$-realization} of $\MSpec(A)$. The affine spaces
$\A^n_k=\Spec(k[t_1,...,t_n])$ of \ref{affine} are useful examples.
The $k$-realization is a faithful functor from monoids to affine
$k$-schemes; a monoid morphism $A\to B$ naturally gives rise to a
morphism $\Spec(k[B])\to\Spec(k[A])$. 

If $X$ is an affine monoid scheme, we write $X_k$ for its realization:
\[
\MSpec(A)_k=\Spec(k[A]).
\]

\goodbreak

We saw in \eqref{pushout} that the $k$-realization functor commutes 
with pullback for affine monoid schemes, because it has a left adjoint
(defined on the category of affine $k$-schemes) sending $\Spec(R)$
to $\MSpec(R,\times)$, where $(R,\times)$ is the multiplicative
monoid whose underlying pointed set is $R$. Thus if $X=\MSpec(A)$ is
an affine monoid, the adjunction
$\Hom(\Spec(R),X_k)\cong\Hom_\MSch(\MSpec(R,\times),X)$
means that $X_k$ represents the functor 
sending $\Spec\,R$ to $\Hom_\MSch(\MSpec(R,\times),X)$.

\begin{defn}
Let $X$ be a monoid scheme and $k$ a ring. Define a contravariant
functor $F_X$ from the category of affine $k$-schemes to sets to be
the Zariski sheafification of the presheaf
\[
\Spec R \mapsto \Hom_{\MSch}(\MSpec(R, \times), X).
\]
If $X$ is affine, the presheaf is already a sheaf since it is
represented by $X_k$.
\end{defn}
\goodbreak

Recall from \cite[VI-14]{EisenbudHarris} that a contravariant
functor $F$ from affine $k$-schemes to sets is represented by a
unique $k$-scheme $X$ if and only if $F$ is a Zariski sheaf and $F$
admits a covering by open subfunctors $F_\alpha$, each of which is
represented by an affine scheme $U_\alpha$. If so, the representing
scheme $X$ is obtained by gluing the $U_\alpha$ together. Here, a
subfunctor $F_\alpha\subseteq F$ is {\em open} if for every
$k$-algebra $R$ and every morphism $\Hom(-,\Spec R) \to F$, i.e.,
for every element of $F(\Spec R)$, the pullback functor
$F_\alpha\times_F\Hom(-,\Spec R)$ is represented by an open
subscheme of $\Spec R$. A collection of subfunctors $\{F_\alpha\}$
of $F$ {\em covers} $F$ if for every $k$-algebra $L$ which is a
field, we have $F(\Spec L) = \bigcup_\alpha F_\alpha(\Spec L)$.

\begin{thm} \label{almostadjoint}
The functor $F_X$ is represented by a scheme $X_k$.
\end{thm}

\goodbreak

\begin{proof}
Suppose that $U=\MSpec(A)$ is any affine monoid subscheme of $X$.
Since sheafification preserves monomorphisms such as
$\Hom(-,U)\subseteq\Hom(-,X)$, $F_U$ is a subfunctor of $F_X$. If
$\Lambda$ is a local $k$-algebra and $L=\Spec(\Lambda)$ then
\addtocounter{equation}{-1}
\begin{subequations}
\begin{equation}\label{eq:FXlocal}
F_X(L)=\Hom_\MSch(\MSpec(\Lambda,\times),X).
\end{equation}
\end{subequations}
Since $\MSpec(\Lambda,\times)$ has a unique  point, each map
$\MSpec(\Lambda,\times)\to X$ factors through an affine open
submonoid $U$. Therefore $F_X$ is covered by the collection of
subfunctors $F_U$, as $U$ ranges over all affine open monoid
subschemes of $X$. We will show that the $F_U$ are open subfunctors
of $F_X$; we have seen that each $F_U$ is represented by the affine
scheme $U_k$. By \cite[VI-14]{EisenbudHarris}, this will prove that
$F_X$ is representable by the $k$-scheme which is obtained by gluing
the affine schemes $U_k$. 

Fix an affine open monoid subscheme $U=\MSpec(A)$. To prove that
$F_U$ is open, fix a $k$-algebra $R$ and consider a morphism
$\Hom(-,\Spec R)\to F_X$ and its corresponding element $\phi\in
F_X(\Spec R)$. We have to show that the pullback
$G=F_U\times_{F_X}\Hom(-,\Spec R)$ is represented by an open
subscheme $V$ of $\Spec(R)$. Since $F_X$ is a sheaf, $\Spec(R)$ has
an affine open covering $\{\Spec R[1/s] \,|\, s\in\cS\}$ such that
the restriction of $\phi$ to $F_X(\Spec R[1/s])$ is represented by a
morphism $\phi_s: \MSpec(R[1/s],\times) \to X$ of monoid schemes. By
Observation \ref{Spec-MSpec} below, there are continuous maps
\[
\Spec(R[1/s]) \into \MSpec(R[1/s],\times) \map{\phi_s} X.
\]

Let $V'_s$ denote the inverse image of $U$ under $\phi_s$ and let
$V_s$ denote the open subspace $V'_s\cap\Spec(R[1/s])$; we regard
$V_s$ as an open subscheme of $\Spec(R[1/s])$ and hence of
$\Spec(R)$. We claim that $G$ is represented by the open subscheme
$V=\cup V_s$
of $\Spec R$. To prove our claim, it suffices to consider a local
$k$-scheme $L=\Spec(\Lambda)$ and prove that $G(L) = \Hom(L, V)$ as
subsets of $\Hom(L,\Spec R)$. Since $L$ is local, we have
$F_U(L)=\Hom(A,(\Lambda,\times))$, and \eqref{eq:FXlocal} holds for
$X$.
Thus $G(L)$ is the set of all $f:L\to\Spec R$ such that
\[
\MSpec(\Lambda,\times) \map{f^\times} \MSpec(R,\times) \map{\phi} X
\]
maps the closed point $\fm$ of $L$ into $U$. If the image of $f$
lies in $V$, $\fm$ lands in some $V_s$ and hence $f^\times$ maps the
closed point $(\fm,\times)$ of $\MSpec(\Lambda,\times)$ into $V'_s$.
It follows that $\phi f^\times(\fm,\times)\in U$, i.e., $f\in G(L)$.
Thus $\Hom(L,V)\subseteq G(L)$.

Conversely, if $f:L\to \Spec(R)$ is in $G(L)$ then $f$ factors
through some $f_s:L\to \Spec(R[1/s])$
and $\phi_sf_s^\times$ maps the closed point $(\fm,\times)$ of
$\MSpec(\Lambda,\times)$ to a point in the subset $U$ of $X$, so
$f_s(\fm)\in V_s$. But since $L$ is local, this implies that
$f_s(L)\subseteq V_s$. The desired equality $G(L) = \Hom(L, V)$
follows.
\end{proof}
\goodbreak

\begin{obs}\label{Spec-MSpec}
Let $R$ be any commutative ring, and $(R,\times)$ its underlying
multiplicative monoid. If $\fp$ is a prime ideal of the ring $R$,
then $(\fp,\times)$ is a prime ideal of the monoid $(R,\times)$. The
resulting inclusion $\Spec(R)\into\MSpec(R,\times)$ is continuous
because if $s\in R$ the open subspace $D(s)$ of $\MSpec(R)$
intersects $\Spec(R)$ in the open subspace $\{\fp\subset
R\,|\,s\not\in\fp\}$.
If $R$ is local, the maximal ideal $\fm$ of $R$ maps to the maximal
prime $(\fm,\times)$ of $\MSpec(R,\times)$.
\end{obs}

\begin{defn}\label{defn:realizn}
Given a commutative ring $k$ and a scheme $(X,\cA)$, we define its
{\em $k$-realization} $X_k$ to be the scheme representing $F_X$.
\end{defn}

\begin{subrem}
Observe that $X_k = X_\Z \times_{\Spec \Z} \Spec k$ for any monoid
scheme $X$ and commutative ring $k$. Those preferring the notion of
a  field with one element ($\bF_1$) might prefer writing $X_k$ as $X
\times_{\Spec \bF_1} \Spec k$ or just $X \times_{\bF_1} k$.
\end{subrem}

\begin{cor}\label{cor:pullbackcommute}
The $k$-realization functor $X \mapsto X_k$ preserves arbitrary
limits (when they exist). In particular, it preserves pullbacks.
\end{cor}

\begin{proof}
Suppose that $\{X_i, i \in I\}$ is a diagram of monoid schemes and
that its limit $X$ exists in the category of monoid schemes. It
suffices to prove the canonical map
$$
F_X \to F = \varprojlim F_{X_i}
$$
is an isomorphism of sheaves on the category of affine $k$-schemes.
Recall that the limit of a diagram of sheaves exists and coincides
with the limit as presheaves. That is, we have $F(\Spec R) =
\varprojlim F_{X_i}(\Spec R)$. When $R$ is local, we have $F_X(\Spec
R)=\Hom(\MSpec(R,\times),X)$ and also
$$
F(\Spec R) = \varprojlim \Hom(\MSpec(R,\times),X_i) \cong
\Hom(\MSpec(R, \times), X),
$$
where the second isomorphism holds since $X = \varprojlim_i X_i$.
Since the sheaf map $F_X \to F$ is an isomorphism on all local
rings, it is an isomorphism of sheaves.
\end{proof}

In Proposition \ref{prop:realizn} below we shall give an explicit
construction of $X_k$ for separated $X$. We need some preliminaries.

\begin{lem}\label{lem:S-and-k}
If $S$ is multiplicatively closed in $A$, $S^{-1}k[A]\cong k[S^{-1}A]$.
\end{lem}

\begin{proof}
The monoid map $A\to S^{-1}A$ is initial among monoid maps $A\to B$ 
that take $S$ to units. Similarly, the map $k[A]\to S^{-1}k[A]$ is 
initial among $k$-algebra homomorphisms $k[A]\to C$ that take $S$ to
units. Being a left adjoint, the functor $k[-]$ preserves initial objects.
\end{proof}

\begin{rem}\label{rem:real-open}
Let $A$ be a monoid. Any affine open monoid subscheme of $\MSpec(A)$
has the form $\MSpec(A_\fp)$ for some prime ideal $\fp$ of $A$, by
Lemma \ref{lem:maxima}, and $A_\fp=A[1/s]$ by Lemma \ref{lem:yesopen}.
Hence $\Spec(k[A_\fp])\to \Spec(k[A])$ is an open immersion, by Lemma \ref{lem:S-and-k}.
\end{rem}

For the next Proposition, let us say that a point $x$ in a monoid
scheme $X$ is {\em nice} if the canonical map $U=\MSpec(\cA_x)\to X$
is an open immersion.
Every closed point is
nice by Lemma \ref{lem:maxima}, but the points of Example
\ref{ex:notopen} are not nice. If $X$ is of finite type, then every
point is nice by Lemma \ref{lem:poset}. The nice points $x\in X$ are
a cofinal subset of the poset underlying $X$ by Lemmas
\ref{lem:yesopen} and \ref{lem:maxima}, because the closed points in
any open subscheme are nice. If $x<y$ are two nice points then
$\Spec(k[\cA_x])\to\Spec(k[\cA_y])$ is an open immersion by Lemma
\ref{lem:S-and-k}.
The criterion for separatedness in Lemma \ref{lem:sep} uses nice points.

\begin{prop}\label{prop:realizn}
Let $k$ be a commutative ring and $(X,\cA)$ a separated monoid
scheme. Then the $k$-realization of $X$ is
\[
X_k = \varinjlim_{x \in X} \Spec(k[\cA_x]).
\]
\end{prop}

\begin{proof}
Put $U_x=\MSpec(\cA_x)$. Because nice points are cofinal in the
poset underlying $X$, the limit can be taken over the nice points.
If $x$ is nice, then $U_x\subset X$ is an open immersion;
set $V_x=(U_x)_k$. If $y$ is also nice, 
then $U_x\cap U_y$ is an affine open, because the
intersection of two affine open subschemes of a separated monoid
scheme is affine open by Corollary \ref{Cor524}. 
By Corollary \ref{cor:pullbackcommute}
we have $(U_x\cap U_y)_k=V_x\times_{X_k}V_y$. Let
$V_{x,y}$ be the image of the projection
$\pi_x:V_x\times_{X_k}V_y\to V_x$. Then $V_{x,y}$ is open in $V_x$ 
and we have an isomorphism $\psi_{x,y}=\pi_y(\pi_x)^{-1}:V_{x,y}\to
V_{y,x}$. Hence the family of schemes $V_x$ indexed by the nice
points of $X$ together with the open subschemes $V_{x,y}\subset V_x$
and the isomorphisms $\psi_{x,y}$ satisfy the hypothesis of
\cite[(4.1.7)]{EGA0} (or \cite[Ex.\,II.2.12]{Hart}). Therefore the
limit of the proposition exists, and is the scheme obtained by
gluing the realizations of the open affine subschemes of $X$. Since
this is also the definition of $X_k$, the proposition follows.
\end{proof}

The $k$-realization functor from monoid schemes to $k$-schemes is
faithful, because it is so locally: $\MSpec(A)_k=\Spec(k[A])$. (This
is clear if $X$ is separated, and follows from Theorem
\ref{almostadjoint} if it is not separated.) It is not full because
$k$-schemes such as $\A^1_k$ have many more endomorphisms than their
monoidal counterparts.

The realization functor loses information, because
distinct monoid schemes can have isomorphic realizations.
This is a well known phenomenon even for toric varieties, where the
additional data of a (faithful) torus action is needed to recover the fan.

\begin{ex}\label{ex:usualtoric}
For a fan $\Delta$ and any field $k$, the variety $X(\Delta)_k$ is the
 usual toric $k$-variety associated to $\Delta$. This is clear from
Construction \ref{ex:fan}.
\end{ex}

\begin{ex}\label{ex:groupring}
Let $T$ be a finite abelian group. The $k$-realization of
$\MSpec(T_*)$ is the cogroup scheme $\Spec(k[T])$. If $|T|$ is a
unit (or nonzerodivisor) in $k$ then $k[T]$ is reduced, but this
fails if $k$ is a field of characteristic $p>0$ and $T$ has
$p$-torsion.
\end{ex}

\goodbreak

\begin{lem}\label{Maschke}
Let $k$ be an integral domain and $A$ a cancellative monoid. Set
$X=\MSpec(A)$ and $U=\MSpec(A^+)$.
\begin{enumerate}
\item
If $A^+$ is torsionfree then $k[A]$ is a domain (i.e., $X_k$ is integral).

\item Suppose that $k$ is a normal domain containing a field; if $\chr(k)=p>0$
assume also that $A^+$ has no $p$-torsion. Then $k[A^+]$ is normal
and its subalgebra $k[A]$ is reduced. That is, $U_k$ is normal and
$X_k$ is reduced.

\item Suppose that $\chr(k)=p>0$ and $A^+$ has $p$-torsion.
Then $k[A]$ is not reduced; $k[A]_{\red}=k[B]$,
where the monoid $B$ is the quotient of $A$
by the congruence relation that $a_1\sim a_2$ if and only if
$a_1^{p^e} = a_2^{p^e}$ for some $e \geq 0$.
\end{enumerate}
\end{lem}

\begin{proof}
Since $A$ is the union of its finitely generated submonoids $A_i$,
and $k[A]=\cup k[A_i]$, we may assume that $A$ is finitely
generated. As noted before \ref{toric-by-finite}, we can write
$A^+=(M\times T)_*$ where $M$ is a free abelian group and $T$ is a
finite torsion group. Since $A$ is a submonoid of $A^+$, $k[A]$ is a
subalgebra of $k[A^+]$. If $T$ is trivial, $k[A]$ is a subring of
$k[M]$, which is manifestly a domain. If $k\supset\Q$ or if
$\chr(k)=p$ and $p\nmid|T|$ then $k\to k[T]$ is a finite \'etale
extension and $k[A]$ is a subring of $k[A^+]=k[T][M]$, which is
manifestly normal if $k$ is normal. Hence $k[A^+]$ and its
subalgebra $k[A]$ are reduced in this case.

Finally, suppose that $\chr(k)=p$ and that the $p$-torsion subgroup
$T_p$ of $T$ is non-trivial. Since $k[T_p]_{\red}=k$ and
$k[A^+/T_p]$ is reduced by (2), we have $k[A^+]_{\red}=k[A^+/T_p]$.
If $B$ is the image of $A\to A^+/T_p$ then $k[A]_{\red}$ is the
image $k[B]$ of $k[A]\to k[A^+/T_p]$. Two elements $a_1,a_2\in A$ go
to the same element of $A^+/T_p$ if and only if their quotient is
$p$-torsion, i.e., if and only if they are congruent under the
relation $\sim$ of the lemma. 
It follows that $B=A/\!\sim$\,; this concludes the proof.
\end{proof}

\begin{subrem}\label{rem:reduction}
If $(X,A)$ is a cancellative monoid scheme of finite type, and $k$
is of characteristic $p>0$, Lemma \ref{Maschke}(3) implies that
$(X_k)_{\red}$ is the $k$-realization of $(X,B)$, where $B=A/\!\sim\,$
is the quotient stalk functor of $A$ defined as in \ref{Maschke}(3).
\end{subrem}

\goodbreak

\begin{prop}\label{prop:real-closed}
If $(Y,\cB)\map{f}(X,\cA)$ is a closed immersion of monoid schemes
then $f_k:Y_k\to X_k$ is a closed immersion of schemes for all rings $k$.
\end{prop}

\begin{proof}
If $V\subseteq X$ is an affine open subscheme, then by Lemma
\ref{lem:maxima} there exists $x\in X$ such that $V=\MSpec(\cA_x)$.
We shall abuse notation and write $V\cap Y$ for $V\times_XY$.
If $V\cap Y=\emptyset$ then $V_k\cap
Y_k=(V\cap Y)_k= \emptyset$. Otherwise $V\cap Y = \MSpec(\cB_y)$ for
some $y$, and $\cA_x\to\cB_y$ is onto, by Definition
\ref{def:closedimm}. Since $k$-realization preserves pullbacks by
Corollary \ref{cor:pullbackcommute}, we have $f_k^{-1}(V_k) =
f^{-1}(V)_k = \Spec(k[\cB_y])$ and the restriction $f_k^{-1}(V_k)
\to V_k = \Spec k[\cA_x]$ of $f$ is induced by the surjection
$k[\cA_x] \to k[\cB_y]$. This proves that the restriction $Y_k\cap
V_k\to V_k$ of $f_k$
is a closed immersion. Since $V$ is an arbitrary affine open
subscheme of $X$, this proves that $Y_k\to X_k$ is a closed immersion.
\end{proof}

A partial converse of this proposition is true.

\begin{lem} \label{L33}
Suppose $i: Y \to X$ is a morphism of monoid schemes such that the
underlying map of topological spaces induces a homeomorphism onto its image.
For any ring $k$, if $i_k: Y_k \to X_k$ is a closed immersion,
then $i$ is a closed immersion of monoid schemes.
\end{lem}

\begin{proof}  It suffices to prove that if
$X =  \MSpec(A)$ is affine, then $Y$ is also affine and the associated
map of monoids is surjective. Let $\cB$ be the sheaf of monoids for
the scheme
$Y$ and set $B = \Gamma(Y, \cB)$.
The map $Y \to X$ factors as
$$
Y \to \MSpec B  \to \MSpec A.
$$
Upon taking $k$-realizations we have $Y_k = \Spec(R)$ and the map
induced by $Y_k \to X_k$ is a surjection:
$k[A] \onto R$.  Since this surjection factors through the map $k[A]
\to k[B]$, which is induced by a map of monoids $A \to B$,
we see that $k[B] \onto R$ is surjection as well.
Let $Y = \cup_j W_j$ be a covering by open
affine subschemes, with $W_j = \MSpec B_j$. Then the map $B \to
\prod_j B_j$ is injective and hence
so is the map $k[B] \to \prod_j k[B_j]$.
Since the latter map
factors as $k[B] \to R \to \prod_J k[B_j]$, it
follows that $k[B] \map{\cong} R$ is an
isomorphism. That is, the $k$-realization of
$$
Y \to \MSpec(B)
$$
is an isomorphism. Moreover, since $k[A] \onto k[B]$ is onto, so is
the map $A \onto B$, and hence
$\MSpec(B) \to X$ is a closed immersion. In particular, the map of
underlying topological spaces is
a homeomorphism onto its image. It follows from this (and our assumption)
that the map of topological spaces underlying
$Y \to \MSpec(B)$ is a homeomorphism onto its image.

We may thus assume that  the $k$-realization $Y_k \to X_k = \Spec(k[A])$
is an isomorphism. We next claim that $Y\to X$ is a surjection on points,
and hence (by our assumption that $Y$ is homeomorphic to its image)
a homeomorphism on underlying topological spaces. To see this, fix a
point $\fp \in X$ and consider the monoid map $i_{\fp}:A\to S^0=\{0,1\}$
sending $\fp$ to $0$ and $A \setminus \fp$ to $1$.
Let $Y'$ denote  the pullback of $Y\to X$ along the map
$\MSpec S^0 \map{i_{\fp}} X$.
By Corollary \ref{cor:pullbackcommute}, the map
$Y'_k \to \left(\MSpec S^0\right)_k = \Spec k$ is an isomorphism,
so in particular $Y'$ is non-empty. By Proposition
\ref{prop:pb}, it follows that $Y \to X$ is onto.

Since $X$ has a unique maximal point, so does $Y$. By
Lemma \ref{lem:maxima}, $Y$ is affine. Since $Y_k\cong\Spec(k[A])$
we conclude that $Y\cong X$.
\end{proof}

\begin{prop}\label{prop:real-sep}
For any ring $k$ and morphism of monoid schemes $f: Y \to X$, the map
$f$ is a separated morphism of monoid schemes
if and only if its $k$-realization $f_k: Y_k \to X_k$ is a separated
morphism of schemes.
\end{prop}

\begin{proof}
One direction is immediate from Corollary \ref{cor:pullbackcommute} and
Proposition \ref{prop:real-closed}.

Assume $f_k$ is separated. Since the underlying topological space of  $Y
\times_X Y$ is given by the pullback in the category of topological
spaces, it follows that $Y \map{\Delta} Y \times_X Y$ is a homeomorphism
onto its image. (Observe that $Y \to\Delta(Y)$ and
$\Delta(Y) \map{\pi_1} Y$ are continuous, and both compositions are
the identity, where
$\Delta(Y) \subset Y \times_X Y$ is given the subspace topology.)
Since $\Delta_k$ is a closed immersion, Lemma \ref{L33} applies to
finish the proof.
\end{proof}

\section{Normal and smooth monoid schemes}\label{sec:smooth}

Throughout this section, $k$ denotes an integrally closed domain
containing a field. The normalization $A_\nor$ of a cancellative
monoid $A$ is defined in Definition \ref{def:norm}; since
$(A_\fp)_\nor=(A_\nor)_{\fp_\nor}$, it makes sense to talk about the
normalization of any cancellative monoid scheme.

The $k$-realization of $X$ cannot be normal unless $X_k$ is reduced.
Lemma \ref{Maschke} shows that $k[A]$ is reduced unless $p>0$
and $A^+$ has $p$-torsion, in which case $k[A]_\red$ is $k[B]$,
where $B$ is a particular quotient of $A$, described there.

\goodbreak
\begin{prop} \label{P:normal}
Let $X = (X, A)$ be a cancellative monoid scheme of finite type
such that its $k$-realization $X_{k}$ is a reduced scheme.
Then
\begin{enumerate}
\item
the normalization of $X_k$ is the $k$-realization of $(X,A_\nor)$.

\item if $X$ is normal, connected and separated, there is a decomposition
\[ X_k = X'_k \times_k \Spec k[T] \]
where $X'_k$ is a toric $k$-variety and $k[T]$ is finite \etale over
$k$.
\end{enumerate}
\end{prop}

\goodbreak
\noindent
As in Remark \ref{not-product},
the decomposition in
Proposition \ref{P:normal}(2) is not natural in $X$.

\begin{proof}
Part (2) is immediate from Proposition \ref{toric-by-finite} and Corollary
\ref{cor:pullbackcommute}.

Since the normalization of a reduced scheme is the scheme
constructed by patching together the normalizations of an affine
cover, we may assume that $X$ is affine, i.e, $X=\MSpec(A)$. Since
$k[A_\nor]$ is integral over $k[A]$, we may assume that $A=A_\nor$.
In this situation, where $A$ is a normal monoid of finite type,
Proposition \ref{toric-by-finite} states that $A\cong A'\smsh T_*$
where $A'$ is torsionfree and $T$ is a finite abelian group. Since
$X_k$ is reduced, we know from Lemma \ref{Maschke}(3) and Example
\ref{ex:groupring} that $T$ has no $p$-torsion and $k[T]$ is finite
\etale over $k$.
Since $k[A]=k[T][A']$, we are reduced
to the case in which $A$ is normal and torsionfree, i.e.,
$X=\MSpec(A)$ is an affine toric monoid scheme. By Theorem
\ref{why_toric}, $X$ is associated to a fan $\Delta$; by Example
\ref{ex:usualtoric}, $X_k$ is the toric variety associated to
$\Delta$, and in particular $X_k$ is normal.
\end{proof}

\begin{subrem}
It is possible to give an elementary proof of this result using that
if $A$ is a torsionfree normal monoid then $k[A]$ is integrally
closed; see \cite[12.6]{Gilmer}.
\end{subrem}

\goodbreak
\paragraph{\em Finite morphisms.}
We will need to know that the normalization of a monoid scheme is
a finite morphism, at least when $X$ is of finite type.

We say that a morphism of monoid schemes $f:Y\to X$ is \emph{affine} if
$X$ can be covered by affine open subschemes $U_i=\MSpec(A_i)$
such that $f^{-1}(U_i)$ is affine. Equivalently, $f$ is affine if
$f^{-1}(U)$ is affine for every affine open subscheme $U\subset X$.

\begin{defn}\label{def:finite}
Let $f:Y\to X$ be a morphism of monoid schemes.
We say that $f$ is \emph{finite} if it is affine and
$\cA_X(U)\to \cA_Y(f^{-1}(U))$ is finite for every
affine subscheme $U\subset X$.
We say that $f$ is {\em integral} if
it is affine and $\cA_X(U)\to \cA_Y(f^{-1}(U))$
is integral for every affine subscheme $U\subset X$.
\end{defn}

If $X$ is cancellative, its normalization $X_\nor \to X$
is an integral morphism. To see this, we may assume $X = \MSpec(A)$ is
affine so that $X_\nor \to X$ is given by $A \into A_\nor$, where
the normalization $A_\nor$ is integral by Definition \ref{def:norm}.
We now show that if $X$ is also of finite type,
then $X_\nor\to X$ is finite.

\begin{prop}\label{finitenorm}
If $X$ is a cancellative monoid scheme of finite type, the
normalization $X_\nor\to X$ is a finite morphism.
\end{prop}

\begin{proof}
It suffices to show that if $A$ is a cancellative monoid
of finite type then $A\to A_\nor$ is finite. Since $A_\nor$ is
integral over $A$ it suffices by Lemma \ref{integralvfinite}(i)
to show that $A_\nor$ is of finite type. Because the group
completion $A^+$ is finitely generated, it has the form $(M\times T)_*$
where $T$ is a finite abelian group and $M$ is free abelian.
Since $A[T]=\bigcup At$ is finite over $A$, we may replace $A$ by $A[T]$
to assume that $T\subset A$. As in the proof of
Proposition \ref{toric-by-finite}, this implies that $A=B\smsh T_*$
where $B=A\cap M_*$ is a finitely generated submonoid of $M$. If
$\beta$ is the rational convex polyhedral cone of $M_{\R}$ spanned by
the generators of $B$, $B_\nor$ is $(\beta\cap M)_*$.
By Gordon's Lemma \cite{Fulton}, $B_\nor$ is finitely generated.
A fortiori, $A_\nor=B_\nor\smsh T_*$ is finitely generated.
\end{proof}

\goodbreak
\medskip
\paragraph{\em Smoothness.}
\goodbreak
\begin{defn}\label{p-smooth}

Let $p$ be a prime. A separated monoid scheme of finite type is {\em
$p$-smooth} if each stalk (equivalently, each maximal stalk) is the
smash product $S\smsh T_*$, where $S=G_*\smsh F$ is the smash
product of a free abelian group with a point adjoined and a free
abelian monoid, and $T$ is a finite abelian group having no
$p$-torsion. A separated monoid scheme is {\em $0$-smooth} if each
stalk has the form $S \smsh T_*$ with $T$ an arbitrary finite
abelian group.

We will say that $X$ is {\em smooth} if it is $p$-smooth for all $p$,
i.e., if each stalk is the product of a free group of finite rank and
a free monoid of finite rank.
\end{defn}
\goodbreak

A cone in a fan $(N,\Delta)$ is said to be {\it nonsingular} if it
is spanned by part of a $\Z$-basis for the lattice $N$, in which
case each monoid $\sigma^\vee\cap M$ is the product of a free
abelian group and a free abelian monoid. A fan is said to be
nonsingular if all its cones are nonsingular.

\begin{prop} \label{smooth-realizn}
Let $X =(X, A)$ be a separated cancellative monoid scheme of finite type.
Its $k$-realization $X_k$ is smooth over a field $k$ of characteristic
$p\ge0$ if and only if $X$ is $p$-smooth.
If $X$ is connected and $p$-smooth then, under the decomposition
\[
X = (X, A') \times \MSpec(T)
\]
of Proposition \ref{P:normal}, the fan underlying $(X, A')$
is nonsingular.
\end{prop}

\begin{proof}
Recall from \cite[2.1]{Fulton} that the toric variety associated to
a fan is smooth if and only if each of its cones is nonsingular.
Therefore the proposition is an immediate corollary of Proposition
\ref{P:normal} and Lemma \ref{Maschke}.
\end{proof}

\begin{subex}\label{smooth-not-psmooth}
The hypothesis in \ref{smooth-realizn}
that $X$ be cancellative is necessary. For example,
consider the monoid $A=\langle t,e \,|\, e=e^2=te \rangle$, which has
$k[A]\cong k[x]\times k$.
Thus $X=\MSpec(A)$ is not $p$-smooth but $X_k$ is smooth for every $k$.
\end{subex}

\medskip\goodbreak

\section{$\MProj$ and Blow-ups}\label{sec:blowup}

An $\N$-grading of a monoid $A$ is a pointed set decomposition
\[
A = \bigvee_{i=0}^\infty A_i
\]
such that $A_i \cdot A_j \subseteq A_{i+j}$; $\Z$-gradings are
defined similarly.  For each nonzero $a$ in $A$,
let $|a|$ denote the unique $i$ such that $a \in A_i$.
For every multiplicative set $S$,
the localization $S^{-1}A$ is $\Z$-graded by $|a/s|=|a|-|s|$.
For example, if $s\in A_i$ is non-zero we have
\[
{A\ad{s}}_0 = \left\{\frac{a}{s^n} \, | \, |a|=|s^n|=ni, n \geq 0\right\} \cup \{0\}.
\]
Let $A_{\ge1}$ denote the ideal
$\bigvee_{i\ge1} A_i = \{a \,|\, |a|>0\} \cup \{0\}$, so that
$A/A _{\ge1}\cong A_0$; the image of the corresponding map
$\MSpec(A_0)\to\MSpec(A)$
consists of the prime ideals of $A$ containing $A_{\ge1}$.

\begin{defn}\label{def:MProj}
If A is an $\N$-graded monoid, we define $\MProj(A) = (X,\cB)$ to be
the following monoid scheme. The underlying topological space is $X
= \MSpec(A)\setminus\MSpec(A_0)$ --- i.e., the open subspace of
those prime ideals of $A$ that do not contain $A_{\ge1}$. The stalks
of $\cB$ on $X$ are defined by sending $\fp \in
\MSpec(A)\setminus\MSpec(A_0)$ to $\cB_\fp=(A_\fp)_0$, the degree
zero part of $A_\fp$. If $\MSpec(A_\fp)\subset X$ is open, that is,
if $A_\fp=A[1/s]$ for some $s\in A_{\ge 1}$, then the map
$\MSpec(A_\fp)\to \MSpec(A_\fp)_0$ is a homeomorphism. Indeed this
follows from the fact that a prime ideal $\fq$ of $A[1/s]$
contains an element $a/s^n$ if and only if 
$\fq\cap (A[1/s])_0$ contains $a^{n|s|}/s^{n|a|}$.
Thus $\MProj(A)$ is
covered by the affine open subschemes $D_+(s)=\MSpec(A\ad{s}_0)$
where $s\in A_{\ge1}$, and moreover, every affine open subscheme is
of this form. Hence $\MProj(A)$ is a monoid scheme of finite type
whenever $A$ is a finitely generated monoid. The maps
$A_0\to(A_\fp)_0$ induce a structure morphism
$\MProj(A)\to\MSpec(A_0)$.
\end{defn}

\begin{subrem}\label{rem:realproj}
The $k$-realization of $A$ is the graded ring $k[A]$,
and $k[A\ad{s}_0]$ is the degree~0 part of the ring $k[A]\ad{s}$,
so the $k$-realization of $\MProj(A)$ is $\Proj(k[A])$.
\end{subrem}

\begin{obs}\label{naturalMProj}
The construction is natural in $A$ for maps $A\to A'$ of
graded monoids such that $A'=A\cdot A'_0$.
For such maps there is a canonical morphism
$\MProj(A')\to\MProj(A)$ induced by the restriction of
$\MSpec(A')\to\MSpec(A)$. If $s\in A_{\ge1}$, the affine open
$\MSpec(A'\ad{s}_0)$ maps to the affine open $\MSpec(A\ad{s}_0)$.

If $S\subset A_0$ is multiplicatively closed, $S^{-1}A$ is graded and
$\MProj(S^{-1}A)$ is $\MProj(A)\times_{\MSpec(A_0)}\MSpec(S^{-1}A_0)$.
It follows that this construction may be sheafified:
for any monoid scheme $(X,\cA_0)$ and any sheaf $\cA$ of graded monoids
on $X$ with $(\cA_x)_0=(\cA_0)_x$ for all $x\in X$, there is a monoid
scheme $\MProj(\cA)$ over $X$ whose stalk at each $x$ is $\MProj(\cA_x)$.
Moreover, if $f:(X',\cA'_0)\to(X,\cA_0)$ is a morphism of monoid schemes,
equipped with sheaves $\cA'$ and $\cA$ of graded monoids as above, any
graded extension $f^{-1}\cA\to \cA'$ of $f^{-1}\cA_0\to \cA'_0$
such that $\cA'=f^{-1}\!\cA\cdot \cA'_0$
induces a canonical morphism $\MProj(\cA')\to\MProj(\cA)$ over $f$.
\end{obs}

\begin{lem}\label{lem:closedproj}
If $f:A\to B$ is a surjective homomorphism of graded monoids, then
the induced map $\MProj(B)\to \MProj(A)$ is a closed immersion.
\end{lem}

\begin{proof}
As noted above, any affine open subscheme $U\subset \MProj(A)$ is of
the form $U=\MSpec(A\ad{s}_0)$ for some $s\in A_{\ge 1}$.
But $U\cap\MProj(B)=\MSpec(B\ad{f(s)}_0)$ is affine,
so we are in the case of Lemma \ref{lem:immersion1}.
\end{proof}

\paragraph{\em Projective monoid schemes.}

For a monoid $A$ and indeterminates $T_0, \dots, T_n$, let $A[T_0,
\dots, T_n]$ denote the monoid freely generated by $A$ and the
$T_i$. It is a graded monoid, where each element of $A$
has degree $0$ and each $T_i$ has degree $1$, and we
define $\bP_A^n$ to be $\MProj(A[T_0,\dots, T_n])$. More generally,
for any monoid scheme $X = (X, \cA)$, define $\bP^n_X$ to be $\MProj(\cB)$
where $\cB$ is the sheaf of graded monoids on $X$ defined by sending an
open subset $U$ to $\cA(U)[T_0,\dots,T_n]$. In other words,
$\bP^n_X$ is defined by patching together
the monoid schemes of the form $\bP^n_A$ as $\MSpec(A)$ ranges over
affine open subschemes of $X$.
If $X$ has finite type, so does $\bP^n_X$.

A morphism of monoid schemes $Y\to X$ is \emph{projective} if,
locally on $X$, it factors as a closed immersion $Y\to \bP^n_X$ for some $n$
followed by the projection $\bP^n_X\to X$.

\begin{lem}\label{Pn-separated}
Projective morphisms are separated.
\end{lem}

Although this follows from Proposition \ref{prop:real-sep},
we give an elementary proof here.

\begin{proof}
Since closed immersions are separated by Lemma \ref{lem:affsep}, it
suffices to show that the morphisms $\bP^n_X\to X$ are separated. We
may assume that $X=\MSpec(A)$, so that
$\bP^n_X=\MProj(A[T_0,\dots,T_n])$. By Definition \ref{def:MProj},
points of $\bP^n_X$ correspond to prime ideals in $A[T_0,\dots,T_n]$
not containing $\{T_0,\dots,T_n\}$. By Lemma \ref{primes-in-smash}
and Example \ref{affine}, every such prime ideal has the form
$P_{S,\fp}=A\smsh\langle S\rangle \cup \fp[T_0,\dots,T_n]$ where
$\fp$ is a prime ideal of $A$ and $\langle S\rangle$ is the prime
ideal generated by a proper subset $S$ of $\{T_0,\dots,T_n\}$;
moreover $\fp$ and $S$ are unique and the projection to $\MSpec(A)$
sends the point $P_{S,\fp}$ to $\fp$.
According to Lemma \ref{lem:sep}, it suffices to observe
that for every $P_{S,\fp}$ and $P_{S',\fp}$ the prime $P_{S\cap
S',\fp}$ is a unique lower bound. (The surjectivity condition of
Lemma \ref{lem:sep} is easy, and left to the reader.)
\end{proof}

\smallskip
\begin{subex}\label{ex:bu-projective}
If $B$ is a finitely generated graded monoid, then
$\MProj(B)\to\MSpec(B_0)$ is projective and hence separated
by Lemma \ref{Pn-separated}.
Indeed, this is a particular case of Lemma \ref{lem:closedproj},
since $B$ is a quotient of some $B_0[T_0,\dots,T_n]$.
\end{subex}

\goodbreak
\medskip
{\em Blow-ups.}
\smallskip

\noindent
Given a monoid $A$ and an ideal $I$, we consider the graded monoid
$A \vee I \vee I^2 \vee \cdots$, where $I^n$ has degree $n$.
It is useful to introduce a variable $t$, and rewrite this as
\[
A[It] = \bigvee\nolimits_{n\ge0}\ I^n t^n \subseteq A\smsh F_1.
\]
If $S$ is multiplicatively closed in $A$, then
$S^{-1}(A[It])\cong(S^{-1}A)[{S^{-1}It]}$. It follows that if
$\cI$ is a quasi-coherent sheaf of ideals in a monoid scheme $(X,\cA)$
then there is a monoid scheme $\MProj(\cA[\cI t])$ over $(X,\cA)$
obtained by patching the $\MProj(A[It])$ in the evident manner.

\begin{defn}\label{def:blowup}
If $X = (X,\cA)$ is a monoid scheme and $Z \subseteq X$ is a
equivariant closed subscheme, given by a quasi-coherent sheaf of
ideals $\cI$, we define the {\it blow-up} of $X$ along $Z$ to be the
monoid scheme $X_Z=\MProj(\cA[\cI t])$.
\end{defn}

\begin{subrem}\label{subrem:blow}
If $X=\MSpec(A)$ is affine and $Z = \MSpec(A/I)$ then
$X_Z=\MProj(A[It])$, together with the structure morphism
$\MProj(A[It])\to\MSpec(A)$. Since $\MProj(A[t])\cong\MSpec(A)$, it
follows that for $U=X\setminus Z$ we have $X_Z \times_X U \cong U$.

\goodbreak

\smallskip

\noindent
The blow-up construction is natural in the pair $(A,I)$ in the
following sense. If $A\to B$ is a morphism of monoids, $I$ is an
ideal of $A$ and $J=IB$, there is a canonical graded morphism
$A[It]\to B[Jt]$ satisfying the hypotheses of \ref{naturalMProj}.
Hence there is a morphism $\MProj(B[Jt])\to\MProj(A[It])$ of the
blow-ups over $\MSpec(B)\to\MSpec(A)$. More generally, if $f:X'\to
X$ is a morphism of monoid schemes, $\cI$ is a quasi-coherent sheaf
of ideals on $X$ and  $\cJ=f^{-1}\cI\cdot\cA'$, then the morphism
$f^{-1}\cA[\cI t]\to \cA'[\cJ t]$ induces a canonical morphism
$\MProj(\cA'[\cJ t])\to\MProj(\cA[\cI t])$ over $f$, described in
\ref{naturalMProj}.
\end{subrem}

\begin{subrem}\label{blowupsareprojective}
The blow-up of $X$ along a quasi-coherent sheaf of ideals $\cI$ is
projective provided $\cI$ is given locally on $X$ by
finitely generated ideals, by \ref{ex:bu-projective}.
For example, if $X$ has finite type then the blowup of $X$ along
any quasi-coherent sheaf of ideals is projective.
\end{subrem}

\begin{ex}\label{ex:toricblowup}
Suppose $N$ is a free abelian group with basis $\{ v_1,\dots,v_n\}$,
and $\{ x_1,\dots,x_n\}$ is the dual basis of $M$. Let $\sigma$ be
the cone in $N_\R$ generated by $\{ v_1,\dots,v_d\}$, the
corresponding affine monoid scheme is $X(\sigma)=\MSpec(A)$, where
$A$ is generated by $x_1,\dots,x_n$ and $x_{d+1}^{-1},\dots,x_n^{-1}$
subject to $x_ix_i^{-1}=1$ for $d<i\le n$. The blow-up of
$X(\sigma)$ along the ideal generated by $x_1,\dots,x_d$ is the
toric monoid scheme $X(\Delta)$, where $\Delta$ is the subdivision
of the fan $\{\sigma\}$ given by insertion of the ray spanned by
$v_0=v_1+...+v_d$. To see this, it suffices to copy the
corresponding argument for toric varieties given in
\cite[p.\,41]{Fulton}.
\end{ex}

\begin{subex}\label{blowup-natch}
If $Z$ is an equivariant closed subscheme of $X$, defined by a
quasi-coherent sheaf of ideals $\cI$,
and $f:X'\to X$ is a morphism, then by naturality of the blow-up construction,
discussed above, there is a canonical morphism over $f$, from the blow-up 
$X'_{Z'}$ of $X'$ along the pullback $Z'=Z\times_X X'$ to the blow-up $X_Z$.
\end{subex}

\begin{lem} \label{L0106}
Let $f:X'\to X$ be a finite morphism of monoid schemes (\ref{def:finite}).
Let $Z$ be an equivariant closed subscheme of $X$, $X_Z$ the
blow-up along $Z$, and $X'_{Z'}$ the blow-up of $X'$ along
the pullback $Z'=Z\times_X X'$.
Then $\tilde{f}:X'_{Z'} \to X_Z$ is a finite morphism.
\end{lem}

\goodbreak
\begin{proof}
We may assume that $X$, and hence $X'$, is affine. Then $f$ is induced
by a map $A\to A'$, $Z$ is defined by an ideal $I\subset A$ and $Z'$
is defined by $J=I\, A'$.
Moreover because $f$ is assumed finite, there are elements
$c_1,\dots,c_r\in B$ such that $B=\bigcup_iAc_i$. If
$a_0,\dots,a_n$ generate $I$, and $b_0,\dots,b_n$ are their images in
$B$, then $\tilde{f}$ restricts to maps $D_+(b_i)\to D_+(a_i)$
induced by the monoid maps
$A_i=A[a_0/a_i,\dots,a_n/a_i]\to B_i=B[b_0/b_i,\dots,b_n/b_i]$.
By inspection, $B_i=\cup_{j=1}^nA_ic_j$.
\end{proof}

\begin{prop}\label{same-blowup}
Let $Z$ be an equivariant closed subscheme of a monoid scheme $X$ of
finite type. Then for any commutative ring $k$ the blow-up of $X_k$
along $Z_k$ is canonically isomorphic to the $k$-realization of the
blow-up of $X$ along $Z$.
\end{prop}

\begin{proof}
It suffices to consider the case
$X=\MSpec(A)$, $Z=\MSpec(A/I)$. In this case $S=k[A[It]]$ is the usual
Rees ring $k[A][Jt]$, $J=k[I]$. Since the blowing-up of $X_k=\Spec(k[A])$
along $Z_k=\Spec(k[A/I])$ is $\Proj(S)$, we have the desired identification
$\Proj(S)=\Proj(k[A[It]])=\MProj(A[It])_k$.
\end{proof}

\smallskip

\goodbreak
We conclude this section by observing that blow-ups of monoid
schemes satisfy a universal property analogous to that for
blow-ups of usual schemes. To state it, we need some
notation. We define a \emph{principal invertible} ideal of $A$ to be
an ideal $I$ such that there is an $x\in I$ such that the map
$A\map{x}I$ ($a\mapsto ax$) is a bijection. If $I$ is a principal
invertible ideal of $A$ then the canonical map 
$\MProj(A[It])\to\MSpec(A)$ is an isomorphism.

A quasi-coherent sheaf of ideals of a monoid scheme $X$ is said to be
\emph{invertible} if $X$ can be covered by affine open subschemes $U$
such that $\cI(U)$ is a principal invertible ideal of $\cA_X(U)$.
If $(X,\cA)$ is a monoid scheme and $\cI\subset\cA$ is a quasi-coherent
sheaf of ideals (see Definition \ref{qc_ideal}),
we say that a morphism $f:Y\to X$
\emph{inverts} $\cI$ if $f^{-1}\cI\cdot\cB$ is an invertible sheaf on $Y$.

\begin{prop}\label{prop:uniblow}
Let $X$ be a monoid scheme of finite type, $Z$ an equivariant closed
subscheme defined by a quasi-coherent sheaf of ideals $\cI$,
and $\pi:\widetilde{X}\to X$ the blow-up of $X$ along $Z$.
Then $\pi$ inverts $\cI$ and is universal with this property
in the sense that if $Y$ is of finite type and $f:Y\to X$ inverts $\cI$,
then the dotted arrow in the diagram below exists and is unique.
\[
\xymatrix{Y\ar@{.>}[r]\ar[dr]_f&\widetilde{X}\ar[d]^\pi\\ &X}
\]
\end{prop}

\begin{proof}
We may assume that $X=\MSpec(A)$ for some finitely generated monoid $A$,
that $\cI$ corresponds to an ideal $I$ of $A$, and that
$\widetilde{X}=\MProj(A[It])$.
The map $\pi$ inverts $I$ because the restriction of $\pi^{-1}\cI$
 to $D_+(s)$ is generated by $s$ for each $s\in I$.
Let $\cB$ be the structure sheaf of $Y$,
and write $\cJ$ for the sheaf of ideals $f^{-1}\cI\cdot\cB$.
By Example \ref{blowup-natch}, there is a unique morphism from the blow-up
$\widetilde{Y}=\MProj(\cB[\cJ t])$ to $\widetilde{X}$ over $f$.
By assumption, $\cJ$ is an invertible sheaf. Hence
$\widetilde{Y}\to Y$ is an isomorphism, because locally
$\cJ$ is a principal invertible ideal $J$ of $B$
and $\MProj(B[Jt])\cong\MSpec(B)$.
\end{proof}

\medskip
\goodbreak

\section{Proper morphisms}\label{sec:proper}
A monoid $V$ is called a {\em valuation monoid} if $V$ is
cancellative and for every non-zero element $\alpha \in V^+$, at least one of
$\alpha$ or $\frac{1}{\alpha}$ belongs to $V$.  For example, if $R$
is a valuation ring, then the underlying multiplicative monoid $(R,
\times)$ is a valuation monoid. Also, the free pointed monoid on one
generator is a valuation monoid. Given a valuation monoid $V$, the
monoid $V \smsh M_*$ is also a valuation monoid for any abelian
group $M$. For example, the monoid $\ip{y_1^{\pm 1}, \dots, y_n^{\pm
1}, x}$ is a valuation monoid.

Given a valuation monoid $V$, the units $U(V)$ are a subgroup of
$V^+ \setminus 0$ and the quotient group $(V^+\setminus0)/U(V)$
is a totally ordered abelian group with the total ordering defined by
$x \geq y$ if and only if $\frac{x}{y}$ belongs to the image of
$V \setminus 0$. To conform to usual custom, we convert the group law
for $(V^+\setminus0)/U(V)$ into $+$. We also adjoin a base point,
written $\infty$, to obtain the
totally ordered pointed (additive) monoid

\[
\Gamma := \left((V^+\setminus0)/U(V)\right)_*.
\]

\noindent
We extend the total ordering to $\Gamma$ by declaring that
$\gamma\le\infty$ for all $\gamma\in\Gamma$.
We call $\Gamma$ the {\em value monoid} of the valuation monoid
$V$. The canonical surjection
\begin{equation} \label{Eordmap}
\ord: V^+\onto \Gamma
\end{equation}
is called the {\em valuation map} of $V$.
The monoid $V$ is then
identified with the set of $x \in V^+$ such that $\ord(x) \ge 0$
(where, recall $0$ is the identity of $\Gamma$),
and the maximal ideal $\fm$ of $V$ is $\{ x\,|\, \ord(x) > 0\}$ (since
$\ord(x)=0$ just in case $x$ is a unit of $V$).

Note that \eqref{Eordmap} satisfies
$\ord(x) \le \infty$, $\ord(xy) = \ord(x) + \ord(y)$ and $\ord(x) = \infty$
if and only if $x = 0$.
Conversely, given an abelian group $M$ and a surjective
morphism $\ord: M_* \to\Gamma$ onto a totally ordered monoid
$(\Gamma,+,0, \infty)$ that satisfies these conditions,
the set $C = \{a \in M \, | \, \ord(a) \ge 0\}$ is a valuation monoid
whose pointed group completion is $M_*$ and whose associated valuation map is $\ord$.

\begin{lem}\label{Vnormal}
A valuation monoid $V$ has no finite extensions contained in $V^+$.
\end{lem}

\begin{proof}
Suppose that $V\subseteq B\subseteq V^+$ with $B$ finite over $V$.
By Lemma \ref{integralvfinite}(ii), $B$ is integral over $V$.
For every nonzero $b\in B$ there is an $n\ge1$ so that $b^n\in V$ and hence
$n\, \ord(b)\ge0$, which implies that $\ord(b)\ge0$ and thus $b\in V$.
\end{proof}

\begin{ex}
A {\em discrete valuation monoid} is a valuation monoid whose value
monoid is isomorphic to $\Z \cup \{\infty\}$ with its canonical ordering.
In this case, a lifting of the generator $1 \in \Z$ to an element $\pi$
in the discrete  valuation monoid $V$ is a
generator of the maximal ideal of $V$ and every non-zero element of
$V^+$ may written uniquely as $u\pi^n$ for $n \in \Z$ and $u \in
U(V)$. Let's call such an element a {\em uniformizing parameter}.

Observe that if $R$ is a discrete valuation ring, then $(R, \times)$ is
a discrete valuation monoid and the notion of a uniformizing parameter
has its usual meaning.

If $V$ is a discrete valuation monoid, its valuation map induces
a surjection $\pi: V^+\setminus 0 \onto \Z$; write $M=\ker\pi$. 
A choice of uniformizing parameter $t$ is equivalent to
a section of $\pi$ and identifies $V^+\setminus 0=M\times\ip{t^{\pm
1}}$. Under this identification, $\pi$ is the evident projection.
Thus, every discrete valuation monoid $V$ is isomorphic to $U(V)_*
\smsh \ip{t}$, where $\ip{t}$ is the free abelian monoid on one
generator and $U(V)$ is the group of units of $V$. Any element of
the form $u \smsh t$ with $u \in U(M)$ is a uniformizing parameter.
\end{ex}

\begin{subrem}\label{rem:zs}
It is well-known that a valuation ring is Noetherian if and only if it
is a discrete valuation ring; see \cite[\S VI.10, Thm.\,16]{ZS} for a
proof. The same argument shows that a valuation monoid is finitely
generated if and only if it is a discrete valuation monoid with a
finitely generated group of units.
\end{subrem}

\begin{defn}\label{def:monoidproper}
A map $f: Y \to X$ of monoid schemes satisfies the {\em valuative
  criterion for properness} if for every
valuation monoid $V$ and commutative square
\addtocounter{equation}{-1}
\begin{subequations}
\begin{equation} \label{E914}
\xymatrix{
\MSpec(V^+) \ar[r] \ar[d] & Y \ar[d]^f \\
\MSpec(V) \ar[r] \ar@{..>}[ur] & X \\
}\end{equation}
\end{subequations}
there is a unique map $\MSpec(V) \to Y$ causing both triangles to commute.

We say $f$ satisfies the {\em valuative criterion of
  separatedness} if each such square has at most one completion.

A map $Y \to X$ of monoid schemes of finite type is said to be {\em
  proper} if it satisfies the valuative
  criterion for properness.
\end{defn}

\begin{subrem} We are not certain what the correct definition of
``proper'' is for monoid schemes not of finite type. (Recall from
Remark \ref{notproper} that ``separated and universally closed'' is
clearly not the correct definition.)
\end{subrem}

Given any  morphism $f:\MSpec(V)\to X$, any affine open $U\subset X$
containing $f(\fm)$ (where $\fm$ is the unique closed point of
$\MSpec(V)$)  will contain the image of $\MSpec(V)$.
Hence the valuative criterion of properness and separatedness are
local on the base: if $Y|_U\to U$ satisfies one of these critera
for every $U$ in a covering of $X$, then so does $Y\to X$.

It is immediate from Definition \ref{def:monoidproper} that the class of
maps satisfying the valuative criterion of properness (resp., separatedness) is closed under
composition and pullback.

\begin{prop}\label{prop:integralproper}
A finite morphism between monoid schemes satisfies the valuative
criterion of properness.
\end{prop}

\begin{proof}
Suppose $Y \to X $ is finite and consider a commutative square \eqref{E914}
with $V$ a valuation monoid. We may assume $Y \to X$ is a map of affine
schemes, say given by a map of monoids
$A \to B$. Then the square \eqref{E914} is associated to the square
$$\xymatrix{
V^+ & B \ar[l] \\
V \ar[u] & A. \ar[l] \ar[u] \\
}$$
of monoids. The image of $B$ in $V^+$ is finite over $V$, but $V$ is
closed under finite extensions in $V^+$, by Lemma \ref{Vnormal}.
It follows that the map $B \to V^+$
actually lands in $V$, which gives the diagonal map we seek.
\end{proof}

\begin{cor}\label{cor:integralproper}
Closed immersions satisfy the
valuative criterion of properness.
\end{cor}

\begin{construction} \label{valconstruction}
To prove Theorem \ref{valuativeThm} below, we need a technical construction:
Let $V$ be a valuation monoid with group completion $V^+$ and
value monoid $(\Gamma,+)$. Recall that totally ordered groups are
necessarily torsion-free, and hence,
for any field $k$, the ring $k[\Gamma]$ is an integral
domain by Lemma \ref{Maschke}.

For an element
$$
\alpha = \sum_\gamma a_\gamma \gamma
$$
of $k[\Gamma]$ (where for this ring we have
rewritten $\Gamma$ using $\cdot$
instead of $+$ notation), define
$$
\ord(\alpha) =
min\{\gamma \in \Gamma \, | \, a_\gamma \ne
0\}.
$$
(For $\alpha = 0$, set $\ord(0) = \infty$.) It is easily verified that
$
\ord: (k[\Gamma], \times) \to \Gamma
$
is a monoid map such that $\ord(\alpha) = \infty$ if and only if
$\alpha = 0$.

It follows that we get an induced map of pointed group completions
$$
\ord: (k(\Gamma), \times) \to \Gamma
$$
where $k(\Gamma)$ denotes the field of fractions of $k[\Gamma]$. Moreover,
the composition
$$
V^+ \to (k(\Gamma), \times) \map{\ord} \Gamma
$$
coincides with the original valuation map $\ord: V^+ \to \Gamma$.

Finally, the pair $(k(\Gamma), \ord)$ is a valuation in the usual
ring-theoretic sense. To prove this,
it remains to show
$$
\ord(\alpha + \beta) \geq
\min\{\ord(\alpha), \ord(\beta)\} \, \text{ for all $\alpha, \beta \in k(\Gamma)$.}
$$
One easily reduces to the case when $\alpha, \beta \in k[\Gamma]$,
where it is obvious from the definition of $\ord$.
\end{construction}

\begin{prop} \label{Prop33}
Given a valuation monoid $V$, with pointed group completion $V^+$ and
value monoid $\Gamma$,  let $\ord$ be the valuation map on the
field $k(\Gamma)$ given in Construction \ref{valconstruction}, and let
$R
\subset k(\Gamma)$ denote the associated valuation
ring. Then the square of affine monoid schemes
$$
\xymatrix{
\MSpec(k(\Gamma), \times) \ar[r] \ar[d] & \MSpec(V^+)  \ar[d] \\
\MSpec(R, \times) \ar[r] & \MSpec(V)
}
$$
is a pushout square in the category of monoid schemes.
\end{prop}

\begin{proof}
For any monoid scheme $T$, suppose morphisms  $f: \MSpec(V^+) \to T$ and
$g: \MSpec(R, \times) \to T$ are given causing the evident square to commute. Let $t \in T$ be the image of
the unique closed point of $\MSpec(R, \times)$ under $g$, and let $U
\subset T$ be any affine open subscheme of $T$ containing $t$. Then $g$
factors through $U$.
Since $\MSpec(k(\Gamma), \times) \to
\MSpec(V^+)$ is a bijection on underlying sets (each is a one-point
set), the unique point of $\MSpec(V^+)$ also lands in $U$ and hence
$f$ too factors  through
$U$. We may thus assume $T = U$ is affine. That is, it suffices to prove
$$
\xymatrix{
V \ar[r] \ar[d] & (R, \times) \ar[d] \\
V^+ \ar[r] & (k(\Gamma), \times)
}
$$
is a pullback square in the category of pointed monoids.  But this is
evident since
$$
\xymatrix{
V^+ \ar[r] \ar[d]^{\ord} & (k(\Gamma), \times) \ar[d]^{\ord} \\
\Gamma \ar[r]^{=} & \Gamma \\
}
$$
commutes, $V = \{\alpha \in V^+ \, | \, \ord(\alpha) \geq 0\}$ and
$R = \{\beta \in k(\Gamma) \, | \, \ord(\beta) \geq 0\}$.
\end{proof}

Recall that a map of (classical) $k$-schemes $Y_k \to X_k$, where $k$
is a field, is said to satisfy the {\em valuative criterion of
  properness} (resp., {\em separatedness}) if
every solid arrow square
$$
\xymatrix{
\Spec(F) \ar[r] \ar[d] & Y_k \ar[d] \\
\Spec(R) \ar@{..>}[ur] \ar[r] & X_k \\
}
$$
has a unique (resp., at most one) completion making both triangles commute, whenever $R$ is a
valuation ring (which is necessarily a $k$-algebra) and $F$ is its field of factions.

\begin{thm} \label{valuativeThm}
Let $f: X \to Y$ be a morphism of monoid schemes and let $k$ be a
field. The morphism $f_k: X_k \to Y_k$ satisfies the valuative
criterion of properness (resp., separatedness)  if and only if $f$
satisfies the valuative criterion of properness (resp., separatedness).
\end{thm}

\begin{proof}
By Theorem \ref{almostadjoint},
for any local $k$-algebra $R$, there is a natural adjunction isomorphism
$$
\Hom_k(\Spec(R),X_k) \cong F_X(\Spec R)=\Hom_{\MSch}(\MSpec(R,\times),X).
$$
Now suppose $R$ is a valuation ring with field of fractions $F$. Then
$V = (R, \times)$ is a valuation monoid with $V^+ = (F, \times)$.
Since $R$ and $F$ are local, a commutative square of the form
$$\xymatrix{
\Spec F \ar[r] \ar[d] & Y_k \ar[d] \\
\Spec R \ar[r] & X_k \\
}$$
corresponds via adjunction to a commutative square
of monoid schemes given by the solid arrows in the diagram
\begin{equation} \label{E3-3}
\xymatrix{
\MSpec(V^+) \ar[r] \ar[d] & Y \ar[d] \\
\MSpec(V) \ar@{..>}[ur] \ar[r] & X. \\
}
\end{equation}
If $Y \to X$ satisfies the valuative criterion of properness
(resp. separatedness), there exists a unique (resp., at most one) morphism of monoid schemes
$\MSpec(V) \to Y$ represented by the dotted arrow above that  causes both triangles to commute.
Again by adjunction, this gives a unique map $\Spec(R)\to Y_k$ causing both
triangles to commute in the first square.

Conversely, say a square \eqref{E3-3} is given. By Construction
\ref{valconstruction}, there is a valuation ring $R$ with field of fractions
$F = k(\Gamma)$ and morphisms $\MSpec(R, \times) \to \MSpec V$ and
$\MSpec(F, \times) \to \MSpec V^+$ fitting into a commutative diagram
\begin{equation} \label{E3-3b}
\xymatrix{
\MSpec(F, \times) \ar[r] \ar[d] & \MSpec(V^+) \ar[r] \ar[d] & Y \ar[d] \\
\MSpec(R, \times) \ar@{..>}[rru] \ar[r] & \MSpec(V) \ar[r] & X. \\
}
\end{equation}
By Proposition \ref{Prop33}, the left-hand square
is a pushout square in the category of monoid
schemes. Using adjunction as above, if $Y_k \to X_k$ satisfies the
valuative criterion of properness (resp., separatedness), there exists
a unique (resp., at most one) map
represented by the dotted arrow in \eqref{E3-3b}  that causes the outer
two triangles to commute. Since the left-hand square is a pushout,
it follows immediately that there exists a unique (resp., at most one) arrow $\MSpec(V) \to Y$
causing both triangles in \eqref{E3-3b} to commute.
\end{proof}

\begin{cor} \label{cor:properft}
For any field $k$, a morphism between monoid schemes of finite type
$Y \to X$ is proper if and only if $Y_k \to  X_k$ is proper.
\end{cor}

\begin{proof} Merely observe that $Y_k$ and $X_k$ are Noetherian,
and apply the valuative criterion of the properness theorem
\cite[II.4.7]{Hart}.
\end{proof}

\begin{subrem}
Say $f: Y \to X$ satisfies the valuative criterion of properness.
If $Y_k$ is quasi-compact,  EGA\,II(7.2.1) implies that $f_k$ is proper. 
\end{subrem}

\begin{cor}\label{Cor:DVM}
A morphism between monoid schemes of finite type is proper
if and only if it satisfies the
valuative criterion of properness
of Definition
\ref{def:monoidproper} for all discrete valuation monoids.
\end{cor}

\begin{proof}
If $f: X \to Y$ satisfies the criterion of Definition
\ref{def:monoidproper} for all discrete valuation monoids, then, for
any field $k$, its $k$-realization $f_k: X_k \to Y_k$ satisfies the
valuative criterion of properness for all DVRs. This follows, using
adjunction,  from the fact that $\MSpec(R,\times)$ is a discrete
valuation monoid if $R$ is a DVR. Since $X_k$ and $Y_k$ are Noetherian
and $f_k$ has finite type, it follows that $f_k$ is proper
(see \cite[Ex.\,II.4.11]{Hart}).
The result now follows from Corollary \ref{cor:properft}.
\end{proof}

\begin{cor}\label{Pn-proper}
A projective morphism $Y\to X$ between monoid schemes of finite type
is proper. In particular,
if $X$ is a monoid scheme of finite type and $X_Z$ is the blow-up along
an equivariant closed subscheme $Z$, then
the map $X_Z\to X$ is proper.
\end{cor}

\begin{proof} Using Proposition \ref{prop:real-closed} and Remark
\ref{rem:realproj}, we see that if $k$ is a field, then $Y_k \to X_k$
is a projective morphism of $k$-schemes and hence is proper. For the
second assertion, recall that $X_Z \to X$ is projective and $X_Z$
has finite type.
\end{proof}

\begin{rem} 
In fact, a projective morphism of arbitrary monoid schemes satisfies 
the valuative criterion of properness. We sketch the proof of this fact.
First one observes that, by Corollary \ref{cor:integralproper}, 
it suffices to check that for any monoid scheme $X$ and $n\ge 1$ 
the projection $\bP^n_X\to X$ satisfies the criterion. Second, 
one reduces further to showing that if $V$ is a valuation monoid then 
any section $\MSpec (V^+)\to \bP^n_{V^+}$ of the canonical projection 
extends to a section $\MSpec(V)\to \bP^n_V$ of $\bP^n_V\to \MSpec (V)$.
Third, one observes that for an affine scheme $\MSpec A$ a section of
$\bP^n_A\to\MSpec A$ is determined by an equivalence class of
$n$-tuples $(b_0,\dots,b_n)$ of elements of $A$ such that at least
one of the $b_i$ is nonzero, modulo the coordinate-wise action of
$U(A)$. Finally, one proves that if $b=(b_0,\dots,b_n)$ determines a
section $\MSpec (V^+)\to \bP^n_{V^+}$ as above, then multiplying the
$b_i$ by an appropriate power of a uniformizing parameter we obtain
an equivalent tuple $b'$ with $b'_i\in V$ for all $i$. Thus the
section extends to $\MSpec(V)$.
\end{rem}

Recall from \ref{def:toric} that a monoid scheme of finite type
is toric if it is separated, connected, torsionfree and normal.
By Theorem \ref{why_toric}, there is a faithful functor from fans
to toric monoid schemes.

\begin{cor} \label{cor:properonfans}
Let $\phi: (N',\Delta') \to (N, \Delta)$ be a morphism of fans. Then
the associated morphism of toric monoid schemes $X' \to X$ is proper
if and only if $\phi$ has the property that for each $\sigma
\in\Delta$, $\phi_\R^{-1}(\sigma)$ is a union of cones in $\Delta'$.
\end{cor}

\begin{proof} 
This follows from the well-known fact that if $k$ is a field, then $X'_k \to
  X_k$ is proper if and only if  $\phi$ has the stated property
(see \cite[p.\,39]{Fulton}).
\end{proof}

\begin{cor}\label{properisseparated}
Every proper map between monoid schemes of finite type is separated.
\end{cor}

\begin{proof}
By Theorem \ref{valuativeThm} and the Valuative Criterion of
Separatedness Theorem for Noetherian schemes, 
then the $k$-realization of a proper map between monoid schemes of
finite type is separated if $k$ is a field. Now use Proposition
\ref{prop:real-sep}.
\end{proof}

\section{Partially cancellative torsion free monoid schemes}\label{sec:pctf}

A monoid $A$ is {\em pctf } if it is isomorphic to a monoid of the
form $B/I$ where $B$ is a cancellative torsion free monoid (i.e., a
cancellative monoid whose group completion is torsion free) and $I$ is
an ideal.  A monoid scheme is {\em pctf} if all of its stalks are.

\goodbreak
\begin{prop} \label{pctfproperties}
We have:
\begin{enumerate}

\item If a pctf monoid is finitely generated, then it is isomorphic to
$A/I$ where $A$ is  a finitely generated torsion free cancellative monoid.

\item All submonoids and localizations of a pctf monoid are pctf.
In particular, for a monoid $A$, $\MSpec(A)$ is pctf if and only if
$A$ is pctf.

\item If $A$ is a pctf monoid and $\fp$ is a
prime ideal, then $A/\fp$ is a cancellative torsionfree monoid.

\item An open subscheme of a pctf monoid scheme is pctf.

\item An equivariant closed subscheme of a pctf monoid scheme is pctf.
\end{enumerate}
\end{prop}

\goodbreak

\begin{proof}
Say $A = B/I$ with $B$ cancellative and torsion free. Pick elements
$b_1, \dots, b_m$ in $B$ that map to a generating set of $A$ and let
$B'$ be the submonoid of $B$ they generate. Then $A = B'/(I \cap B')$,
proving the first assertion.

For the second, say $A = C/I$ with $C$ cancellative and torsion free.
If $B$ is a submonoid of $A$, let $B'$ denote the inverse image of
$B$ in $C$ and set $I' = I \cap B'$. Then $B = B'/I'$, and so $B$ is
pctf. The assertion concerning localizations holds since
$S^{-1}(C/I) \cong S^{-1}C/S^{-1}I$. The remaining assertion of part
(2) is clear.

If $A=B/I$ then $A/\fp=B/\fp'$ for some prime ideal of $B$, so
(3) follows from the elementary observation that if $A$ is cancellative
and torsionfree then so is $A/\fp$.

Assertion (4) is local and follows from (2); hence (5) is local,
and is then easy.
\end{proof}

\begin{prop}\label{prop:pctfblows}
The blow-up of a pctf monoid scheme along an equivariant closed subscheme
is pctf.
\end{prop}

\begin{proof}
Let $Y \to X$ be the blow-up of a pctf monoid scheme $X$
along an equivariant closed subscheme. Since the question is local on $X$,
we may assume that $X$ is affine, say $X = \MSpec(A)$ with $A$ pctf.
Then $Y$ is
$\MProj(A[It])$ for an ideal $I$. For each $s\in I$, we get an
affine open subset of $Y$ given by the monoid
$\{\frac{f}{s^n} \,|\, f \in I^n,\ n\ge0 \}$. This is a submonoid of
$A[\frac{1}{s}]$ and hence is pctf. The collection of such open subsets
as $s$ varies over all elements of $I$ form an open cover of $Y$. Thus
$Y$ is pctf.
\end{proof}

\begin{prop} \label{prop:closure}
Let $X=(X,\cA)$ and $Y=(Y, \cB)$ be monoid schemes
and let $f:Y\to X$ be a morphism. There is a unique closed subscheme
$Z$ of $X$ which is minimal with respect to the property that
$f$ factors through $Z \subset X$.

If $U \subset X$ is an affine open subscheme of $X$, then $Z\cap U$
is the affine scheme $\MSpec(C)$, where the monoid 
$C$ is the image of $\cB(U) \to \cA(U\times_X Y)$.
In particular, if $X$ is of finite type then so is $Z$.
\end{prop}

\begin{proof}
If $f$ factors through two different closed subschemes $W_1$
and $W_2$ of $X$, then it factors through $W_1 \times_X W_2$, which is
(canonically isomorphic to) a closed subscheme of $X$ (see Example
\ref{intersection}). So, we define $Z$ to be the
inverse limit taken over the partially ordered set of closed subschemes $W$
of $X$ such that $f$ factors through $W$.

For the local description of $Z$, we may assume that $X = U =\MSpec(B)$
is affine. Any closed subscheme of $X$ has the form $W=\MSpec(D)$ with
$B \to D$ a surjection of monoids. Then $f$ factors through $W$
if and only if $B\to \cA(Y)$ factors through $D$,
that is, if and only if $B \to C$ factors as $B\to D\to C$; in other
words, if and only if $Z\subseteq W$.
\end{proof}

\begin{defn}\label{defn:closure}
The subscheme $Z$ of \ref{prop:closure} is called the 
{\em scheme-theoretic image} of $f.$ 
If $f$ is an open immersion, we write $\overline{Y}$ for $Z$ and
(by abuse) call it the {\em closure} of $Y$\!.
\end{defn}

\begin{prop} \label{prop:closureispctf}
Let $Y$ be a monoid scheme and suppose $U \subset Y$ is an open
subscheme that is pctf. Then the scheme-theoretic image
$\overline{U}$ of $U$ in $Y$ is pctf. Moreover, if $Y$ is separated,
then $\overline{U}$ is separated.
\end{prop}

\begin{proof}
The first assertion is local on $Y$ and so we may assume
$Y = \MSpec(B)$ for a monoid $B$ and $U = \MSpec(S^{-1}B)$ for a
multiplicative subset $S$.  Then $\overline{U}$ is the affine scheme
associated to the image $\overline{B}$ of $B \to S^{-1}B$. The monoid
$S^{-1}B$ is pctf by assumption and
\ref{pctfproperties}(4), and hence so is $\overline{B}$
by \ref{pctfproperties}(2).

The second assertion is just the observation that
a closed subscheme of a separated scheme
is also separated by Lemma \ref{lem:affsep}.
\end{proof}

\medskip\goodbreak
\section{Birational morphisms}\label{sec:birat}

A morphism $p:Y\to X$ of monoid schemes is {\em birational} if
there is an open dense subscheme $U$ of $X$ such that $p^{-1}(U)$
is dense in $Y$ and $p$ induces an isomorphism from $p^{-1}(U)$ to $U$.

\begin{prop}[Birational maps]
Let $p: (Y,B) \to (X,A)$ be a map between monoid schemes of finite type,
Then $p$ is birational if and only if the following conditions hold:
\begin{enumerate}
\item $p$ maps the generic points of $Y$ bijectively onto the
generic points of $X$
\item A point $y\in Y$ is generic if (and only if) $p(y)\in X$ is generic
\item for each generic point $y\in Y$ the induced map
    $A(p(y))\to B(y$) on stalks is an isomorphism.
\end{enumerate}
\end{prop}

\begin{proof}
If $p$ is birational and $U$ is as in the definition above, then
$U$ contains all of the generic points of $X$ and
$p^{-1}(U)$ contains all the generic points of $Y$
as well as every point of $y$ that maps to a generic point of $X$.
The conditions are then clearly satisfied.

Conversely,
take $U$ to be the (dense open) set of generic points
of $X$. By hypothesis, $p^{-1}(U)$ is the set of generic points of $Y$
and the map $p: p^{-1}(U) \to U$ is bijective. Hence $p^{-1}(U)$ is
open and dense. Since the map $p^{-1}(U) \to U$ is bijective and
induces an isomorphism on all stalks, it is an isomorphism.
\end{proof}

\begin{cor} \label{properbirational}
If $p: X' \to X$ is a proper map of toric monoid schemes that is
birational, then $p$ is given by a map of fans $\phi: (N', \Delta')
\to (N, \Delta)$ such that $\phi: N'\map{\cong} N$ and the image of
$\Delta'$ under the isomorphism $\phi_\R$ is a subdivision of
$\Delta$. Conversely, any such map $\phi$ induces a proper birational map of
monoid schemes.
\end{cor}

\begin{proof}
From \ref{why_toric}(2), $p$ comes from a morphism of fans such that
$\phi: N'\map{\cong} N$, and such a morphism is a subdivision by
Corollary \ref{cor:properonfans}. Conversely, if $p$ is induced by a
morphism of fans $\phi:(N', \Delta') \to (N, \Delta)$ such that
$\phi_\R$ is a subdivision of $\Delta$, then $p_k$ is proper by
\cite[\S2.4]{Fulton}; hence $p$ is proper by Theorem
\ref{valuativeThm}.
\end{proof}

\goodbreak
\begin{ex}\label{ex:birat-components}
If $X$ is a monoid scheme of finite type, let $X_\eta$ denote the
equivariant closure of a generic point $\eta$ (in the sense of
\ref{equivclosure}).
Then each $X_\eta$ has a unique generic point,
namely $\eta$.
If $X$ is pctf, then each $X_\eta$ is
cancellative and torsionfree by \ref{pctfproperties}(3), and hence pctf.
If $X$ is reduced, the morphism $\coprod_\eta X_\eta\to X$
is birational.
\end{ex}

\begin{prop} \label{pb-birational}
If $Y \to X$ is a birational map
and $X' \to X$ is a morphism such that $X'$ is of finite type and
every generic point of $X'$ maps to a generic point of $X$,
then the pullback $Y \times_X X' \to X'$ is birational.
\end{prop}

\begin{proof} 
The poset underlying $Y \times_X X'$ is given by the pullback
of the underlying posets (by \ref{prop:pb}).
Since $Y \to X$ is birational, a point $(y,x')$
in $Y \times_X X'$ is generic if and only if $x'$ is a generic point
of $X'$, and in this case $y$ and $x'$ map to the same point $x$ of $X$,
which is generic. Hence the map
$Y\times_X X'\to X'$ is a bijection on sets of generic points.
Writing $A',A$ and $B$ for the stalk functors of $X',X$ and $Y$\!,
the map on generic stalks is of the form
$A'(x') \to A'(x') \smsh_{A(x)} B(y)$. This is an isomorphism,
since the map $A(x)\to B(y)$ is an isomorphism.
\end{proof}

Define the {\em height} of a point $x$ in a monoid scheme $X$
to be the dimension of $\cA_x$; i.e., it is
the largest integer $n$  such that there exists a strictly
decreasing chain $x = x_n > \cdots > x_0$ in the poset underlying $X$.
We write this as $ht(x)$ or $ht_X(x)$.

For example,
if $X=X(N,\Delta)$ is the monoid scheme associated to a fan, then
$ht(\sigma)=\dim(\sigma)$ for each cone $\sigma \in \Delta$.
Here $\dim(\sigma)$ refers to the dimension of the
real vector subspace of $N_\R$ spanned by $\sigma$.

\begin{lem}\label{lem:heightgoesup}
Suppose $p: Y \to X$ is a proper, birational map of separated
pctf schemes of finite type.
Then for any $y \in Y$, we have $ht_Y(y) \leq ht_X(p(y))$.
\end{lem}

\begin{proof}
Suppose $ht_Y(y) = m$, so that we have a chain of points
$y = y_m > \cdots > y_0$ in $Y$. Clearly $y_0$ must be minimal,
and thus generic. Let $\eta= p(y_0)$, and define $X_\eta$ to be the
equivariant closure of $\{\eta\}$ in $X$. As pointed out in
Example \ref{ex:birat-components}, $X_\eta$ is cancellative and
torsionfree. The pullback $Y_\eta = X_\eta\times_X Y$ is an equivariant
closed subscheme of $Y$ containing $y_0$ as its unique generic point, and
hence each $y_i$.  By Proposition \ref{pctfproperties}(5),
$Y_\eta$ is also pctf, and $Y_\eta\to X_\eta$ is birational by
Proposition \ref{pb-birational}.

Let $Y'$ denote the equivariant closure of $y_0$ in $Y_\eta$.
By Example \ref{ex:birat-components}, $Y'\to Y_\eta$ is birational,
$Y'$ contains all the $y_i$ and $Y'$ is cancellative and torsionfree.
Replacing $X$ and $Y$ by $X_\eta$ and $Y'$, we may assume that both
$X$ and $Y$ are connected, cancellative and torsionfree.
Hence the normalization maps $X_\nor \to X$ and $Y_\nor\to Y$ exist and are
homeomorphisms (by \ref{normzn-homeo}), and both $X_\nor$ and $Y_\nor$ are
torsionfree. Since $Y\to X$ is birational, it induces a
birational morphism $Y_\nor\to X_\nor$.
The map $Y_\nor\to Y$ is finite by \ref{integralvfinite} and
hence proper by \ref{prop:integralproper}. Thus $Y_\nor\to X$ and
hence $Y_\nor\to X_\nor$ are proper. Thus we may assume that
$X$ and $Y$ are separated, normal and torsionfree.

By Proposition \ref{toric-by-finite} and Corollary \ref{factor},
we have reduced to the case where $Y\to X$ is a proper birational map
of toric monoid schemes, given by a map of fans
$\phi: (N',\Delta') \to (N,\Delta)$.
The birational hypothesis means that $\phi: N' \to N$ is an isomorphism.
By Corollary \ref{properbirational}, the proper hypothesis means that
$\Delta'$ is a subdivision of $\Delta$.
Since $\phi(\sigma)$ is the smallest cone in $\Delta$ containing the
image of $\sigma$ under $\phi_\R$ and since height corresponds to
dimension of cones, the result is now clear.
\end{proof}

\medskip
\goodbreak
\section{Resolutions of singularities for toric varieties}\label{sec:ROS}

The purpose of this section is to establish some properties for
monoid schemes that are analogous to those known to hold for
arbitrary varieties in characteristic zero. These properties will be
used in Section \ref{sec:cdh} to prove that certain presheaves of
spectra satisfy the analogue of ``smooth $cdh$ descent'' for monoid
schemes.

\goodbreak
\begin{thm}\label{thm:resolution1}
Let $X$ be a separated cancellative pctf monoid scheme of finite type.
Then there is a birational proper morphism $Y \to X$
such that $Y$ is smooth.
\end{thm}

\begin{proof}
We may assume that $X$ is connected. Since the normalization map
is proper birational by Propositions \ref{finitenorm} and
 \ref{prop:integralproper}, we may assume that $X$ is normal.
Since $X$ is pctf it is torsionfree by Proposition \ref{pctfproperties}(3).
By Proposition \ref{toric-by-finite},
$X$ is toric and $X \cong X(\Delta)$ for some fan $\Delta$.
There exists a subdivision $\Delta'$ of $\Delta$ such that
$X(\Delta')$ is smooth, and it follows from Corollary \ref{properbirational} that the morphism $X(\Delta')\to X(\Delta)$
is proper birational.
\end{proof}
\goodbreak\smallskip

Let $N$ be a free abelian group of finite rank. Recall
(from \cite[page ~34]{Fulton}, e.g.) that a
cone in $N_\R$ is called \emph{simplicial} if it is generated by
linearly independent vectors, and that a fan is simplicial if every
cone in it is simplicial. We will need the notion of the barycentric
subdivision of a simplicial fan $\Delta$ in $N_\R$: For a simplicial
cone $\sigma$ in $N_\R$ of dimension $d$, let $v_1, \dots,v_d$ be
the minimal lattice points along the one-dimensional faces of $\sigma$, 
also called the rays of $\sigma$. For each non-empty
subset $S$ of $\{1, \dots, d\}$, let $v_S = \sum_{i\in S} v_i$. The
{\em barycentric subdivision} of $\sigma$, which we write as
$\sigma^{(1)}$, is defined as the collection of $2^d$ cones given as
the span of vectors of the form $v_{S_1}, \dots v_{S_e}$, where $0
\leq e \leq d$ and $S_1 \subset \cdots \subset S_e$ is a chain of
proper subsets of $\{1,\dots,d\}$. It is clear that if $\tau$ is a
face of $\sigma$, then the set of cones in $\sigma^{(1)}$ that are
contained in $\tau$ form the fan $\tau^{(1)}$. It follows that
$$
\Delta^{(1)} := \left\{ \sigma^{(1)} \, | \, \sigma \in \Delta \right\}
$$
is again a simplicial fan. We inductively define
$\Delta^{(i)} = (\Delta^{(i-1)})^{(1)}$ for $i \geq 2$.

\begin{lem} \label{L:subdivide}
If $\Delta'$ is any subdivision of a simplicial fan $\Delta$ in $N_\R$,
then for $i\gg0$, the fan $\Delta^{(i)}$ is a subdivision of $\Delta'$.
\end{lem}

\begin{proof}
It suffices to show that any ray of $\Delta'$, that is, any
$1$-dimensional cone of $\Delta'$, is a ray of some
$\Delta^{(i)}$. Given a positive integer combination
$v=\sum n_iv_i$ of the vertices in a cone, we may reorder the
vertices to assume the $n_i$ are in decreasing order. Then $v$ is in
the cone of $\Delta^{(1)}$ spanned by the $v_{S_i}$, where
$S_i=\{1,\dots,i\}$, and (if $v\ne v_1$) we can write $v=\sum n'_i
v_{S_i}$ with $\sum n'_i<\sum n_i$. The result follows by induction
on $\sum n_i$.
\end{proof}

\begin{lem} \label{L510b}
If $\Delta$ is a smooth fan, then for all $i \geq 1$, the toric
monoid scheme
$X(\Delta^{(i)})$ is obtained from $X(\Delta)$
via a sequence of blow-ups along smooth centers.
\end{lem}

\begin{proof}
We may assume $i=1$. If $\Delta$ is smooth already, then
$\Delta^{(1)}$ is also smooth. In
general, the fan $\Delta^{(1)}$ is obtained from $\Delta$ via a
series of steps of the following sort: starting with a smooth fan
$\Delta$, we form a subdivision $\Delta'$ by picking a cone
$\sigma$, letting $v_1,\dots,v_d$ be the minimal lattice points
along its rays, and defining $\Delta'$ to be the subdivision of
$\Delta$ given by insertion of the ray spanned by $v_1+\cdots+v_d$.
By Example \ref{ex:toricblowup}, $X(\Delta') \to X(\Delta)$ is the
blow-up along the smooth, closed equivariant subscheme defined by
$x_1=\cdots=x_d=0$.
\end{proof}

\begin{thm}\label{thm:resolution2}
For a morphism $\pi: Y \to X$ between separated cancellative pctf
monoid schemes of  finite type,
assume $X$ is smooth and $\pi: Y \to X$ is proper and  birational.
Then there exists a sequence of blow-ups
along smooth closed equivariant centers,
$$
X^n \to \cdots \to X^1 \to X_0 = X,
$$
such that $X^n \to X$ factors through $\pi:Y\to X$.
\end{thm}
\goodbreak

\begin{proof}
By Theorem \ref{thm:resolution1},
there is a proper birational morphism $Z \to Y$ with $Z$ smooth.
We may therefore assume
that $Y$ is smooth. We may also assume that $X$ and $Y$
are connected, so that they have unique generic points.

Thus, by Corollary \ref{properbirational},
$Y\to X$ is given by a morphism $(N^\prime, \Delta^\prime) \to (N,\Delta)$ of fans
that is an isomorphism of lattices and such that
$\Delta'$ is a subdivision of $\Delta$.
Lemmas \ref{L:subdivide}  and \ref{L510b} complete the proof.
\end{proof}

\medskip\goodbreak

\section{cd structures on monoid schemes.}\label{sec:cdh}

Let $\cMpctf$ denote the category of monoid schemes of finite type
that are separated and pctf.
In this section, we will be concerned with cartesian squares
of the form
\begin{equation} \label{sq:cd}
\xymatrix{
D \ar[r] \ar[d] & Y \ar[d]^p \\
C \ar[r]^e & X.
}
\end{equation}

\begin{defn} \label{def:abstractblowup}
An {\em abstract blow-up} is a cartesian square
of monoid schemes of finite type of the form \eqref{sq:cd}
such that $p$ is proper, $e$ is an equivariant closed immersion, and
$p$ maps the open complement $Y \setminus D$ isomorphically onto
$X \setminus C$.
The square with $Y=\emptyset$ and $C=X_\red$ is such a square.
\end{defn}

\begin{prop} \label{blowupisabstract}
If $X$ is of finite type, $C$ is an equivariant closed subscheme of $X$
and $p:Y \to X$ is the blow-up of $X$ along $C$, then the resulting cartesian
square is an abstract blow-up.
If $X$ belongs to $\cMpctf$, so do $Y$, $C$ and $D$.
\end{prop}

\begin{proof}
By Corollary \ref{Pn-proper}, $p$ is proper.
As noted in Definition \ref{def:blowup},
$p$ maps $Y \setminus D$ isomorphically to $X\setminus C$
(because $D = C \times_X Y$).
The second assertion follows from Propositions \ref{pctfproperties}
and \ref{prop:pctfblows}.
\end{proof}

\begin{prop} \label{pctf-abstblowup}
Suppose an abstract blow-up square \eqref{sq:cd} is given with $X$
in $\cMpctf$. Let $\Bar{Y}$ be the scheme-theoretic image of
$Y\setminus D$ in $Y$, and define $\Bar{D} = C \times_X \Bar{Y}$.
Then
$$\xymatrix{
\overline{D} \ar[r] \ar[d] & \overline{Y} \ar[d]^p \\
C \ar[r]^e & X
}$$
is an abstract blow-up square in $\cMpctf$.
\end{prop}

\begin{proof}
By Proposition \ref{pctfproperties}, $X\setminus C$ and
hence $Y\setminus D$ is pctf, and so by
Proposition \ref{prop:closureispctf}, $\overline{Y}$ is pctf as well.
Since equivariant closed subschemes of pctf schemes are
pctf, $C$ and $\overline{D}$ also belong to $\cMpctf$.
The map $\overline{Y} \to X$ is a composition of proper maps and
hence is proper. Finally, $\Bar{Y}\setminus \Bar{D} = Y \setminus D$.
\end{proof}

Recall from \cite[2.1]{VVcdh} that a $cd$ structure on a category $\cC$
is a collection of distinguished commutative squares in $\cC$.
If $\cC$ has an initial object $\emptyset$, any $cd$ structure defines
a topology: the smallest Grothendieck topology such that
for each distinguished square \eqref{sq:cd} the sieve generated by $\{ p,e\}$
is a covering sieve (and the empty sieve is a covering of the initial object).
The coverings $\{ p,e\}$ are called {\em elementary}.

\begin{defn} \label{def:cdstructures}
The {\em blow-up $cd$ structure} on $\cMpctf$ is given by
the collection of all abstract blow-up squares with $X,Y,C,D$ all
belonging to $\cMpctf$. The {\em Zariski $cd$ structure} on $\cMpctf$
is given by all cartesian squares associated to a covering
of $X$ by two open subschemes.

The $cdh$ topology on $\cMpctf$ is the topology generated by the union
of these two $cd$ structures.
\end{defn}

Following \cite[2.3, 2.4]{VVcdh},
we say that a $cd$ structure is {\it complete} if
$\cC$ has an initial object $\emptyset$
and any pullback of an elementary covering contains a
sieve which can be obtained by iterating elementary coverings.
We say that a $cd$ structure is {\em regular} (see \cite[2.10]{VVcdh})
if each distinguished
square \eqref{sq:cd} is a pullback, $e$ is a monomorphism and the
morphism of sheaves
\begin{equation}\label{eq:reg}
\rho(D) \times_{\rho(C)} \rho(D) \amalg \rho(Y) \to \rho(Y)
\times_{\rho(X)} \rho(Y)
\end{equation}
is onto, where $\rho(T)$ denotes the sheafification of the presheaf
represented by $T$.

\begin{thm} \label{thm:completeregular}
The blow-up and Zariski $cd$ structures on $\cMpctf$ are
complete and regular.
\end{thm}

\begin{proof}
The completeness property for Zariski squares is clear since they are
preserved by pullback, and the regularity property is even clearer.
For the blow-up $cd$ structure, consider an abstract blow-up square
$$
\xymatrix{
D \ar[r] \ar[d] & Y \ar[d]^p \\
C \ar[r]^e & X.
}
$$
Let $X' \to X$ be any morphism in $\cMpctf$ and consider the square
involving $X', C', Y'$ and $D'$ formed by pullback. The scheme $Y'$
might not belong to $\cMpctf$, but the scheme-theoretic image $Y''$
of $Y' \setminus D'$ in $Y'$ does by Proposition
\ref{pctf-abstblowup}. The resulting square involving $C', X', Y''$
and $D'' := C' \times_{X'} Y''$ is an abstract blow-up by the same
result, and hence by \cite[Lemma 2.4]{VVcdh} the blow-up $cd$
structure is complete.

For the regularity property, we need to show that
\eqref{eq:reg} is onto.
 Every object admits a covering in
this topology by affine, cancellative monoids, and it suffices to
prove surjectivity of the map given by the underlying presheaves
evaluated at such an affine cancellative  $U$. That is, say $f: U \to Y$,
$g:U \to Y$ are given with $p \circ f = p \circ g$. We need to prove
either $f=g$ or they both factor through $D$ and coincide as maps to $C$.
Let $u$ be the
unique generic point of $U$. If either $f(u)$ or $g(u)$ lands in $Y
\setminus D$, then they both must land there. Since $Y \setminus D
\cong X \setminus C$, it follows that $f$ and $g$ coincide
generically. But since $U$ is cancellative, it follows $f=g$ on all of
$U$. (To see this, one may work locally: If $h,l: A \to B$ are two maps of monoids with
$B$ cancellative and the compositions of $h,l$ with the inclusion $B
\into B^+$ coincide, then $h=l$.)
Otherwise, we have that the generic point, and hence every point, of
$U$ is mapped by both $f$ and $g$ to points in the closed subset $D$ of
$Y$. Again using that $U$ is cancellative, it follows that $f,g$
factor through $D \into Y$. (This is also proven by working locally.)
Finally, the compositions of these maps $f,g : U \to D$ with $D \to C$
coincide since $C \to X$ is a closed immersion.
\end{proof}

We define the {\em standard density structure} on $\cMpctf$ as follows:
The set $\cD_i(X)$ consists of those open immersions $U\subset X$
such that every point in $X \setminus U$ has height at least $i$.
It is clear that this satisfies the axioms required of a density
structure of finite dimension in \cite[2.20]{VVcdh}.

A $cd$ structure is said to be {\it bounded} for a given density
structure if any distinguished square has a refinement which is
reducing for the density structure in the sense of \cite[2.21]{VVcdh}.

\begin{thm} \label{thm:bounded}
The blow-up and Zariski $cd$ structures on $\cMpctf$ are both
bounded for the standard density structure.
\end{thm}

\begin{proof}
To see that the blow-up $cd$ structure is bounded, we need to show
that any abstract blow-up square \eqref{sq:cd} in $\cMpctf$ has a
refinement that is reducing for $\cD_*$. Consider the square
obtained by replacing $Y$ by the monoid scheme-theoretic image of $Y
\setminus D$ (in the sense of Definition \ref{defn:closure}),
and $D$ by the pullback. This is also an
abstract blow-up square, and it refines \eqref{sq:cd}. This
refinement has the features that $p^{-1}(X \setminus C)$ is dense in
$Y$, $Y$ maps birationally onto the scheme-theoretic image of
$X\setminus C$ in $X$, and $D$ does not contain any generic points
of $Y$.

To show that this square is reducing, we assume given
$C_0 \in \cD_i(C), Y_0 \in \cD_i(Y)$ and $D_0 \in \cD_{i-1}(D)$.
Define $X'$ to be the open subscheme $X \setminus Z$ of $X$, where
$Z \subset X$ is the equivariant closure (in the sense of
\ref{equivclosure}) of the union of the images of each of
$C \setminus C_0$, $D\setminus D_0$ and $Y\setminus Y_0$ in $X$.
We need to show that $X'$ belongs to $\cD_i(X)$ and that the pullback
of the original square \eqref{sq:cd} along $X' \into X$
gives an abstract blow-up square.

If $y \in Y$ is a point of height at least $i$, then $p(y)$ has
height at least $i$ in the scheme-theoretic image of $X \setminus
C$, by Lemma \ref{lem:heightgoesup}. Hence $p(y)$ has height at
least $i$ in $X$ itself (since a closed immersion is an injection on
underlying posets). If $d \in D$ has height at least $i-1$, then its
height in $Y$ is at least $i$ (since $D$ contains no generic points
of $Y$) and hence its image in $X$ has height at least $i$ too.
Since $C$ is an equivariant closed subscheme, if $c\in C$ has height
at least $i$, it has height at least $i$ in $X$.

Recall that $Z \subset X$ is the equivariant closure of the union of the images
of each of $C \setminus C_0$, $D \setminus D_0$ and $Y \setminus Y_0$ in $X$.
Each of these images consists of points of height at least $i$
and hence every point in $Z$ has height at least $i$ in $X$
by Remark \ref{rem:equiclose}.
Therefore $X'$ belongs to $\cD_i(X)$
and the pullback of the above square along $X' \into X$ gives an
abstract blow-up square that proves our original square is reducing.

The argument in the previous paragraphs applies {\em mutatis mutandis}
to show that every Zariski square is reducing.
\end{proof}

\begin{cor}\label{cor:bounded}
Let $\mathscr{S}$ be a presheaf of abelian groups on
$\cMpctf$; let $t$ be either the Zariski or the $cdh$-topology, and
write $a_t\mathscr{S}$ for the sheafification with respect to $t$.
If $X\in\cMpctf$ is of dimension $d$, then
\[
H^n_t(X,a_t\mathscr{S})=0 \mbox{ for } n>d.
\]
\end{cor}
\begin{proof}
Immediate from Theorem \ref{thm:bounded} and \cite[Thm.
2.26]{VVcdh}.
\end{proof}

\smallskip\goodbreak

The category of spectra we use in this paper will not be critical.
In order to minimize technical issues, we will use the terminology
that a {\it spectrum} $E$ is a sequence $E_n$ of simplicial sets
together with bonding maps $b_n: E_n \to \Omega E_{n+1}$. We say
that $E$ is an {\it $\Omega$-spectrum} if all bonding maps are weak
equivalences. A map of spectra is a strict map. We will use the
model structure on the category of spectra defined in \cite{BF}.
Note that in this model structure, every fibrant spectrum is an
$\Omega$-spectrum. Given a Grothendieck topology, the category of
contravariant functors $\cF$ from $\cMpctf$ to spectra ({\it
presheaves} of spectra) has a closed model structure, in which a
morphism $\phi:\cF\to\cF'$ is a cofibration when $\cF(X)\to\cF'(X)$
is a cofibration for every monoid scheme $X$ in $\cMpctf$; $\phi$ is
a weak equivalence if it induces isomorphisms between the sheaves of
stable homotopy groups (see \cite{JardineSPS}, \cite{JardineGen}).
We write $\bH_{\cdh}(-, \cF)$ for the fibrant replacement of $\cF$
using this model structure for the $cdh$ topology, as in \cite{chsw}.

A presheaf of spectra $\cF$ on $\cMpctf$ satisfies the
{\em Mayer-Vietoris property} for some family $\cC$ of cartesian squares
if $\cF(\emptyset)=\ast$ and the application of $\cF$ to each member of the
family gives a homotopy cartesian square of spectra.

\begin{prop} \label{cdh-MV}
Let $\cF$ be a presheaf of spectra on $\cMpctf$. Then the canonical
map $\cF(X)\to\bH_{\cdh}(X,\cF)$ is a weak equivalence of
spectra for all $X$ if and only if it has the Mayer-Vietoris
property for every abstract blow-up square and every Zariski square
of $pctf$ monoid schemes.
\end{prop}

\begin{proof}
By Theorems \ref{thm:completeregular} and \ref{thm:bounded},
the $cdh$ $cd$-structure is complete, regular and bounded.
Now the assertion follows from \cite[Theorem 3.4]{chsw}.
\end{proof}

Given Proposition \ref{cdh-MV}, the definition of
$cdh$-descent given in \cite[3.5]{chsw} becomes:

\begin{defn}\label{def:cdh-descent}
Let $\cF$ be a presheaf of spectra on $\cMpctf$. We say that $\cF$
satisfies {\em $cdh$ descent} if the canonical map
$\cF(X)\to\bH_{\cdh}(X,\cF)$ is a weak equivalence of spectra for
all $X$.
\end{defn}

\begin{subrem}
Writing $\bH_{\zar}(-, \cF)$ for the fibrant replacement with respect to
the model structure for the Zariski topology, we obtain the notion
of Zariski descent. The proof of Proposition \ref{cdh-MV} applies to show
that $\cF$ satisfies Zariski descent if and only if it has the
Mayer-Vietoris property for every Zariski square. It follows that
$cdh$-descent implies Zariski descent.
\end{subrem}
\smallskip

It is useful to restrict to the full subcategory $\cS$ of smooth
monoid schemes (see Definition \ref{p-smooth}). By Proposition
\ref{smooth-realizn}, these are the cancellative, torsionfree,
separated monoid schemes of finite type whose $k$-realizations are
smooth for any commutative ring $k$. (This condition is independent
of $k$, by \ref{smooth-realizn}.)

\begin{defn}\label{def:scdh}
We define the {\em smooth blow-up $cd$ structure} on $\cS$ to
consist of squares \eqref{sq:cd} such that $X$ is smooth, $e$ is the
inclusion of an equivariant, smooth closed subscheme and $Y$ is the
blow-up of $X$ along $C$. (These assumptions ensure, by
\eqref{same-blowup}, that $Y$ and $D$ are also smooth.)

The Zariski $cd$ structure is given by all cartesian squares in $\cS$
associated to a covering of $X$ by two open subschemes.

We define the $scdh$ topology on $\cS$ to be the Grothendieck
topology associated to the union of the smooth blow-up $cd$-structure
and the Zariski $cd$-structure on $\cS$. For a presheaf of spectra on
$\cS$, we define $\bH_{\scdh}(-, \cF)$ just as $\bH_{\cdh}$ was defined
above. We say such a presheaf $\cF$
satisfies $scdh$ descent if the canonical fibrant replacement map
$$
\cF(X) \to \bH_{\scdh}(X,\cF)
$$
is a weak equivalence for all $X \in \cS$.
\end{defn}

\begin{prop} \label{smooth-cdstructure}
The smooth blow-up $cd$-structure and the Zariski $cd$ structure on
$\cS$ are regular, bounded, and complete. Consequently,
a presheaf of spectra defined on $\cS$ satisfies $scdh$ descent
if and only if it has the Mayer-Vietoris property for
each smooth blow-up square and each Zariski square in $\cS$.
\end{prop}
\goodbreak

\begin{proof}
That the smooth blow-up $cd$-structure is complete can be proved
exactly as Voevodsky did for smooth $k$-schemes in \cite[Lemma 4.3]{VVcomp},
replacing resolution of singularities by our Theorem \ref{thm:resolution2}.
Regularity is proved exactly as in Theorem \ref{thm:completeregular}
for the non-smooth case. The proof that the smooth blow-up $cd$-structure
is bounded works exactly as in Theorem \ref{thm:bounded},
keeping in mind that open subschemes of smooth monoid schemes are smooth.
The proof that the Zariski $cd$-structure is complete, regular and bounded
is again the same as in the non-smooth category.
It follows that the $scdh$ topology is generated by a complete, regular,
bounded $cd$-structure and so \cite[Theorem 3.4]{chsw}
applies to prove the second assertion.
\end{proof}

\begin{prop} \label{P:cdh-scdh}
For any $X \in \cS$ and any presheaf of spectra $\cF$ defined on
$\cMpctf$, we have a weak equivalence
$$
\bH_{{\cdh}}(X, \cF)\ \map{\sim}\
\bH_{\scdh}(X, \cF|_{\cS}) .
$$
\end{prop}

\begin{proof}
In this proof we write $\cF_\cdh$ for the restriction of the
presheaf $\bH_{{\cdh}}(-, \cF)$ to $\cS$. By Proposition \ref{cdh-MV},
$\cF_\cdh$ satisfies the Mayer-Vietoris property for smooth blow-up and
Zariski squares. Therefore $\cF_\cdh$ satisfies $scdh$ descent
(Definition \ref{def:scdh}).

By Theorems \ref{thm:resolution1} and \ref{thm:resolution2}, every
covering sieve for the $cdh$ topology on $\cMpctf$ has a refinement
containing a sieve generated by a cover consisting of objects of
$\cS$. It follows that $\cF|_{\cS} \to \cF_\cdh$ is an $scdh$-local
weak equivalence.
Therefore $\bH_{\scdh}(-,\cF|_{\cS}) \to \bH_{\scdh}(-,\cF_\cdh)$ is an 
objectwise weak equivalence (see \cite[page ~561]{chsw}). Together, the 
two objectwise weak equivalences exhibited in the proof give the assertion.
\end{proof}

\section{Weak $cdh_k$ descent}

Throughout this section, we fix a commutative ring $k$.

\begin{defn}  Let $X_k$ be a scheme of finite type over $k$ and assume
$Z_k\subset X_k$ is a closed subscheme. We say $Z_k$ is {\em regularly
embedded} in $X_k$ if the sheaf of ideals defining $Z_k$ is locally
generated by a regular sequence --- that is, if for all $x \in Z_k$,
the kernel $I_x$ of $\cO_{X_k,x} \to \cO_{Z_k,x}$ is generated by a
$\cO_{X_k,x}$-regular sequence of elements.
\end{defn}

\begin{defn} \label{def:quasi}
A presheaf of spectra $\cF$ defined on $\cMpctf$
has {\em weak $cdh_k$ descent}
if $\cF$ has the Mayer-Vietoris property for each cartesian square
$$\xymatrix{
D \ar[r] \ar[d] & Y \ar[d]^p \\
C \ar[r]^e & X
}$$
in $\cMpctf$ satisfying one of  the following conditions:
\begin{enumerate}

\item It is member of the Zariski cd structure.

\item It is a finite abstract blow-up --- i.e., it is a member of
the abstract blow up cd structure having the
additional property that $p$ is a finite morphism.

\item $C$ is an equivariant closed subscheme,
$Y\to X$ is the blow-up of $X$ along $C$, and
$C_k$ is a regularly embedded closed subscheme of $X_k$.
\end{enumerate}
\end{defn}

\begin{subrem} Theorems \ref{T:descent} and \ref{MainThm} below suggest (but do not prove)
  that the definition of weak $cdh_k$ descent
  is actually independent of the choice of $k$.
\end{subrem}

Since a smooth blow-up square is an example of a blow-up along a
regularly embedded subscheme,
Propositions \ref{smooth-cdstructure} and \ref{P:cdh-scdh}
imply the following theorem.

\begin{thm} \label{T:descent}
If $\cF$ is a presheaf of spectra on $\cMpctf$ that satisfies
weak ${cdh_k}$ descent, then $\cF$ satisfies $scdh$ descent.
That is, the canonical map
$$
\cF(X) \to \bH_{{\cdh}}(X, \cF)
$$
is a weak equivalence for every smooth monoid scheme $X$.
\end{thm}

The main goal of this paper, realized in  the next section, is to
establish a partial generalization of
Theorem \ref{T:descent} to all schemes in $\cMpctf$.
The goal of the rest of this section is to establish some technical
properties needed in the next. We first introduce a
slightly stronger notion than that of weak $cdh_k$ descent.

Recall from \cite[6.10.1]{EGAIV}
that given a closed subscheme $C_k$ of a $k$-scheme $X_k$,
defined by an ideal sheaf $\cI$, $X_k$ is said to be {\em normally flat}
along $C_k$ if the restriction of each $\cI^n/\cI^{n+1}$ to $C_k$ is flat.

\begin{subrem}
Here is a monoid-theoretic condition on a sheaf $I$ of ideals on a
monoid scheme $(X,A)$ which guarantees that, for all $k$,
the $k$-realization of $X$ is normally flat along the $k$-realization
of the equivariant closed submonoid $C$ defined by $I$:
at each point $x$ of $C$, under the natural action of
the monoid $A_x/I_x$ on each of the pointed sets $L_n=I_x^n/I_x^{n+1}$,
each $L_n$ is a bouquet of copies of $A_x/I_x$.
We do not know if this condition is necessary.
\end{subrem}

We will say that a cartesian square of
schemes in $\cMpctf$,
\begin{equation*}
\xymatrix{
D \ar[r] \ar[d] & Y \ar[d]^p \\
C \ar[r]^e & X,
}
\end{equation*}
is a
{\em nice blow-up square} if $C$ is an equivariant closed subscheme of
$X$, $Y$ is the blow-up of $X$ along $C$
and there exists a cartesian square in $\cMpctf$ of the form
\begin{equation} \label{sq:nflat}
\xymatrix{
C \ar[r]_{e} \ar[d] & X \ar[d] \\
B \ar[r] & Z \\
}\end{equation}
such that $Z$ is cancellative, $X \to Z$ is the normalization of $Z$ and
$B$ is an equivariant closed smooth subscheme of $Z$ such that
$Z_k$ is normally flat along $B_k$.

\begin{defn}\label{def:nicecdh}
A presheaf of spectra on $\cMpctf$ satisfies {\em weak+nice $cdh_k$
descent} provided it satisfies weak $cdh_k$ descent
and, in addition, it has the Mayer-Vietoris property for all nice
blow-up squares in $\cMpctf$.
\end{defn}

\begin{prop} \label{P:nice+weak}
If $\cF$ is a presheaf of spectra on $\cMpctf$ that satisfies $cdh$
descent, then $\cF$ satisfies weak+nice $cdh_k$ descent for any
commutativ ring $k$.
\end{prop}

\begin{proof}
This is immediate from Proposition \ref{cdh-MV},
since each of the squares appearing in the definition of weak+nice $cdh_k$
descent is a member of the $cdh$ cd structure.
\end{proof}

We will need  the following technical result about local domains.
Recall that if $I$ is an ideal in a commutative ring $R$ then an ideal
$J\subseteq I$ is called a {\em reduction} of $I$ if $JI^{n-1}=I^n$ for
some $n>0$; a {\em minimal reduction} of $I$ is a reduction which
contains no other reduction of $I$.

\begin{lem} \label{Lem0105}
Let $R$ be a noetherian local domain
with infinite residue field $k$, let
$\fp$ be a prime ideal, and assume $R$ is normally flat along $R/\fp$.
Let  $J$ be a minimal reduction  of $\fp$
that is generated by $h := ht(\fp) = ht(J)$ elements.
(Given $R$ and $\fp$ with these properties, such a $J$ exists by
\cite[5.2, 5.3]{HKH}.)
Let $\tilde{R}$ be the normalization of $R$ and assume $\tilde{R}$ is
Cohen-Macaulay.

Then $J\tilde{R}$ is a reduction of $\fp\tilde{R}$ generated by $h$ elements
and $\Spec(\tilde{R}/J\tilde{R})$ is regularly embedded in
$\Spec(\tilde{R})$.
\end{lem}

\begin{proof}
We have that $J\fp^{n-1}\tilde R=\fp^n\tilde R$, and so the first
assertion is clear.

Since $R \into \tilde{R}$ is an integral extension of
domains, we have $h = ht(J) = ht(J\tilde{R})$.  For any
maximal ideal  $\tilde{\fm}$ of $\tilde{R}$, we have that  $J\tilde{R}_{\tilde{\fm}}$ is
a height $h$ ideal generated by $h$ elements in the local
ring $\tilde{R}_{\tilde{\fm}}$. Since $\tilde{R}_{\tilde{\fm}}$ is
Cohen-Macaulay by assumption, these generators necessarily form a
regular sequence.
\end{proof}

The following is the evident analogue of the notion of weak $cdh_k$
descent for presheaves of spectra on the category of $k$-schemes.

\begin{defn} \label{def:weak/k}
For a commutative ring $k$, let $Sch/k$ be the category of separated
schemes essentially of finite type over $k$. A presheaf of spectra
defined on $Sch/k$ satisfies {\em weak $cdh$ descent} if it has the
Mayer-Vietoris property for each cartesian square
$$\xymatrix{
D \ar[r] \ar[d] & Y \ar[d]^p \\
C \ar[r]^e & X
}$$
of schemes satisfying one of the following conditions:
\begin{enumerate}
\item $e$ and $p$ are open immersions whose images cover $X$.

\item It is a finite abstract blow-up --- i.e., $e$ is a closed
  immersion, $p$ is finite, and $p$ maps $Y \setminus D$
  isomorphically onto $X \setminus C$

\item $e$ is a regular closed immersion and $p$ is the blow-up
  of $X$ along $C$.
\end{enumerate}
\end{defn}

\begin{lem}\label{keylemma}
Assume $k$ is a commutative regular noetherian 
domain containing an infinite field and
$\cG_k$ is a presheaf of spectra on $Sch/k$ that satisfies weak
$cdh$ descent. Let $\cG$ be the presheaf of spectra on $\cMpctf$
defined by $\cG(X) := \cG_k(X_k)$.

Then $\cG$ satisfies weak+nice $cdh_k$ descent on $\cMpctf$.
\end{lem}

\begin{proof}
Since the $k$-realizations of the squares involved in
the definition  of weak $cdh_k$ descent for $\cMpctf$
(Definition \ref{def:quasi}) are squares
involved in the definition of weak $cdh$ descent
for $Sch/k$ (Definition \ref{def:weak/k}),
it follows that $\cG$ satisfies weak $cdh_k$ descent.
Say $X, Y, C, D, Z$, and $B$ are as in the definition
of a nice blow-up square.  We need to prove that the square
\begin{equation} \label{E57}
\xymatrix{
\cG_k(X_k) \ar[r] \ar[d] & \cG_k(C_k) \ar[d] \\
\cG_k(Y_k) \ar[r] & \cG_k(D_k)
}
\end{equation}
is homotopy cartesian.

Let $R$ be any local ring of $Z_k$ and let $\fp$ be the
prime ideal of $R$ cutting out $B_k$ locally.
Let $V = \Spec(\tilde{R}_{\tilde{\fm}})$ where $\tilde{R}$ is the
normalization of $R$ and $\tilde{\fm}$ is any of the maximal ideals of
$\tilde{R}$. Then, since $X_k$ is the normalization of $Z_k$ by
Proposition \ref{P:normal}, $V$ is the spectrum of a local ring of
$X_k$, and for various choices of $R$ and $\tilde{\fm}$, every local
ring of $X_k$ arises in this manner.

By Corollary \ref{cor:pullbackcommute}, $C_k=X_k\times_{Z_k}B_k$, so
the closed subscheme $V\times_{X_k}C_k$ of $V$ is cut out by
$\fq=\fp\tilde{R}_{\tilde{\fm}}$. As $X$ is the normalization of the
separated cancellative, torsionfree  monoid scheme $Z$, Proposition
\ref{P:normal} implies that $X_k$ is a toric variety. By
\cite{Hochster}, all toric schemes over $k$ are Cohen-Macaulay;
hence so are $X_k$ and $V$.

By Lemma  \ref{Lem0105},
$\fq = \fp\tilde{R}_{\tilde{\fm}}$ admits a reduction $I \subset \fq$
such that $\Spec(\tilde{R}_{\tilde{\fm}}/I) \into V$ is a regular embedding.
Since $V\times_{X_k}Y_k$ is the blow-up of $V_k$ along $V\times_{X_k}C_k$
(by Proposition \ref{same-blowup}), and the exceptional divisor is
$V\times_{X_k}D_k$ (by \ref{cor:pullbackcommute}),
the proof of \cite[5.6]{HKH} (with $\cKH$ replaced by $\cG$) gives that
\[ \xymatrix{
\cG_k(V) \ar[r] \ar[d] & \cG_k(V \times_{X_k} C_k) \ar[d] \\
\cG_k(V \times_{X_k} Y_k) \ar[r] & \cG_k(V \times_{X_k} D_k)
}\]
is homotopy cartesian. Since $\cG_k$
satisfies the Mayer-Vietoris property for Zariski covers
and the $V$ occurring here is an arbitrary local scheme of $X_k$,
the proof of \cite[5.7]{HKH} (with $\cKH$ replaced by $\cG_k$)
shows that \eqref{E57} is homotopy cartesian.
\end{proof}

\begin{ex}  \label{Ex0106}
Let $\cKH$ denote Weibel's homotopy algebraic $K$-theory \cite{Wei}.
We may view $\cKH$ as a presheaf of spectra on $Sch/k$. By abuse of notation,
we also write $\cKH$ for the
presheaf of spectra on $\cMpctf$ defined by $\cKH(X) = \cKH(X_k)$.

By \cite{TT}, \cite[4.9]{Wei} and \cite{Th},
$\cKH$ satisfies weak $cdh$ descent on $Sch/k$ (\ref{def:weak/k});
by Lemma \ref{keylemma},
$\cKH$ satisfies weak+nice $cdh_k$-descent on $\cMpctf$.
\end{ex}

\section{Main Theorem}\label{sec:main}

In this section, we prove our main theorem (Theorem \ref{MainThm}),
which gives a condition for $\cF$ to satisfy $cdh$ descent on $\cMpctf$.
We will need the Bierstone-Milman Theorem, which we extract from the
embedded version \cite[Thm.\,1.1]{BM}.

\begin{thm}\label{thm:BM}
Let $X$ be a separated cancellative torsionfree monoid scheme
of finite type, embedded as a closed subscheme (see  Definition
\ref{def:closedimm}) in a smooth toric monoid scheme $M$
(see Definition \ref{def:toric}). For any commutative ring $k$ containing a field, there  is a sequence of blow-ups
along smooth equivariant centers $Z_i \subset X_i$, $0 \leq i \leq n-1$,
$$
Y = X_n \to \cdots \to X_0 = X
$$
such that $Y$ is smooth, and each
$(X_i)_k$ is normally flat along $(Z_i)_k$.
\end{thm}

\begin{proof}
Since normal flatness is stable under flat extension of the base, 
and $k$ is flat over a field, we may assume that $k$ is a field.
Let $\bar k$ denote the algebraic closure of $k$, and let
$T$ be the torus acting on $M_{\bar k}$.
The Bierstone-Milman Theorem (\cite[Thm.\,1.1]{BM}) tells us that
we can find a sequence of blow-ups $M_n\to\cdots\to M_0=M_{\bar k}$ of
smooth toric $\bar k$-varieties, the blow-up of $M_i$ being taken along a
smooth $T$-invariant center $N_i$, with the following properties.
Setting $X'_0=X_{\bar k}$, we inductively define $Z'_i=N_i\cap X'_i$;
then $Z'_i$ is a smooth equivariant $\overline{k}$-variety, $X'_i$ is normally flat
along  $Z'_i$, and $X'_{i+1}$ is the strict transform of $X'_i$.

The $\bar k$-realization functor from fans to (normal) toric
$\bar k$-varieties (and equivariant morphisms) is well known to be
an equivalence.
It follows that each of the $N_i$ and $M_i$ and the morphisms between
them come from fans, and hence by Theorem \ref{why_toric} are
$\bar k$-realizations of toric monoid schemes (which by abuse of notation,
we will call $N_i$ and $M_i$), and morphisms of such.

Inductively we define monoid schemes $X_i$ and $Z_i$, starting from
$X_0=X$ and $Z_0=N_0\cap X$, to be the blow-up of the monoid scheme
$X_{i-1}$ along $Z_{i-1}$ in the sense of \ref{def:blowup}.
By Proposition \ref{same-blowup}
and Corollary \ref{cor:pullbackcommute},
$Z'_i=(Z_i)_{\bar k}$ and
$X'_i=(X_i)_{\bar k}$. In particular, $(X_n)_{\bar k}=Y$ is a
smooth toric variety and therefore the monoid scheme $X_n$ is smooth
by Proposition \ref{smooth-realizn}.
Finally, faithfully flat descent implies that
$(X_i)_k$ is normally flat along $(Z_i)_k$ if and only if
$(X_i)_{\bar k}$ is normally flat along $(Z_i)_{\bar k}$.
\end{proof}

\begin{thm} \label{Prelude}
Suppose $\cG$ is a presheaf of spectra on $\cMpctf$ satisfying
weak+nice $cdh_k$ descent for some commutative ring $k$ containing a field.  
If  $\cG(X) \simeq\ast$
for all $X$ in $\cS$ then $\cG(X) \simeq\ast$ for all $X$ in $\cMpctf$.
\end{thm}

\begin{proof}
We proceed by induction on the dimension of $X$.
Given $X$, let $x_1,\dots,x_l$ be its generic points, and let
$Y_i=\overline{\{x_i\}}{}^\eq$ be their equivariant closures
(see Lemma \ref{equivclosure}).
We have a cover $X = Y_1 \cup \cdots \cup Y_l$ by
equivariant closed subschemes each of which is cancellative by
Example \ref{ex:birat-components}.
Moreover, each $Y_i\times_X Y_j$ is equivariant and closed, hence pctf.
Since $\cG$ has the Mayer-Vietoris property for closed covers, and $\cG$
vanishes on the $Y_i\times_X Y_j$ for all $i \ne j$ by the induction
hypothesis, we get
$$
\cG(X) = \prod_i \cG(Y_i).
$$
We may thus assume that $X$ is cancellative. (This also establishes
the base case $\dim(X)=0$, since in that case the $Y_i$ are in $\cS$.)

Since $\cG$ satisfies Mayer-Vietoris for open covers, we may assume
$X$ is affine. In  particular, we may assume $X$ can be embedded in a
smooth toric monoid scheme, for example, by choosing a surjection from a
free abelian monoid onto $A$ where $X = \MSpec(A)$. This will allow
us to apply the Bierstone-Milman Theorem \ref{thm:BM} to obtain
a sequence of blow-ups along smooth monoid schemes $Z_i$,
$$
Y = X_n \to \cdots \to X_0 = X.
$$

We claim that $\cG(X_i) \simeq \cG(X_{i+1})$ for all $i$.
Since $\cG(Y)\simeq\ast$, this will finish the inductive step
and hence the proof of the theorem.
To simplify the notation, fix $i$ and write 
$Z$ for $Z_i \subset X_i$ and $X_Z$ for $X_{i+1}$,
the blow-up of $X_i$ along $Z$, so that our goal is to prove that
$\cG(X_i)\to\cG(X_Z)$ is a weak equivalence.
Let $\tilde{X}$ denote the normalization $(X_i)_\nor$ of $X_i$ and set
$\tilde{Z} = Z \times_{X_i} \tilde{X}$. Write $\tilde{X}_{\tilde{Z}}$
for the blow-up of $\tilde{X}$ along $\tilde{Z}$.
By naturality of blow-ups (see \ref{def:blowup}),
there is a commutative square
$$
\xymatrix{
\tilde{X}_{\tilde{Z}}\ar[r] \ar[d] & X_Z \ar[d] \\
\tilde{X} \ar[r] & X_i \\
}
$$
(that need not be cartesian). Since the map $\tilde{X}\to X_i$ is finite,
the map $\tilde{X}_{\tilde{Z}} \to X_Z$ is also finite, by Lemma \ref{L0106}.
Applying $\cG$ gives a commutative square
of spectra
$$
\xymatrix{
\cG(\tilde{X}_{\tilde{Z}}) & \cG(X_Z) \ar[l] \\
\cG(\tilde{X}) \ar[u] & \cG(X_i).  \ar[l] \ar[u]\\
}
$$
To prove that the right-hand vertical arrow is a weak equivalence,
it suffices to prove the other three are.

The finite map $\tilde{X} \to X_i$ is an isomorphism on
the generic points.
Consider the equivariant closure $E \subset X_i$ of the
finitely many height~1 points of $X_i$; by  Remark \ref{rem:equiclose},
every point in $E$ has height $\ge1$ in $X_i$, so
$E$ is the complement of the generic point of $X_i$.
Since $E$ is pctf, $\cG(E)\simeq\ast$ by our inductive assumption.
Since the pullback $\tilde{E} := E\times_{X_i}\tilde{X}$ is an
equivariant closed subscheme of $\tilde{X}$, it is pctf by
Proposition \ref{pctfproperties}, and hence $\cG(\tilde{E})\simeq\ast$
as well, by induction. Using the finite abstract blow-up square
involving $X_i$, $\tilde{X}$, $E$ and $\tilde{E}$, we have a weak equivalence
$$
\cG(X_i) \map{\simeq} \cG(\tilde{X}).
$$
The map $\tilde{X}_{\tilde{Z}} \to X_Z$ is also finite and birational,
and so the same argument shows
$$
\cG(X_Z) \map{\simeq} \cG(\tilde{X}_{\tilde{Z}}).
$$
is a weak equivalence.
Finally, observe that
$$\xymatrix{
\tilde{Z} \times_{\tilde{X}} \tilde{X}_{\tilde{Z}} \ar[r] \ar[d] &
\tilde{X}_{\tilde{Z}} \ar[d] \\
\tilde{Z} \ar[r] &
\tilde{X}
}$$
is a nice blow-up square, because the bottom row may be compared with
$Z\to X_i$ and $(X_i)_k$ is normally flat along $Z_k$. Because
$\cG$ has descent for nice blow-up squares, and
$\cG(\tilde{Z})\simeq
\cG(\tilde{Z}\times_{\tilde{X}}\tilde{X}_{\tilde{Z}})
\simeq\ast$
by the induction hypothesis, we get a weak equivalence
$$
\cG(\tilde{X}) \map{\simeq} \cG(\tilde{X}_{\tilde{Z}}).
$$
It follows that $\cG(X_i)\simeq\cG(X_Z)$, as claimed.
This completes the proof.
\end{proof}

We now state and prove the main theorem of this paper, which
gives a partial generalization of  Theorem \ref{T:descent} to all objects in
the category $\cMpctf$.

\begin{thm} \label{MainThm}
Let $\cF_k$ be a presheaf of spectra on $Sch/k$ for some commutative regular
noetherian ring $k$ containing an infinite field, and define $\cF$ to be the
presheaf of spectra on $\cMpctf$ defined by $\cF(X) = \cF_k(X_k)$.

If $\cF_k$ satisfies weak $cdh$ descent on $Sch/k$, then
$\cF$ satisfies $cdh$ descent on $\cMpctf$.
\end{thm}

\begin{proof}
Let $\cG$ be the homotopy fiber of $\cF\to\bH_{\cdh}(-,\cF)$ --- i.e.,
for all $X$ in $\cMpctf$, $\cG(X)$ is the homotopy fiber of
$\cF(X)\to\bH_{{\cdh}}(X,\cF)$.
By Lemma \ref{keylemma} and Proposition \ref{P:nice+weak},
both $\cF$ and $\bH_{\cdh}(-, \cF)$ satisfy
weak+nice $cdh_k$ descent, and hence $\cG$ satisfies weak+nice
$cdh_k$ descent too. Theorem \ref{T:descent} gives that $\cG(X) \simeq\ast$
for all $X \in \cS$. Now we apply Theorem \ref{Prelude} to
conclude $\cG(X) \simeq\ast$ for all $X$ in $\cMpctf$.
\end{proof}

The following corollary is the Theorem announced in the introduction.

\begin{cor} \label{MainCor}
Assume $k$ is a commutative regular noetherian ring 
containing an infinite field and let
$\cF_k$ be a presheaf of spectra on $Sch/k$ that satisfies the
Mayer-Vietoris property for Zariski covers, finite abstract blow-up
squares, and blow-ups along regularly embedded subschemes. 

Then $\cF_k$ satisfies the Mayer-Vietoris property for all abstract
blow-up squares of toric $k$-schemes obtained from subdividing a fan.
\end{cor}

\begin{proof}
By Definition \ref{def:weak/k}, $\cF_k$ satisfies weak $cdh$ descent
on $Sch/k$.
      By Theorem  \ref{MainThm}, $\cF$ satisfies $cdh$ descent in $\cMpctf$.
      Now use Proposition \ref{cdh-MV}.
\end{proof}

\begin{cor} \label{Cor0106}
Let $k$ be a commutative regular noetherian ring containing a field. 
The presheaf of spectra $\cKH$ on $\cMpctf$, 
defined as $\cKH(X) = \cKH(X_k)$, satisfies $cdh$ descent. 
Moreover, both natural maps
$$
\cKH(X) \to \bH_{\cdh}(X, \cKH) \leftarrow \bH_{\cdh}(X, \cK)
$$
are weak equivalences for all $X$ in $\cMpctf$.
\end{cor}

\begin{proof}
We first reduce to the case when $k$ is of finite type over a field.
We can express $k$ as a filtered colimit of rings $k_i$, all regular of
finite type over a field (by Popescu's theorem \cite[2.5]{Popescu}). 
The functor $\cKH$ is the
homotopy colimit of the corresponding functors defined by
$k_i$-realization. By Proposition \ref{cdh-MV}, 
we can check descent by showing that certain
squares of monoid schemes are transformed by $\cKH$ into homotopy
co-cartesian squares of spectra (a square of spectra is
homotopy cartesian if and only if it is homotopy co-cartesian);
since homotopy colimits of homotopy co-cartesian squares are 
homotopy co-cartesian, we may assume that 
$k$ is of finite type over its field of constants.  

Now if the regular ring $k$ does not contain an infinite field,
it is smooth over the (perfect) field of constants it contains and 
hence stays regular under base change from its field of constants 
to any algebraic extension. We
can therefore apply the standard transfer argument and may assume that
$k$ contains an infinite field. 

By Example \ref{Ex0106} and
Theorem \ref{MainThm}, $\cKH$ satisfies $cdh$ descent on $\cMpctf$.
For any $X$ in $\cMpctf$, consider the commutative square of spectra:
$$
\xymatrix{
\cK(X) \ar[r] \ar[d] & \cKH(X) \ar[d] \\
\bH_{\cdh}(X, \cK) \ar[r] & \bH_{\cdh}(X, \cKH), \\
}$$
where $\cK$ is algebraic $K$-theory, regarded as a presheaf of spectra
on $Sch/k$ and hence on $\cMpctf$. Since $\cKH$ satisfies $cdh$ descent,
the right-hand vertical map is a weak equivalence for all $X$.
This is the first assertion of the corollary.

If $X$ is smooth, then the top horizontal map is a weak equivalence
by \cite{Wei} (since $X_k$ is smooth over $k$ hence regular by \ref{smooth-realizn}).
By fibrant replacement and Proposition \ref{P:cdh-scdh},
the bottom map is also a weak equivalence for all $X$ in $\cS$.
By induction on $\dim(X)$ and Theorem \ref{thm:resolution1},
this implies that $\bH_{\cdh}(-,\cK) \to \bH_{\cdh}(-,\cKH)$
is a local weak equivalence and, as observed (for any site)
in \cite[p.\,561]{chsw}, this implies that
$\bH_{\cdh}(X,\cK) \to \bH_{\cdh}(X, \cKH)$
is a weak equivalence for all $X$ in $\cMpctf$.
\end{proof}

\begin{rem}\label{rem:kdim}
It follows from Corollaries \ref{cor:bounded} and \ref{Cor0106} and
a $cdh$-descent argument that if $X\in\cMpctf$ is of dimension $d$
and $k$ is a commutative regular ring containing a field, then $KH_n(X_k)=0$ for
$n<-d$ (cf. \cite[Thm. 8.19]{HKH}). The analogous statement for
$K$-theory is also true, at least if $X$ is cancellative and
torsion-free. Indeed for affine $X$, $K_n(X_k)=0$ for $n<0$, by
\cite[Thm. 1.3]{gubelaff}; the general case follows from this by a
Zariski descent argument, using \ref{cor:bounded}.
\end{rem}

\medskip
In order to apply Corollary \ref{Cor0106} to the relation between
$K$-theory and topological cyclic homology, we need to recall some terms.
Fix a prime $p$ and a commutative regular ring $k$ of characteristic $p$. To each
scheme $X$ essentially of finite type over $k$, there is a pro-spectrum
$\{ TC^\nu(X,p)\}_{\nu=0}^\infty$ and the cyclotomic trace is a compatible
family of morphisms $tr^\nu:\cK(X)\to TC^\nu(X,p)$.
Define $\cF_k^\nu$ to be the presheaf of spectra on $Sch/k$ given as the
homotopy fiber of $\cK(X)\to TC^\nu(X,p)$. Then Geisser and Hesselholt
observe in the proof of \cite[Thm.\,B]{GHvanish} that each $\cF_k^\nu$
takes elementary Nisnevich squares and regular blow-up squares to
homotopy cartesian squares of pro-spectra. 

Following Geisser-Hesselholt \cite{GHvanish}, a strict map
of pro-spectra $\{ X^\nu\}\to \{ Y^\nu\}$ is said to be a
{\it weak equivalence} if for every $q$ the induced map
$\{ \pi_q(X^\nu)\} \to \{\pi_q(Y^\nu)\}$
is an isomorphism of pro-abelian groups. A square diagram of strict maps
of pro-spectra is said to be {\it homotopy cartesian} if the canonical map
from the upper left pro-spectrum
to the level-wise homotopy limit of the other terms is a weak equivalence.

Given a class $\cC$ of squares we will say that a pro-presheaf of spectra
satisfies the pro-analogue of $\cC$-descent if it sends each square in $\cC$
to a homotopy cartesian square of pro-spectra.

Define $\{\cF^\nu\}$ to be the pro-presheaf of spectra on $\cMpctf$
given as the family of homotopy fibers of the maps $\cK(-) \to TC^\nu(-,p)$.
That is, $\cF^\nu(X)=\cF^\nu_k(X_k)$
is the homotopy fiber of $\cK(X_k)\to TC^\nu(X_k,p)$ for each $X$ and $\nu$.

\begin{prop}\label{K-TC}
Assume $k$ is a commutative regular noetherian ring containing an infinite field of
characteristic $p
> 0$. Then $\{\cF^\nu\}$ satisfies $cdh$ descent on $\cMpctf$ in the
sense that $\{\cF^\nu\} \to\{\bH(-,\cF^\nu)\}$ is a weak equivalence
of pro-spectra.
\end{prop}

\begin{proof}
Fix $\nu$ and let $\cG^\nu$ be the homotopy fiber of
$\cF^\nu\to\bH_{\cdh}(-,\cF^\nu)$. It suffices to prove that for each $X$
and $q$ the pro-abelian group $\{\pi_q\cG^\nu(X)\}$ is pro-zero.
We will do so by modifying the proof of Theorem \ref{MainThm}.

For each $\nu$, $\bH_{\cdh}(-,\cF^\nu)$ satisfies weak+nice $cdh_k$ descent
by Proposition \ref{P:nice+weak}.
By \cite[Thm.\,1]{GHbirel} and \cite[Thms.\,B, D]{GHbound},
$\{\cF_k^\nu\}$ sends finite abstract blow-up squares to
homotopy cartesian squares of pro-spectra. Thus $\{\cF_k^\nu\}$
satisfies the pro-analogue of weak $cdh$ descent (Definition \ref{def:weak/k}).
In the proof of Lemma \ref{keylemma}, the reduction ideals used are
reduction ideals on affine neighborhoods of the maximal ideal $\fm$ of $R$.
By the argument used in the proof of \cite[Thm.\,1.1]{GHvanish}, the proof of
our Lemma \ref{keylemma} now applies {\it mutatis mutandis} to show that
the pro-presheaf of spectra $\cF^\nu$ satisfies the Mayer-Vietoris
property for nice blow-up squares.
It now follows that $\{\cG^\nu\}$ satisfies
the pro-analogue of weak+nice $cdh_k$ descent.

For each $\nu$, $\cF^\nu$ satisfies Zariski descent and also has the
Mayer-Vietoris property for regular blow-ups, so $\cF^\nu$ satisfies
$scdh$ descent by \ref{smooth-cdstructure}.
By definition, this means that for each smooth $X$ the spectrum
$\cG^\nu(X)$ is contractible.
Now the proof of Theorem \ref{Prelude} applies {\it verbatim}
to finish the proof.
\end{proof}

\begin{cor}\label{cor:pro-hcart}
Assume $k$ is any commutative regular noetherian ring of characteristic $p>0$. For any
monoid scheme $X$ in $\cMpctf$, the following square of pro-spectra
is homotopy cartesian.
$$\xymatrix{
\cK(X) \ar[r] \ar[d] &  \cKH(X)  \ar[d] \\
\{TC^\nu(X,p)\} \ar[r] &  \{\bH_{{\cdh}}(X, TC^\nu(-,p))\}.
}$$
\end{cor}

\begin{proof}
By a standard transfer argument as in Corollary \ref{Cor0106}, 
we may assume that $k$ contains an
infinite field. By Proposition \ref{K-TC}, the homotopy fiber
$\{\cF^\nu(X)\}$ of the left vertical map is weakly equivalent to
$\{\bH_{\cdh}(X,\cF^\nu)\}$.
By Corollary \ref{Cor0106}, this coincides up to weak equivalence
with the homotopy fiber of the right vertical map.
\end{proof}

\begin{rem}\label{rem:notgh}
As explained in Remark \ref{rem:kdim}, if $k$ is any commutative regular ring
containing a field, and $X\in\cMpctf$ is cancellative and
torsion-free then $K_n(X_k)=0$ for $n<-\dim X$. To extend this
result to all $X\in\cMpctf$ it would suffice to prove that the
bottom horizontal map in the diagram in Corollary
\ref{cor:pro-hcart} induces an isomorphism (resp. an epimorphism) of
homotopy groups in degrees $n < -\dim(X)$ (resp. $n=-\dim(X)$). Geisser
and Hesselholt proved the analogue statement for schemes essentially of finite
type over a field of positive characteristic which admits resolution of
singularities (\cite[Thm. C]{GHvanish}). Adapting their methods to our situation 
seems rather hard.
\end{rem}

\subsection*{Acknowledgements}
The authors would like to thank the referee for a careful reading,
for suggesting the notion of a monoid poset and for the current proof
of Lemma \ref{lem:S-and-k}.

\bibliographystyle{plain}

\end{document}